\documentclass[a4paper,fleqn]{cas-sc}

\usepackage[numbers]{natbib}

\def\tsc#1{\csdef{#1}{\textsc{\lowercase{#1}}\xspace}}
\tsc{WGM}
\tsc{QE}

\usepackage{amssymb}
\usepackage{latexsym}
\usepackage{amsmath}

\usepackage{subcaption}
\usepackage{float}

\usepackage{caption}
\usepackage{booktabs}
\usepackage{tabularx}

\usepackage{siunitx}

\usepackage{placeins}

\DeclareMathOperator{\grad}{\nabla}
\DeclareMathOperator{\dive}{\nabla\cdot}
\newcommand{\ave}[1]{\left\{\!\left\{#1\right\}\!\right\}}
\newcommand{\jump}[1]{\left[\!\left[#1\right]\!\right]}
\newcommand{\tjump}[1]{\left<\!\left<#1\right>\!\right>}
\newcommand{\pad}[2]{\frac{\partial{#1}}{\partial{#2}}}

\usepackage{amsthm}
\newtheorem{theorem}{Theorem}[section]
\newtheorem{definition}[theorem]{Definition}

\usepackage[normalem]{ulem}

\definecolor{dartmouthgreen}{rgb}{0.05, 0.5, 0.06}

\begin{document}

\let\WriteBookmarks\relax
\def\floatpagepagefraction{1}
\def\textpagefraction{.001}
	
\shorttitle{AP and AA IMEX methods for Euler equations: non-ideal gases}    
	
\shortauthors{G. Orlando et al.}  
	
\title[mode = title]{A quantitative comparison of high-order asymptotic-preserving and asymptotically-accurate IMEX methods for the Euler equations with non-ideal gases}  
	
\author[1]{Giuseppe Orlando}[orcid=0000-0002-7119-4231]
\cormark[1]
\ead{giuseppe.orlando@polytechnique.edu}
	
\author[2]{Sebastiano Boscarino}[]
\ead{sebastiano.boscarino@unict.it}

\author[2]{Giovanni Russo}[]
\ead{giovanni.russo1@unict.it}
	
\affiliation[1]{organization={CMAP, CNRS, \'{E}cole polytechnique, Institut Polytechnique de Paris},
addressline={Route de Saclay}, 
city={Palaiseau},
postcode={91120}, 
country={France}}	
	
\affiliation[2]{organization={Department of Mathematics and Computer Science, University of Catania},
addressline={}, 
city={Catania},
postcode={95125}, 
country={Italy}}
	
\cortext[1]{Corresponding author}
	
\begin{abstract}
We present a quantitative comparison between two different Implicit-Explicit Runge-Kutta (IMEX-RK) approaches for the Euler equations of gas dynamics, specifically tailored for the low Mach limit. In this regime, a classical IMEX-RK approach involves an implicit coupling between the momentum and energy balance so as to avoid the acoustic CFL restriction, while the density can be treated in a fully explicit fashion. This approach leads to a mildly nonlinear equation for the pressure, which can be solved according to a fixed point procedure. An alternative strategy consists of employing a semi-implicit temporal integrator based on IMEX-RK methods (SI-IMEX-RK). The stiff dependence is carefully analyzed, so as to avoid the solution of a nonlinear equation for the pressure also for equations of state (EOS) of non-ideal gases. The spatial discretization is based on a Discontinuous Galerkin (DG) method, which naturally allows high-order accuracy. The asymptotic-preserving (AP) and the asymptotically-accurate (AA) properties of the two approaches are assessed on a number of classical benchmarks for ideal gases and on their extension to non-ideal gases.
\end{abstract}
	
	
\begin{highlights}
    \item Quantitative comparison between IMEX and SI-IMEX asymptotic-preserving methods
    \item Asymptotically-accurate schemes up to fourth order
    \item Extension of SI-IMEX method for non-ideal gases
    \item Analysis of the stiff dependence for non-ideal gases so as to avoid the solution of a nonlinear equation
\end{highlights}
	
	
\maketitle

\section{Introduction}
\label{sec:intro}

The Euler equations of gas dynamics represent the standard mathematical model for a wide range of applications in fluid mechanics, mechanical engineering, and environmental engineering \cite{patankar:2018}. Several numerical methods have been developed over the years, which can generally be divided into two categories, depending on a dimensionless parameter called Mach number. The Mach number $M$ represents the ratio between the local fluid velocity and the speed of sound in the medium \cite{graebel:2018}. For moderate to high Mach numbers, compressible effects have to be taken into account and numerical discretization strategies typically rely on Godunov-type shock-capturing schemes \cite{godlewski:2013, harten:1983, leveque:2002, munz:1994, shu:2009, toro:2009}. On the other hand, in the low Mach number regime, the flow can be considered weakly compressible or even incompressible. 

Explicit time discretization methods are very popular for high Mach number flows \cite{munz:1994, osher:1982, shu:2009}. The time step has to satisfy a Courant-Friedrichs-Lewy (CFL) condition, given by the mesh size divided by the fastest wave speed \cite{boscarino:2022, courant:1928}. For moderate to high Mach numbers, this restriction is not a problem, since one is interested in resolving all the waves. However, for flows characterized by low values of the Mach number, severe time step restrictions may be required by these schemes. In this regime, acoustic waves usually carry a negligible amount of energy and, therefore, one may not be interested in resolving them. Hence, the system becomes stiff and stability limitations on the time step are much stricter than the restrictions imposed by accuracy.

The use of implicit and semi-implicit time discretization methods, so as to avoid the acoustic CFL restriction, has a long tradition in low Mach number flows \cite{arun:2020, boscarino:2022, casulli:1984, klein:1995, orlando:2022, orlando:2023b, orlando:2025}. Several numerical methods for weakly compressible flows have been proposed in the literature, see, e.g., \cite{chalons:2013, cordier:2012, dimarco:2017, kucera:2022, noelle:2014, thomann:2019} and the references therein. Since the seminal paper \cite{casulli:1984}, an effective approach to deal with low Mach flows is given by pressure-based algorithms. Indeed, an implicit treatment of the pressure gradient term within the momentum equation and of the pressure work term in the energy equation is sufficient to remove the acoustic CFL restriction and to decouple acoustic and transport effects \cite{casulli:1984}. 

Another aspect to consider is that weakly compressible flows are characterized by multiple length and time scales. Under specific assumptions, the compressible Euler equations converge to the incompressible ones as the Mach number goes to zero \cite{feireisl:2016, klainerman:1981}. Hence, the density remains constant along the fluid particle trajectories and the pressure acts as a Lagrange multiplier to enforce incompressibility of the flow \cite{boscarino:2022, klein:1995}. Robust numerical methods for the Euler equations should recover the incompressible limit for vanishing Mach number. For this purpose, the concept of asymptotic-preserving (AP) schemes has been introduced \cite{haack:2012}. We also refer to \cite{boscarino:2024a} for a review of asymptotic-preserving methods for quasilinear hyperbolic systems with stiff relaxation. A numerical method for the compressible Euler equations is said to be asymptotic-preserving if its stability condition does not depend on the Mach number $M$ and if it provides a consistent discretization of the incompressible Euler equations as $M \to 0$. The AP property of the aforementioned approach was proven in \cite{orlando:2025}.

Preserving incompressibility and resolving vortex dynamics are among the main purposes of numerical discretizations for weakly-compressible flows and high-order methods can help to reach this goal. The aim of the present work is to provide a quantitative comparison between two different Implicit-Explicit Runge-Kutta (IMEX-RK) approaches for the Euler equations, specifically tailored for the low Mach limit. The first approach uses an IMEX-RK solver, as proposed in \cite{orlando:2022} and validated in \cite{orlando:2023b} for atmospheric applications. A key feature of this approach is the implicit treatment of acoustic waves, while material waves are handled explicitly. This method combines IMEX-RK schemes, carefully designed for stability and accuracy, with a time-stepping size that is independent of the Mach number $M$. The spatial discretization is based on the Discontinuous Galerkin (DG) method \cite{giraldo:2020}, which naturally allows for high-order accuracy and has proven highly effective for a wide range of computational fluid dynamics problems, across various flow regimes (see, e.g.,  \cite{bassi:1997, cockburn:1989b, cockburn:1989a, karniadakis:2005}). Additionally, the IMEX-DG method can handle a general equation of state (EOS) \cite{orlando:2022, orlando:2025}, for which only a few studies have been devoted \cite{abbate:2019, cordier:2012}, and that, as we will discuss later, poses additional numerical challenges. However, this approach leads to a mildly nonlinear equation for the pressure also for ideal gases. In order to avoid this nonlinearity, an alternative strategy consists of employing a semi-implicit temporal integrator based on IMEX-RK schemes (SI-IMEX-RK), similar to the one adopted in \cite{boscarino:2016}. This method leads to a linearized equation for both the pressure and the EOS. We refer to this scheme as the SI-IMEX-DG.

Furthermore, in this work, we show that the time discretization methods satisfy the AP property and the asymptotically accuracy (AA) property, i.e. they maintain their high-order accuracy also in the case of $M \to 0$. Moreover, if specific boundary conditions are considered, the limit model can differ from the incompressible Euler equations and depends on the employed EOS \cite{orlando:2025}. Specific numerical treatments for non-ideal EOS in the framework of the IMEX-DG scheme were presented in \cite{orlando:2022}. In \cite{boscheri:2021a}, for non-ideal EOS, a nonlinear equation is solved through a Newton method. We propose here a novel different strategy, so as to avoid any nonlinear equation for the semi-implicit approach. 

The paper is structured as follows. In Section \ref{sec:model}, we briefly recall the mathematical model and its limit as $M \to 0$. In Section \ref{sec:num_disc}, we present the numerical method. More specifically, we outline the IMEX and the SI-IMEX time discretization methods. Moreover, we provide suitable strategies to deal with a general class of EOS. Some details of the DG formulation will also be discussed, specifying some advantages and disadvantages of this method for low Mach number flows. In Section \ref{sec:num_res}, some numerical results to assess the properties of the two methods and to compare them are presented. Finally, some conclusions and perspectives for future work are discussed in Section \ref{sec:conclu}.

\section{The mathematical model}
\label{sec:model}

Let $\Omega \subset \mathbb{R}^{d}, 1 \le d \le 3$ be a connected open bounded set with a sufficiently smooth boundary $\partial\Omega$ and denote by $\mathbf{x}$ the spatial coordinates and by $t$ the temporal coordinate. The mathematical model consists of the fully compressible Euler equations of gas dynamics, written in non-dimensional form as follows \cite{boscarino:2022, boscheri:2021a, klein:1995, orlando:2025}:
\begin{eqnarray}\label{eq:euler_energy}
    \frac{\partial\rho}{\partial t} + \dive\left(\rho\mathbf{u}\right) &=& 0 \nonumber \\
    \frac{\partial\rho\mathbf{u}}{\partial t} + \dive\left(\rho\mathbf{u} \otimes \mathbf{u}\right) + \frac{1}{M^{2}}\grad p &=& \mathbf{0} \\
    \frac{\partial\rho E}{\partial t} + \dive\left[\left(\rho E + p\right)\mathbf{u}\right] &=& 0. \nonumber 
\end{eqnarray}
Here, $\rho$ is the density, $\mathbf{u}$ is the fluid velocity, $p$ is the pressure, and $E$ is the total energy per unit of mass. Moreover, $M \equiv u_{0}/\sqrt{p_{0}/\rho_{0}}$, where $u_0$, $p_0$ and $\rho_0$ are reference fluid speed, pressure and density, respectively. It is related to the Mach number $M_0$, i.e.\ the ratio between a typical fluid velocity and a typical speed of sound. For a $\gamma$-law gas, for example, it is $M = M_{0}\sqrt{\gamma}$. We are mainly interested in the so-called low Mach regime, i.e.\ $M \ll 1$, for which material waves are much slower than acoustic waves. The previous set of equations has to be completed by an equation of state (EOS). Further details on the EOS will be discussed in the upcoming Section \ref{ssec:eos}. The total energy $\rho E$ can be rewritten as $\rho E = \rho e + M^{2}\rho k$, where $e$ denotes the internal energy and $k = \left|\mathbf{u}\right|^{2}/2$ the kinetic energy. For the sake of convenience, we also introduce the specific enthalpy $h = e + p/\rho$, and we notice that the energy flux can be rewritten as
\begin{equation}
    \left(\rho E + p\right)\mathbf{u} = \left(e + M^{2}k + \frac{p}{\rho}\right)\rho\mathbf{u} = \left(h + M^{2}k\right)\rho\mathbf{u}.
\end{equation}
Hence, system \eqref{eq:euler_energy} reads equivalently as follows:
\begin{eqnarray}\label{eq:euler}
    \frac{\partial\rho}{\partial t} + \dive\left(\rho\mathbf{u}\right) &=& 0 \nonumber \\
    \frac{\partial\rho\mathbf{u}}{\partial t} + \dive\left(\rho\mathbf{u} \otimes \mathbf{u}\right) + \frac{1}{M^{2}}\grad p &=& \mathbf{0} \\
    \frac{\partial\rho E}{\partial t} + \dive\left[\left(h + M^{2}k\right)\rho\mathbf{u}\right] &=& 0. \nonumber 
\end{eqnarray}

\subsection{The equation of state}
\label{ssec:eos}	

System \eqref{eq:euler} has to be completed with an equation of state (EOS). In this work, we focus on the classical ideal gas law, the stiffened gas EOS (SG-EOS) \cite{metayer:2016}, and the general cubic EOS \cite{vidal:2001}. The equation that links together pressure, density, and internal energy for an ideal gas is given by \cite{vidal:2001}
\begin{equation}\label{eq:ideal_gas}
    p = \left(\gamma - 1\right)\rho e = \left(\gamma - 1\right)\left(\rho E - \frac{1}{2}M^{2}\rho\mathbf{u} \cdot \mathbf{u}\right). 
\end{equation}
Notice that \eqref{eq:ideal_gas} is valid only for a constant value of the specific heats ratio $\gamma$ \cite{vidal:2001}. The analogous relation for the SG-EOS reads as follows:
\begin{equation}\label{eq:sg_eos}
    p = \left(\gamma - 1\right)\left(\rho e - \rho q_{\infty}\right) - \gamma\pi_{\infty} = \left(\gamma - 1\right)\left(\rho E - \frac{1}{2}M^{2}\rho\mathbf{u} \cdot \mathbf{u} - \rho q_{\infty}\right) - \gamma\pi_{\infty},
\end{equation} 
with $q_{\infty}$ and $\pi_{\infty}$ representing constant parameters that determine the characteristics of the fluid. Notice that, for $q_{\infty} = \pi_{\infty} = 0$ in \eqref{eq:sg_eos}, we recover \eqref{eq:ideal_gas}. The last relation that we consider is the so-called general cubic EOS, for which the link between pressure, density, and temperature can be expressed as follows \cite[p.~221]{sandler:2017}, \cite[p.~119]{vidal:2001}:
\begin{equation}\label{eq:general_cubic_eos}
    p = \frac{\rho R_{g} T}{1 - \rho b} - \frac{a(T)\rho^{2}}{\left(1 - \rho b r_{1}\right)\left(1 - \rho b r_{2}\right)}.
\end{equation}
After some algebraic manipulations, \eqref{eq:general_cubic_eos} can be expressed as \cite[p.~222]{sandler:2017}
\begin{equation}
    z^{3} + c_{0}z^{2} + c_{1}z + c_{2} = 0,
\end{equation}
which is a cubic polynomial for the compressibility factor $z = p/(\rho R_{g}T)$. The parameters $c_{0}, c_{1}$, and $c_{2}$ depend on the thermodynamic state. Moreover, $R_{g} = R/m^{*}$ denotes the specific gas constant, with $R$ being the gas constant and $m^{*}$ the molar mass of the gas. Notice that for $a = b = 0$, the expression of the pressure of an ideal gas is recovered. As discussed in \cite{sandler:2017, vidal:2001}, it is convenient to express thermodynamic functions such as the internal energy or the enthalpy as the sum of a contribution due to the ideal gas and a \textit{residual} contribution due to non-ideality. Hence, after some manipulations, the equation linking together internal energy, density, and temperature, is given by \cite[p.~231]{sandler:2017}, \cite[p.~116]{vidal:2001}
\begin{equation}\label{eq:int_energy_general_cubic_eos}
    e = e^{\#}(T(p,\rho)) + \frac{a(T(p,\rho)) - T(p,\rho)\frac{da(T(p,\rho))}{dT}}{b}\frac{1}{r_{1} - r_{2}}\log\left(\frac{1 - \rho b r_{1}}{1 - \rho b r_{2}}\right).
\end{equation}
Here, $e^{\#}(T(p,\rho))$ denotes the internal energy of an ideal gas, which is function solely of the temperature $T$, $r_{1}$ and $r_{2}$ are suitable constants, whereas the parameters $a(T)$ and $b$ determine fluid characteristics \cite{vidal:2001}. More specifically, $a(T)$ is related to intermolecular forces, whereas $b$, the so called co-volume, takes into account the volume occupied by the molecules. Notice that, for $r_{1} \to 0$ and $r_{2} \to 0$, then $\frac{1}{r_{1} - r_{2}}\log\left(\frac{1 - \rho b r_{1}}{1 - \rho b r_{2}}\right) \to -\rho b$, which corresponds to the van der Waals EOS. For $r_{1} = -1 -\sqrt{2}, r_{2} = -1 + \sqrt{2}$, we get the Peng-Robinson EOS \cite{peng:1976}, \cite[p.~118]{vidal:2001}. If $c_{v} = \frac{de^{\#}}{dT}$ is constant, as in the case of calorically perfect gas, relation \eqref{eq:int_energy_general_cubic_eos} can be simplified
\begin{equation}\label{eq:int_energy_general_cubic_eos_constant}
    e = c_{v}T(p,\rho) + \frac{a(T(p,\rho)) - T(p,\rho)\frac{da(T(p,\rho))}{dT}}{b}\frac{1}{r_{1} - r_{2}}\log\left(\frac{1 - \rho b r_{1}}{1 - \rho b r_{2}}\right).
\end{equation}

\subsection{Asymptotic expansion}
\label{ssec:ap}

In this Section, we analyze the asymptotic limit of \eqref{eq:euler} as $M \to 0$. Consider the following expansion for density, velocity, and pressure, respectively:
\begin{eqnarray}
    \rho(\mathbf{x},t) &=& \bar{\rho}(\mathbf{x},t) + M\rho{'}(\mathbf{x},t) + M^{2}\rho{''}(\mathbf{x},t) + o(M^{2}) \label{eq:rho_expansion} \\
    \mathbf{u}(\mathbf{x},t) &=& \bar{\mathbf{u}}(\mathbf{x}, t) + M\mathbf{u}{'}(\mathbf{x},t) + M^{2}\mathbf{u}{''}(\mathbf{x},t) + o(M^{2}) \label{eq:u_expansion} \\
    p(\mathbf{x},t) &=& \bar{p}(\mathbf{x},t) + Mp{'}(\mathbf{x},t) + M^{2}p{''}(\mathbf{x},t) + o(M^{2}) \label{eq:p_expansion}
\end{eqnarray}
From now on, to simplify the notation, we omit the explicit dependence on space and time for all the variables. Substituting \eqref{eq:rho_expansion} and \eqref{eq:u_expansion} into the continuity equation in \eqref{eq:euler}, the leading order term relation is
\begin{equation}\label{eq:continuity_limit}
    \frac{\partial\bar{\rho}}{\partial t} + \dive\left(\bar{\rho}\bar{\mathbf{u}}\right) = 0.
\end{equation}
The leading order term relation for the momentum balance reduces to
\begin{equation}\label{eq:momentum_limit}
    \grad\bar{p} = \mathbf{0}.
\end{equation}
Hence, $\bar{p}$ is a function solely of time. Analogously, for the first order term, we obtain
\begin{equation}\label{eq:momentum_limit_first_order}
    \grad p' = \mathbf{0}.
\end{equation}
In addition, the second order term relation reads as follows:
\begin{equation}\label{eq:momentum_limit_second_order}
    \frac{\partial\bar{\rho}\bar{\mathbf{u}}}{\partial t} + \dive\left(\bar{\rho}\bar{\mathbf{u}} \otimes \bar{\mathbf{u}}\right) + \grad p{''} = \mathbf{0},
\end{equation}
where $p{''}$ represents a dynamical pressure \cite{cordier:2012, klein:1995, thomann:2019}. Finally, the leading order term relation for the energy equation reads as follows:
\begin{equation}\label{eq:energy_limit}
    \frac{\partial\bar{\rho}\bar{e}}{\partial t} + \dive\left(\bar{\rho}\bar{h}\bar{\mathbf{u}}\right) = 0,
\end{equation}
where $\bar{e} = e\left(\bar{\rho}, \bar{p}\right)$ and $\bar{h} = h\left(\bar{\rho}, \bar{p}\right)$. Since $\bar{\rho}\bar{e} = \bar{\rho}\bar{h} - \bar{p}$, we obtain
\begin{equation}
    \frac{\partial\bar{\rho}\bar{h}}{\partial t} - \frac{\partial\bar{p}}{\partial t} + \bar{\mathbf{u}} \cdot \grad\left(\bar{\rho}\bar{h}\right) + \bar{\rho}\bar{h}\left(\dive\mathbf{u}\right) = 0,
\end{equation}
or equivalently, considering $\bar{\rho}\bar{h} = \left(\bar{\rho}\bar{h}\right)\left(\bar{\rho}, \bar{p}\right)$,
\begin{equation}
    \frac{\partial\bar{\rho}\bar{h}}{\partial\bar{\rho}}\left(\frac{\partial\bar{\rho}}{\partial t} + \bar{\mathbf{u}} \cdot \grad\bar{\rho}\right) + \frac{\partial\bar{\rho}\bar{h}}{\partial\bar{p}}\left(\frac{\partial\bar{p}}{\partial t} + \bar{\mathbf{u}} \cdot \grad\bar{p}\right) - \frac{\partial\bar{p}}{\partial t} + \bar{\rho}\bar{h}\left(\dive\bar{\mathbf{u}}\right) = 0. 
\end{equation}
Thanks to \eqref{eq:continuity_limit} and \eqref{eq:momentum_limit}, we obtain
\begin{equation}
    \left(\bar{\rho}\bar{h} - \frac{\partial\bar{\rho}\bar{h}}{\partial\bar{\rho}}\bar{\rho}\right)\left(\dive\bar{\mathbf{u}}\right) + \left(\frac{\partial\bar{\rho}\bar{h}}{\partial\bar{p}} - 1\right)\frac{d\bar{p}}{dt} = 0, 
\end{equation}
or, since
$$\bar{\rho}\bar{h} - \frac{\partial\bar{\rho}\bar{h}}{\partial\bar{\rho}}\bar{\rho} = -\bar{\rho}^{2}\frac{\partial\bar{h}}{\partial\bar{\rho}}$$
and
$$\frac{\partial\bar{\rho}\bar{h}}{\partial\bar{p}} - 1 = \frac{\partial\bar{\rho}\bar{e}}{\partial\bar{p}},$$
equivalently,
\begin{equation}\label{eq:energy_limit_incomp}
    -\bar{\rho}^{2}\frac{\partial\bar{h}}{\partial\bar{\rho}}\left(\dive\bar{\mathbf{u}}\right) + \frac{\partial\bar{\rho}\bar{e}}{\partial\bar{p}}\frac{d\bar{p}}{dt} = 0. 
\end{equation}
We assume $\frac{\partial\bar{h}}{\partial\bar{\rho}} \neq 0$, as it holds away from vacuum. If $\frac{d\bar{p}}{dt} = 0$, we recover the incompressibility constraint
\begin{equation}
    \dive\bar{\mathbf{u}} = 0.
\end{equation}
Summing up, the asymptotic limit of \eqref{eq:euler} is
\begin{eqnarray}\label{eq:euler_ap_incomp}
    \frac{\partial\bar{\rho}}{\partial t} + \dive\left(\bar{\rho}\bar{\mathbf{u}}\right) &=& 0 \nonumber \\
    \grad\bar{p} &=& \mathbf{0} \nonumber \\
    \grad p{'} &=& \mathbf{0} \\
    \frac{\partial\bar{\rho}\bar{\mathbf{u}}}{\partial t} + \dive\left(\bar{\rho}\bar{\mathbf{u}} \otimes \bar{\mathbf{u}}\right) + \grad p{''} &=& \mathbf{0} \nonumber \\
    -\bar{\rho}^{2}\frac{\partial\bar{h}}{\partial\bar{\rho}}\left(\dive\bar{\mathbf{u}}\right) + \frac{\partial\bar{\rho}\bar{e}}{\partial\bar{p}}\frac{\partial\bar{p}}{\partial t} &=& 0. \nonumber
\end{eqnarray}
System \eqref{eq:euler_ap_incomp} represents a more general asymptotic limit with respect to the EOS and the boundary conditions of the compressible Euler equations for vanishing Mach number \cite{dellacherie:2010, orlando:2025}. We notice that we can rewrite the last equation of system \eqref{eq:euler_ap_incomp} as follows:
\begin{equation}\label{eq:dive_u}
    \dive\bar{\mathbf{u}} = \left(\frac{\partial\bar{\rho}\bar{e}}{\partial\bar{p}}\right)
    \left({\bar{\rho}^{2}\frac{\partial\bar{h}}{\partial\bar{\rho}}}\right)^{-1}
    \frac{d\bar{p}}{dt} = -\frac{1}{\bar{\rho}c^{2}\left(\bar{\rho}, \bar{p}\right)}\frac{d\bar{p}}{dt},
\end{equation}
where $c$ denotes the speed of sound. Indeed, from the definition of $\bar{h}$, one has
$$\left({\bar{\rho}^{2}\frac{\partial\bar{h}}{\partial\bar{\rho}}}\right) \left(\frac{\partial\bar{\rho}\bar{e}}{\partial\bar{p}}\right)^{-1} = \rho\left(\pad{\bar{e}}{\bar{\rho}} - \frac{\bar{p}}{\bar{\rho}^2}\right)\left(\frac{\partial\bar{e}}{\partial\bar{p}}\right)^{-1}.$$
From the first principle of thermodynamics, denoting by $s$ the specific entropy, one has
$$T ds = de + p d\frac{1}{\rho} = \pad{e}{\rho}d\rho + \pad{e}{p}dp - \frac{p}{\rho^2}d\rho = \left(\pad{e}{\rho} - \frac{p}{\rho^2}\right)d\rho + \pad{e}{p}dp.$$
As $c^{2} = \frac{\partial p}{\partial \rho}\bigg\rvert_{s}$, one has
\begin{equation}\label{eq:speed_sound}
    c^{2} = 
    \frac{\partial p}{\partial \rho}\bigg\rvert_{s} = 
    \left(\frac{\partial e}{\partial p}\right)^{-1}
    \left(\frac{p}{\rho^{2}} - \frac{\partial e}{\partial\rho}\right)
     =
    -\left(\frac{\partial e}{\partial p}\right)^{-1}
    {\frac{\partial h}{\partial \rho}},
\end{equation}
and therefore 
$$\frac{1}{\rho c^2} = -\frac{1}{\rho}\left(\pad{e}{p}\right)^{-1}
\pad{h}{\rho},$$
which proves \eqref{eq:dive_u} when applied to the lowest order terms in the asymptotic expansions \eqref{eq:rho_expansion}, \eqref{eq:p_expansion}. For more details consult \cite{orlando:2025, vidal:2001, whitham:2011}. Under periodic or free-slip boundary conditions, thanks to the divergence theorem, we have
$$\int_{\Omega} \dive\bar{\mathbf{u}}d\Omega = 0,$$
so that, by integrating \eqref{eq:dive_u} on $\Omega$, we find $\frac{d\bar{p}}{dt} = 0$. On the other hand, as one can easily notice from \eqref{eq:dive_u}, a time-dependent pressure with large amplitude variations imposed by a Dirichlet outflow boundary condition leads to a non-incompressible flow, i.e. $\dive\bar{\mathbf{u}} \neq 0$ and depending on the specific EOS. Consider, e.g., the ideal gas law \eqref{eq:ideal_gas}. We get
\begin{equation}
    \frac{\partial\bar{\rho}\bar{e}}{\partial\bar{p}} = \frac{1}{\gamma - 1} \qquad \bar{\rho}^{2}\frac{\partial\bar{h}}{\partial\bar{\rho}} = -\frac{\gamma}{\gamma - 1}\bar{p},
\end{equation}
so that \eqref{eq:dive_u} reduces to
\begin{equation}\label{eq:dive_u_IG}
    \dive\bar{\mathbf{u}} = -\frac{1}{\gamma}\frac{d\log\bar{p}}{dt}.
\end{equation}
Hence, the compressibility of a fluid described by the ideal gas law \eqref{eq:ideal_gas} is uniform in space and changes only in time. This is no longer valid for a general EOS \cite{orlando:2025}.

\section{The numerical framework}
\label{sec:num_disc}

In the low Mach number limit, pressure gradients terms, which are proportional to $1/M^{2}$, yield stiff components for the resulting semi-discretized ODE system \cite{casulli:1984, munz:2003, orlando:2022}. Implicit-Explicit Runge-Kutta (IMEX-RK) methods \cite{boscarino:2016, kennedy:2003} are widely employed for ODE systems that include both stiff and non-stiff components, to which the implicit and explicit schemes are applied, respectively. Therefore, an implicit coupling between the energy equation and the momentum one is appropriate, while the continuity equation can be treated in a fully explicit fashion. The spatial discretization is based on the Discontinuous Galerkin (DG) method, which easily allows for high-order accuracy. We refer to \cite{giraldo:2020} for a general introduction to the method. In this Section, we review some well-known concepts of IMEX-RK schemes. Then, we present the IMEX and the SI-IMEX time discretization for system \eqref{eq:euler}, respectively. Finally, some details concerning the spatial discretization will be also discussed.

\subsection{IMEX Runge-Kutta schemes}
\label{ssec:IMEX_review}

Implicit-Explicit Runge–Kutta (IMEX-RK) methods find extensive application in the numerical solution of PDEs, such as hyperbolic systems with relaxations \cite{boscarino:2024b, pareschi:2005}, convection–diffusion equations \cite{ascher:1997} and convection–diffusion–reaction equations \cite{kennedy:2003, kennedy:2019}. Let us start considering the following initial value problem for a system of ODE's
\begin{equation}\label{eq:IMEX_problem}
    \frac{d\mathbf{y}}{dt} = \mathbf{f}_{E}\left(\mathbf{y}, t\right) + \mathbf{f}_{I}\left(\mathbf{y}, t\right), \quad \mathbf{y}(0) = \mathbf{y}_{0},
\end{equation}
where $\mathbf{y}(t) \in \mathbb{R}^{M}, M \ge 1$ and we assume that $\mathbf{f}_{E}$ and $\mathbf{f}_{I} : \mathbb{R}^{M} \times \mathbb{R} \to \mathbb{R}^{M}$ are Lipschitz functions of $\mathbf{y}(t)$. We assume that the term $\mathbf{f}_{I}$ is {\em stiff} and the term $\mathbf{f}_{E}$ is {\em non-stiff}.

An $s$-stage IMEX-RK scheme applied to system \eqref{eq:IMEX_problem} takes the form:
\begin{subequations}\label{eq:IMEX_form}
\begin{eqnarray}
    \mathbf{v}^{(l)} &=& \mathbf{v}^{n} + \Delta t\sum_{m=1}^{l - 1}\tilde{a}_{lm}\mathbf{f}_{E}\left(t^{n} + \tilde{c}_{m}\Delta t, \mathbf{v}^{(m)}\right) + \Delta t\sum_{m=1}^{s}a_{lm}\mathbf{f}_{I}\left(t^{n} + c_{m}\Delta t, \mathbf{v}^{(m)}\right), \label{eq:stage_imex} \\
    \mathbf{v}^{n+1} &=& \mathbf{v}^{n} + \Delta t\sum_{l = 1}^{s}\tilde{b}_{l}\mathbf{f}_{E}\left(t^{n} + \tilde{c}_{l}\Delta t, \mathbf{v}^{(l)} \right) + \Delta t\sum_{l = 1}^{s}b_{l}\mathbf{f}_{I}\left( t^{n} + c_{l}\Delta t, \mathbf{v}^{(l)}\right). \label{eq:update_imex}
\end{eqnarray}
\end{subequations}
where the quantities $\mathbf{v}^{(l)}$ for $l = 1, \dots, s$, are called internal stages and approximate the exact solution $\mathbf{y}(t)$, at time $t = t^{n} + c_{l} \Delta t$, whereas $\mathbf{v}^{n+1}$ is the numerical solution that approximates the exact solution $\mathbf{y}(t)$ at time $t = t^{n} + \Delta t$. An $s$-stage IMEX-RK method is defined by two $s \times s$ real matrices $\tilde{\mathbf{A}} = \left\{\tilde{a}_{lm}\right\}$ and $\mathbf{A} = \left\{a_{lm}\right\}$, where the matrix $\tilde{\mathbf{A}}$ corresponds to the explicit method and is a lower triangular matrix with zero diagonal, i.e., $\tilde{a}_{lm}= 0$ for $l \le m$, while $\mathbf{A}$ is the one corresponding to the implicit scheme. We consider Diagonally Implicit Runge-Kutta (DIRK) methods for the implicit scheme so that $a_{lm} = 0$ for $l < m$. The use of a DIRK method for the treatment of $\mathbf{f}_{I}$ provides a sufficient condition to guarantee that the function $\mathbf{f}_{E}$ is always evaluated explicitly. The method is also characterized by the quadrature nodes $\tilde{\mathbf{c}} = \left(0, \tilde{c}_{2}, \dots \tilde{c}_{s}\right)^{\top}$, $\mathbf{c} = \left({c}_{1}, {c}_{2}, ..., {c}_{s}\right)^{\top}$, given by the usual relation
\begin{equation}\label{eq:imex_compatibility}
    \sum_{m=1}^{s}a_{lm} = c_{l} \qquad \sum_{m=1}^{s}\tilde{a}_{lm} = \tilde{c}_{l}, \quad l = 1,\dots,s,
\end{equation}
and by the weights $\tilde{\mathbf{b}}^{\top}= \left(\tilde{b}_{1}, \tilde{b}_{2}, \dots, \tilde{b}_{s}\right)$ and ${\mathbf{b}}^{\top}= \left(b_{1}, b_{2}, \dots, b_{s}\right)$ in $\mathbb{R}^{s}$. IMEX-RK methods can be represented in the usual Butcher notation \cite{butcher:2008}
\begin{center}
    \begin{tabular}{c|c}
	$\tilde{\mathbf{c}}$ & $\tilde{\mathbf{A}}$  \\
	\hline \\[-.3cm]
	& $\tilde{\mathbf{b}}^{\top}$
    \end{tabular}
    \qquad
    \begin{tabular}{c|c}
	$\mathbf{c}$ & $\mathbf{A}$ \\
	\hline \\[-.3cm]
	& $\mathbf{b}^{\top}$
    \end{tabular}.
\end{center}
Notice that the relation \eqref{eq:imex_compatibility} is a usual assumption for Runge-Kutta methods \cite{hairer:1993}.

It is useful to characterize the different IMEX-RK methods presented in the literature in two main types according to the structure of the matrix of the DIRK method. Following \cite{boscarino:2024b}, we have

\begin{definition}
    An IMEX-RK method is said to be of \textbf{type I} \cite{boscarino:2016, pareschi:2005} if the matrix $\mathbf{A}$ is invertible. It is said to be of \textbf{type II} \cite{boscarino:2016, kennedy:2003} if the matrix $\mathbf{A}$ can be written in the form
    $$\mathbf{A} = 
    \begin{pmatrix}
	0 & 0 \\
	\mathbf{a} & \boldsymbol{\mathcal{A}}
    \end{pmatrix},$$
    with $\mathbf{a} = (a_{21}, \dots, a_{s1})^{\top} \in  \mathbb{R}^{s-1}$ and the matrix $\boldsymbol{\mathcal{A}} \in \mathbb{R}^{(s-1) \times (s-1)}$ is invertible. In the special case $\mathbf{a} = 0$, $b_{1} = 0$, the method is said of \textbf{type ARS} (see \cite{ascher:1997}) and the DIRK method is reducible to a method using $s-1$ stages.
\end{definition}

Schemes of type II allow some simplifying assumptions, that make order conditions easier to treat for the construction of higher order schemes \cite{kennedy:2019}. On the other hand, schemes of type I are more suited to a theoretical analysis \cite{boscarino:2007, boscarino:2009} because of the invertibility of $\mathbf{A}$.

\begin{definition}
    We call an IMEX-RK method {\em stiffly accurate} (SA) if the corresponding DIRK method is stiffly accurate, namely \cite{wanner:1996}
    \begin{equation}
	a_{sl} = b_{l}, \quad i = 1, \dots, s.
    \end{equation}
\end{definition}
\noindent
All the IMEX-RK schemes employed in this work for the numerical simulations are stiffly accurate (see Appendix \ref{app:IMEX_coeffs}). 

The asymptotic properties of IMEX-RK methods are strongly related to the \texttt{L}-stability of the implicit part of the scheme. An implicit Runge-Kutta scheme is said to be \texttt{L}-stable \cite{wanner:1996} if it is \texttt{A}-stable and $R(z) \to 0$ as $z \to \infty$, where $R(z)$ is the stability function of the DIRK scheme. Following the result in \cite{wanner:1996}, \texttt{L}-stability is typically obtained combining the \texttt{A}-stability property with the SA property. However, for methods of type II, this combination does not necessarily lead to a \texttt{L}-stable scheme for the implicit part, because the matrix $\mathbf{A}$ is not invertible \cite{boscarino:2009}. For SA schemes of type II, a supplementary condition is required to obtain  \texttt{L}-stability, i.e. \cite{boscarino:2009}
\begin{equation}\label{eq:Lstab_type_II}
    \mathbf{e}^{\top}_{s}\boldsymbol{\mathcal{A}}^{-1}\mathbf{a}  = \sum_{m=2}^{s}\hat{w}_{sm}{a}_{m1} = 0,
\end{equation}
where $\mathbf{e}^{\top}_{s} = \left(0, \dots, 0,1\right)^{\top}$ and $\hat{w}_{lm}$ denotes the elements of the inverse of $\boldsymbol{\mathcal{A}}$. One can easily verify that all the implicit companion methods reported in Appendix \ref{app:IMEX_coeffs} are \texttt{L}-stable. 

In the sequel, to identify the different IMEX-RK schemes, we shall use the notation $\left(s, \sigma, p\right)$, where $s$ is the number of function evaluations of the implicit companion method, $\sigma$ is the number of function evaluations of the explicit companion method, and $p$ is the order of the IMEX scheme. In this work, we employ second, third, and fourth order time discretization schemes (see Appendix \ref{app:IMEX_coeffs}).

\subsection{IMEX time discretization for the Euler equations}
\label{ssec:IMEX}

In this Section, we outline the IMEX time discretization for the Euler equations \eqref{eq:euler}. Following \cite{casulli:1984, dumbser:2016a}, we consider an implicit treatment of the pressure gradient term within the momentum equation and of the pressure work term in the energy equation, while the continuity equation is discretized in a fully explicit fashion. A generic stage reads therefore as follows \cite{orlando:2022, orlando:2025}:
\begin{eqnarray}\label{eq:stage_euler_IMEX}
    \rho^{(l)} &=& \rho^{n} - \Delta t\sum_{m=1}^{l-1}\tilde{a}_{lm}\dive\left(\mathbf{q}\right)^{(m)} \nonumber \\
    \mathbf{q}^{(l)} + \frac{1}{M^{2}}a_{ll}\Delta t\grad p^{(l)} &=& \mathbf{q}^{n} - \frac{\Delta t}{M^{2}}\sum_{m=1}^{l-1}a_{lm}\grad p^{(m)} - \Delta t\sum_{m=1}^{l-1}\tilde{a}_{lm}\dive\left(\mathbf{q} \otimes \mathbf{u}\right)^{(m)} \\
    \mathcal{E}^{(l)} + a_{ll}\Delta t\dive\left(h\mathbf{q}\right)^{(l)} &=& \mathcal{E}^{n} - \Delta t\sum_{m=1}^{l-1}a_{lm}\dive\left(h\mathbf{q}\right)^{(m)} - \Delta t M^{2}\sum_{m=1}^{l-1}\tilde{a}_{lm} \dive\left(k\mathbf{q}\right)^{(m)}, \nonumber
\end{eqnarray}
where
\begin{equation}
    \mathbf{q}^{(l)} = (\rho\mathbf{u})^{(l)} \qquad \mathbf{u}^{(l)} = \frac{\mathbf{q}^{(l)}}{\rho^{(l)}} \qquad \mathcal{E}^{(l)} = (\rho E)^{(l)}.
\end{equation}
Notice that, substituting formally $\mathbf{q}^{(l)}$ into the energy equation and taking into account the definitions $\rho E = \rho e + M^{2}\rho k$ and $h = e + p/\rho$, the following nonlinear Helmholtz-type equation for the pressure is obtained:
\begin{eqnarray}\label{eq:pressure_elliptic_IMEX}
    &&\hspace{0.6cm} \rho^{(l)}\left[e(p^{(l)},\rho^{(l)}) +  M^{2}k^{(l)}\right] - a_{ll}^{2}\frac{\Delta t^{2}}{M^{2}}\dive\left[\left(e(p^{(l)},\rho^{(l)}) + \frac{p^{(l)}}{\rho^{(l)}}\right)\nabla p^{(l)}\right] \nonumber \\
    &&+\hspace{0.2cm} a_{ll}\Delta t\dive\left[\left(e(p^{(l)},\rho^{(l)}) + \frac{p^{(l)}}{\rho^{(l)}}\right)\mathbf{m}^{(l)}\right]
    = \hat{e}^{(l)}, 
\end{eqnarray}
where
\begin{subequations}
\begin{eqnarray}
    \mathbf{m}^{(l)} &=& \mathbf{q}^{n} - \frac{\Delta t}{M^{2}}\sum_{m=1}^{l-1}a_{lm}\grad p^{(m)} - \Delta t \sum_{m=1}^{l-1}\tilde{a}_{lm}\dive\left(\mathbf{q} \otimes \mathbf{u}\right)^{(m)},\\
    \hat{e}^{(l)} &=& \mathcal{E}^{n} - \Delta t\sum_{m=1}^{l-1}a_{lm}\dive\left(h\mathbf{q}\right)^{(m)} - M^{2}\Delta t\sum_{m=1}^{l-1}\tilde{a}_{lm} \dive\left(k\mathbf{q}\right)^{(m)}.
    \end{eqnarray}
\end{subequations}
Equation \eqref{eq:pressure_elliptic_IMEX} is solved through a fixed point procedure \cite{dumbser:2016a, orlando:2022}. More specifically, setting $ \xi^{(0)} = p^{(l-1)}, k^{(l,0)} = k^{(l-1)}, $ one solves for $\tilde{k} = 0, \dots, L$ the equation
\begin{eqnarray}\label{eq:pressure_elliptic_IMEX_fixed}
    \rho^{(l)}e(\xi^{(\tilde{k}+1)},\rho^{(l)})   
    &-& a_{ll}^{2}\frac{\Delta t^{2}}{M^{2}}\dive\left[\left(e(\xi^{(\tilde{k})},\rho^{(l)})  + \frac{\xi^{(\tilde{k})}}{\rho^{(l)}}\right)\nabla\xi^{(\tilde{k}+1)}\right] \nonumber \\
    &=& \hat{e}^{(l)} - M^{2}\rho^{(l)}k^{(l,\tilde{k})} - a_{ll}\Delta t\dive\left[\left(e(\xi^{(\tilde{k})},\rho^{(l)}) + \frac{\xi^{(\tilde{k})}}{\rho^{(l)}}\right)\mathbf{m}^{(l)}\right] \nonumber
\end{eqnarray}
and then updates the velocity as 
$$\mathbf{u}^{(l,\tilde{k}+1)} + \frac{a_{ll}\Delta t}{\rho^{(l)}M^2}\nabla \xi^{(\tilde{k}+1)} = \frac{\mathbf{m}^{(l)}}{\rho^{(l)}}.$$
As already discussed in \cite{dumbser:2016a}, solving directly \eqref{eq:pressure_elliptic_IMEX} keeping a full implicit treatment of the enthalpy as in a classical Newton method yields a system strongly nonlinear and difficult to control. For this purpose, one adopts a Picard iteration technique in which the contribution of the enthalpy is computed at the previous fixed point iteration so as to reduce the nonlinearity of \eqref{eq:pressure_elliptic_IMEX}. Moreover, this choice is justified by the fact that two/three iterations are typically sufficient to obtain a satisfactory solution, as already observed in \cite{casulli:2010, dumbser:2016a} and as further confirmed by our numerical experiments (see in particular Section \ref{ssec:fixed_point}).

\subsection{Semi-Implicit IMEX (SI-IMEX) time discretization for the Euler equations}
\label{ssec:SI_IMEX}

The procedure outlined in the previous Section always requires the solution of a nonlinear system at each stage. In this Section, we outline the semi-implicit IMEX (SI-IMEX) time discretization for the Euler equations \eqref{eq:euler} \cite{boscarino:2022}, which provides similar results with less computational time. The governing partial differential equations \eqref{eq:euler} can be cast into a compact and general form as
\begin{equation}\label{eq:SI_system}
    \frac{\partial\mathbf{U}}{\partial t} = \mathbf{H}\left(\mathbf{U}, \mathbf{U}\right).
\end{equation}
Here $\mathbf{U} = \left(\rho, \rho\mathbf{u}, \rho E\right)^{\top}$ is the vector of conserved variables, and
\begin{equation}
    \mathbf{H}\left(\mathbf{U}, \mathbf{U}\right) = -\dive \left(
    \begin{array}{c}
 	\rho\mathbf{u} \\
 	\rho\mathbf{u} \otimes \mathbf{u} \\
 	M^{2}\rho k\mathbf{u} 
    \end{array}
    \right)  
    -\dive \left(
 	\begin{array}{c}
 		0 \\
 		\frac{1}{M^{2}}p \\
 		h \rho \mathbf{u}
 	\end{array}
    \right),
\end{equation}
where $\mathbf{H}: \mathbb{R}^{d + 2}  \times \mathbb{R}^{d + 2} \to \mathbb{R}^{d + 2}$ is a sufficiently regular mapping. Following \cite{boscarino:2016}, the governing partial differential equations \eqref{eq:euler} are written under the form of an autonomous system \eqref{eq:SI_system}, for all $t > t_{0}$ with the initial condition $\mathbf{U}(t_{0}) = \mathbf{U}_{0}$.

We refer to semi-implicit schemes as numerical methods that address problems of the form \eqref{eq:SI_system}, wherein the variable $\mathbf{U}$, appearing as first argument of $\mathbf{H}$, is treated explicitly and will be denoted by $\mathbf{U}_{E}$, while the variable $\mathbf{U}$ appearing as the second argument is treated implicitly and denoted by $\mathbf{U}_{I}$. Thus, we obtain a partitioned system of the form
\begin{equation}\label{eq:SI_RK}
    \begin{split}
        \frac{\partial\mathbf{U}_{E}}{\partial t} &= \mathbf{H}(\mathbf{U}_{E}, \mathbf{U}_{I}), \\
        \frac{\partial\mathbf{U}_{I}}{\partial t} &= \mathbf{H}(\mathbf{U}_{E}, \mathbf{U}_{I}).
    \end{split}
\end{equation}
with
\begin{equation}\label{eq:Hfun}
    \mathbf{H}(\mathbf{U}_{E}, \mathbf{U}_{I}) = -\dive\mathbf{F}_{E} - \dive\mathbf{F}_{SI},
\end{equation}
and
\begin{equation}
    \mathbf{F}_{E} = \left(
    \begin{array}{c}
	(\rho\mathbf{u})_{E} \\
	(\rho\mathbf{u} \otimes \mathbf{u})_{E} \\
	M^{2}(\rho k\mathbf{u})_{E} 
    \end{array}
    \right), \quad 
    \mathbf{F}_{SI} = \left(
    \begin{array}{c}
	0 \\
	\frac{1}{M^{2}}p_{I} \\
	h_{E}(\rho\mathbf{u})_{I}
    \end{array}
    \right).
\end{equation}
Moreover, $\mathbf{U}_{E} = \left(\rho_{E}, (\rho\mathbf{u})_{E}, (\rho E)_{E}\right)^{\top}$ and $\mathbf{U}_{I} = \left(\rho_{I}, (\rho\mathbf{u})_{I}, (\rho E)_{I}\right)^{\top}$. Subscripts $E$ and $SI$ in \eqref{eq:Hfun} indicate the explicit and semi-implicit treatment of the first and the second term, respectively. Notice that the number of unknowns in \eqref{eq:SI_RK} has been doubled. However, when specific time discretizations are chosen for autonomous systems, this doubling is only apparent \cite{boscarino:2016}. Finally, the kinetic energy in the total energy definition splits into an explicit and an implicit contribution, namely:
\begin{equation}\label{eq:total_energy_SI}
    (\rho E)_{I} = (\rho e)_{I} + \frac{M^{2}}{2}\mathbf{u}_{E} \cdot (\rho \mathbf{u})_{I}.
\end{equation}

High-order time discretization is achieved making use of IMEX-RK schemes. More specifically, we adopt methods for which $\tilde{b}_{l} = b_{l}, \hspace{0.1cm} l = 1, \dots, s$. We observe that, since $\tilde{b}_{l} = b_{l}$, the numerical solutions are the same, i.e., if $\mathbf{U}_{E}^{0} = \mathbf{U}_{I}^{0}$, then $\mathbf{U}_{E}^{n} = \mathbf{U}_{I}^{n}$ for all $n > 0$. Hence, the duplication of the system is only apparent.

Under the assumption that system \eqref{eq:SI_RK} is autonomous, a SI-IMEX-RK method is obtained as follows. First, set $\mathbf{U}_{E}^{n} = \mathbf{U}_{I}^{n} = \mathbf{U}^{n}$. Then, the internal stage values read
\begin{equation}
    \begin{split}
        \mathbf{U}_{E}^{(l)} &= \mathbf{U}_{E}^{n} + \Delta t\sum_{m=1}^{l-1} \tilde{a}_{lm}\mathbf{H}\left(\mathbf{U}_{E}^{(m)}, \mathbf{U}_{I}^{(m)}\right) \\
        \mathbf{U}_{I}^{(l)} &= \mathbf{U}_{I}^{n} + \Delta t\sum_{m=1}^{l-1} a_{lm}\mathbf{H}\left(\mathbf{U}_{E}^{(m)}, \mathbf{U}_{I}^{(m)}\right) + \Delta t a_{ll}\mathbf{H}\left(\mathbf{U}_{E}^{(l)}, \mathbf{U}_{I}^{(l)}\right), 
    \end{split}
\end{equation}
for $l = 1, \dots, s$. Finally, the numerical solution is updated with
\begin{equation}\label{eq:update_SI_IMEX}
    \mathbf{U}^{n+1} = \mathbf{U}^{n} + \Delta t\sum_{l=1}^{s}b_{l}\mathbf{H}\left(\mathbf{U}_{E}^{(l)}, \mathbf{U}_{I}^{(l)}\right).
\end{equation}

For the sake of clarity in the notation, we denote
$$\mathbf{q} = \rho\mathbf{u} \qquad \mathcal{E} = \rho E$$
as done in the previous Section. As an example, we present the first order semi-implicit scheme solving system \eqref{eq:euler} to compute the numerical solution $\mathbf{U}^{n+1} = \left(\rho^{n+1}, \mathbf{q}^{n+1}, \mathcal{E}^{n+1}\right)^{\top}$. We focus on the time discretization, while keeping the space continuous. We consider the first order IMEX-RK scheme
\begin{center}
    \begin{tabular}{c|c}
	0 & 0  \\
	\hline \\[-.3cm]
	& 1
    \end{tabular}
    \qquad
    \begin{tabular}{c|c}
	1& 1 \\
	\hline \\[-.3cm]
	& 1
    \end{tabular}.
\end{center}
Formally applying the above tableau to the partitioned system \eqref{eq:SI_RK}, it reads
\begin{equation}
    \begin{split}
	\mathbf{U}_{E}^{(1)} &= \mathbf{U}^{n} \\
        \mathbf{U}_{I}^{(1)} &= \mathbf{U}^{n} + \Delta t\mathbf{H}\left(\mathbf{U}^{n}, \mathbf{U}_{I}^{(1)}\right), \\
        \mathbf{U}^{n+1} &= \mathbf{U}_{I}^{(1)},
    \end{split}
\end{equation}
and explicitly we get
\begin{equation}
    \mathbf{U}^{n+1} = \mathbf{U}^{n} - \Delta t\dive\mathbf{F}_{E}(\mathbf{U}^{n}) - \Delta t\dive\mathbf{F}_{SI}(\mathbf{U}^{n}, \mathbf{U}^{n+1})
\end{equation}
with
$$\mathcal{E}^{n+1} = (\rho e)^{n+1} + \frac{M^{2}}{2}\mathbf{u}^{n} \cdot \mathbf{q}^{n+1}.$$

A generic stage $l$ for the Euler equations using a SI-IMEX-RK scheme reads therefore as follows.
\paragraph{Explicit step.} Set
\begin{eqnarray}\label{eq:stage_explicit_euler_SI_IMEX}
    \rho_{E}^{(l)} &=& \rho^{n} - \Delta t\sum_{m=1}^{l-1}\tilde{a}_{lm}\dive\mathbf{q}_{E}^{(m)} \nonumber \\
    \mathbf{q}_{E}^{(l)} &=& \mathbf{q}^{n} - \Delta t\sum_{m=1}^{l-1}\tilde{a}_{lm}\left[\dive\left(\mathbf{q} \otimes \mathbf{u}\right)_{E}^{(m)} + \frac{1}{M^{2}}\grad p_{I}^{(m)}\right] 
    \\
    \mathcal{E}_{E}^{(n,l)} &=& \mathcal{E}^{n} \nonumber - \Delta t\sum_{m=1}^{l-1}\tilde{a}_{lm}\left[M^{2}\dive\left(k\mathbf{q}\right)_{E}^{(m)} + \dive\left(h_{E}\mathbf{q}_{I}\right)^{(m)}\right]. \nonumber
\end{eqnarray}
\paragraph{Implicit step.} Solve $\mathbf{U}_{I}^{(l)}$:
\begin{eqnarray}\label{eq:stage_implicit_euler_SI_IMEX}
    {\rho}_{I}^{(l)} &=& \bar{\rho}_{I}^{(l)} - {a}_{ll}\Delta t \dive\mathbf{q}_{E}^{(l)}, \nonumber\\
    \mathbf{q}_{I}^{(l)} + \frac{1}{M^{2}}{a}_{ll}\Delta t\grad p_{I}^{(l)} &=& \mathbf{m}_{I}^{(l)} - a_{ll}\Delta t\dive\left(\mathbf{q} \otimes \mathbf{u}\right)_{E}^{(l)} \\
    \mathcal{E}_{I}^{(l)} + {a}_{ll}\Delta t\dive\left(h_{E}\mathbf{q}_{I}\right)^{(l)} &=& \hat{e}_{I}^{(l)} - a_{ll}\Delta t M^{2}\dive\left(k\mathbf{q}\right)_{E}^{(l)}, \nonumber 
\end{eqnarray}
where 
\begin{eqnarray}
    \bar{\rho}_{I}^{(l)} &=& \rho^{n} - \Delta t\sum_{m=1}^{l-1}{a}_{lm}\dive\mathbf{q}_{E}^{(m)} \nonumber \\
    \mathbf{m}_{I}^{(l)} &=& \mathbf{q}^{n} - \Delta t\sum_{m=1}^{l-1}{a}_{lm}\left[\dive\left(\mathbf{q} \otimes \mathbf{u}\right)_{E}^{(m)} + \frac{1}{M^{2}}\grad p_{I}^{(m)}\right] \\
    \hat{e}_{I}^{(l)} &=& \mathcal{E}^{n} - \Delta t\sum_{m=1}^{l-1}{a}_{lm} \left[M^{2}\dive\left(k\mathbf{q}\right)_{E}^{(m)} - \dive\left(h_{E}\mathbf{q}_{I}\right)^{(m)}\right]. \nonumber
\end{eqnarray}
To solve system \eqref{eq:stage_implicit_euler_SI_IMEX}, we substitute $\mathbf{q}_{I}^{(l)}$ in the energy equation, so as to obtain an elliptic equation for $p_{I}^{(l)}$ \cite{boscarino:2022}, which reads as follows:
\begin{eqnarray}\label{eq:pressure_elliptic_SI_IMEX}
    \rho_{I}^{(l)}e_{I}(p_{I}^{(l)},\rho_{I}^{(l)})   
    &-& a_{ll}^{2}\frac{\Delta t^{2}}{M^{2}}\dive\left(h_{E}^{(l)}\nabla p_{I}^{(l)}\right) - \frac{1}{2} a_{ll}\Delta t \hspace{0.05cm} \mathbf{u}_{E}^{(l)} \cdot \grad p_{I}^{(l)} \nonumber \\
    &=& \hat{e}_{I}^{(l)} - a_{ll}\Delta t M^{2}\dive\left(k\mathbf{q}\right)_{E}^{(l)} - \frac{M^{2}}{2}\mathbf{u}_{E}^{(l)} \cdot \left[\mathbf{m}_{I}^{(l)} - a_{ll}\Delta t\dive\left(\mathbf{q} \otimes \mathbf{u}\right)_{E}^{(l)}\right] \nonumber \\
    &-& a_{ll}\Delta t\dive\left[h_{E}^{(l)}\left[\mathbf{m}_{I}^{(l)} - a_{ll}\Delta t\dive\left(\mathbf{q} \otimes \mathbf{u}\right)_{E}^{(l)}\right]\right].
\end{eqnarray}	
Next, one computes $\mathbf{q}_{I}^{(l)}$ and $\mathcal{E}_{I}^{(l)}$ from \eqref{eq:stage_implicit_euler_SI_IMEX}. Finally, one updates the numerical solution $\mathbf{U}^{n+1}$ from \eqref{eq:update_SI_IMEX}.

\subsection{Impact of the EOS}
\label{ssec:num_eos}

Equations \eqref{eq:pressure_elliptic_IMEX_fixed} and \eqref{eq:pressure_elliptic_SI_IMEX} need the relation between the internal energy and the pressure before being solved. In the case of the ideal gas law \eqref{eq:ideal_gas}, since
$$\rho e = \frac{1}{\gamma - 1}p,$$
equations \eqref{eq:pressure_elliptic_IMEX_fixed} and \eqref{eq:pressure_elliptic_SI_IMEX} constitute a linear system for $\xi^{(\tilde{k} + 1)}$ and $p_{I}^{(l)}$, respectively. Analogous considerations hold for the SG-EOS \eqref{eq:sg_eos}, since
$$\rho e = \frac{p}{\gamma - 1} + \frac{\gamma\pi_{\infty}}{\gamma - 1} + \rho q_{\infty}.$$
Hence, only a supplementary term depending on the already updated density is present. 

On the other hand, the use of a more general equation of state, such as the general cubic EOS \eqref{eq:general_cubic_eos}, leads to a nonlinear relation between internal energy and pressure and therefore a nonlinear equation for the pressure should be solved \cite{boscheri:2021a}. We rewrite the term $\rho e$ as $\frac{\rho e}{p}p$, so that, following the discussion in \cite{brugnano:2008}, the nonlinear equation is solved by the following Picard iteration 
\begin{eqnarray}\label{eq:pressure_elliptic_IMEX_fixed_mod}
    \frac{\rho^{(l)}e(\xi^{(\tilde{k})},\rho^{(l)})}{\xi^{(\tilde{k})}}\xi^{(\tilde{k}+1)}   
    &-& a_{ll}^{2}\frac{\Delta t^{2}}{M^{2}}\dive\left[\left(e(\xi^{(\tilde{k})},\rho^{(l)})  + \frac{\xi^{(\tilde{k})}}{\rho^{(l)}}\right)\nabla\xi^{(\tilde{k}+1)}\right] \nonumber \\
    &=& \hat{e}^{(l)} - M^{2}\rho^{(l)}k^{(l,\tilde{k})} - a_{ll}\Delta t\dive\left[\left(e(\xi^{(\tilde{k})},\rho^{(l)}) + \frac{\xi^{(\tilde{k})}}{\rho^{(l)}}\right)\mathbf{m}^{(l)}\right].
\end{eqnarray}
It is worth noting that in the case of the ideal gas law \eqref{eq:ideal_gas},
$$\frac{\rho^{(l)}e(\xi^{(\tilde{k})},\rho^{(l)})}{\xi^{(\tilde{k})}}\xi^{(\tilde{k}+1)} = \frac{1}{\gamma - 1}\xi^{(\tilde{k}+1)} = \rho^{(l)}e(\xi^{(\tilde{k} + 1)}),$$
so that \eqref{eq:pressure_elliptic_IMEX_fixed_mod} reduces to \eqref{eq:pressure_elliptic_IMEX_fixed}. Notice also that \eqref{eq:pressure_elliptic_IMEX_fixed_mod} corresponds to a slightly different linearization with respect to the one proposed in \cite{orlando:2022}, which was tailored for the general cubic EOS, while \eqref{eq:pressure_elliptic_IMEX_fixed_mod} is applicable to a general EOS. 

Analogous considerations hold for \eqref{eq:pressure_elliptic_SI_IMEX}: one can rewrite
\begin{equation}
    (\rho e)_{I} = \frac{(\rho e)_{I}}{p_{I}}p_{I}
\end{equation}
and solve the resulting mild nonlinear equation according to the Picard iteration described in \cite{brugnano:2008}. Note that the fixed point procedure proposed in \cite{brugnano:2008} corresponds to a Newton-type method which can also be applied to non-differentiable relations, like those that could be obtained for tabulated EOS, unlike the standard Newton method. We refer to \cite{brugnano:2008} for a detailed description and analysis of the algorithm. Nevertheless, in the case of the semi-implicit time discretization, in order to avoid the solution of a nonlinear equation and the use of a fixed point loop, we approximate 
\begin{equation}
	\frac{(\rho e)_{I}}{p_{I}}p_{I} \approx \frac{(\rho e)_{E}}{p_{E}}p_{I},
\end{equation}
so that \eqref{eq:pressure_elliptic_SI_IMEX} modifies as
\begin{eqnarray}\label{eq:pressure_elliptic_SI_IMEX_mod}
    \frac{(\rho e)_{E}^{(l)}}{p_{E}^{(l)}}p_{I}^{(l)}   
    &-& a_{ll}^{2}\frac{\Delta t^{2}}{M^{2}}\dive\left(h_{E}^{(l)}\nabla p_{I}^{(l)}\right) - \frac{1}{2} a_{ll}\Delta t \hspace{0.05cm} \mathbf{u}_{E}^{(l)} \cdot \grad p_{I}^{(l)} \nonumber \\
    &=& \hat{e}_{I}^{(l)} - a_{ll}\Delta t M^{2}\dive\left(k\mathbf{q}\right)_{E}^{(l)} - \frac{M^{2}}{2}\mathbf{u}_{E}^{(l)} \cdot \left[\mathbf{m}_{I}^{(l)} - a_{ll}\Delta t\dive\left(\mathbf{q} \otimes \mathbf{u}\right)_{E}^{(l)}\right] \nonumber \\
    &-& a_{ll}\Delta t\dive\left[h_{E}^{(l)}\left[\mathbf{m}_{I}^{(l)} - a_{ll}\Delta t\dive\left(\mathbf{q} \otimes \mathbf{u}\right)_{E}^{(l)}\right]\right].
\end{eqnarray}
It is to be noted that, in the case of the ideal gas law \eqref{eq:ideal_gas},
$$\frac{(\rho e)_{E}^{(l)}}{p_{E}^{(l)}}p_{I}^{(l)} = \frac{1}{\gamma - 1}p_{I}^{(l)} = (\rho e)_{I}^{(l)},$$
so that \eqref{eq:pressure_elliptic_SI_IMEX_mod} reduces to \eqref{eq:pressure_elliptic_SI_IMEX}.

\subsection{The spatial discretization strategy}
\label{ssec:space_disc}

\begin{figure}[h!]
    \centering
    \includegraphics[width = 0.5\textwidth]{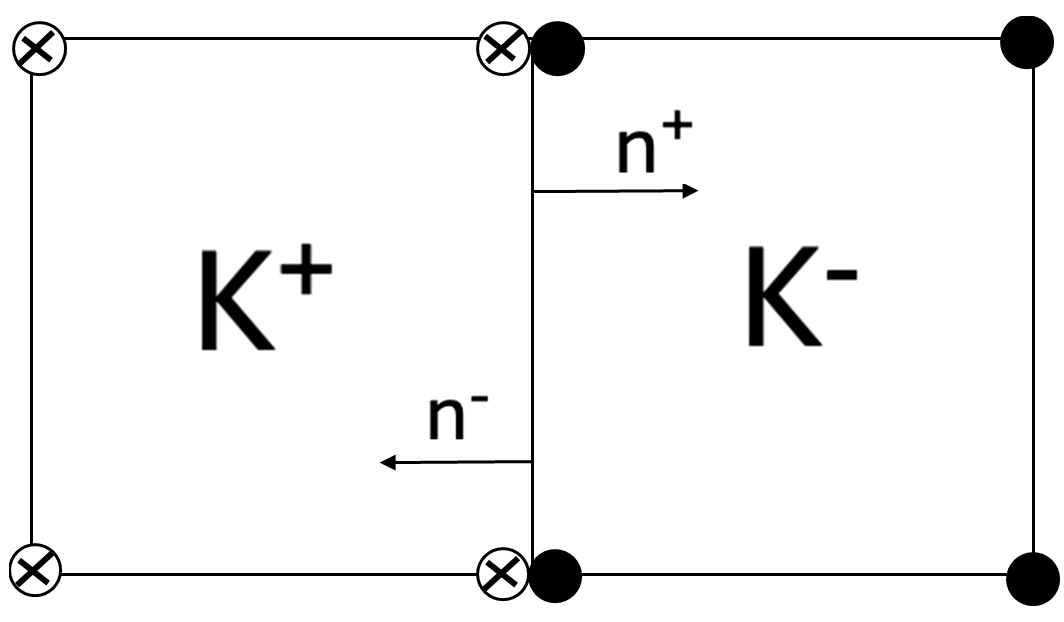}
    \caption{Example of two neighboring elements for a nodal DG formulation based on Lagrange polynomials. The nodes correspond to the support points of $\left(r + 1\right)$-order Gauss-Lobatto quadrature rule (in the image $r= 1$).}
    \label{fig:DG_mesh}
\end{figure}

In this Section, we briefly outline the spatial discretization adopted for \eqref{eq:stage_euler_IMEX} and \eqref{eq:stage_explicit_euler_SI_IMEX}-\eqref{eq:stage_implicit_euler_SI_IMEX}, which is based on the Discontinuous Galerkin (DG) method \cite{giraldo:2020} as implemented in the \texttt{deal.II} library \cite{arndt:2023, bangerth:2007}. We consider a decomposition of the domain $\Omega$ into a family of quadrilaterals $\mathcal{T}_{h}$ and denote each element by $K$. The skeleton $\mathcal{E} = \mathcal{E}^{I} \cup \mathcal{E}^{B}$ denotes the set of all the element faces, with $\mathcal{E}^{I}$ and $\mathcal{E}^{B}$ being the subset of interior and boundary faces, respectively. A face $\Gamma \in \mathcal{E}_{I}$ shares two elements, $K^{+}$ with outward unit normal $\mathbf{n}^{+}$, and $K^{-}$ with outward unit normal $\mathbf{n}^{-}$, while we simply denote by $\mathbf{n}$ the outward unit normal for a face $\Gamma \in \mathcal{E}^{B}$ (see Figure \ref{fig:DG_mesh}). For a scalar function $\varphi$, the jump is defined as
\begin{equation}
    \jump{\varphi} = \varphi^{+}\mathbf{n}^{+} + \varphi^{-}\mathbf{n}^{-} \text{ if } \Gamma \in \mathcal{E}^{I} \qquad \jump{\varphi} = \varphi\mathbf{n} \text{ if } \Gamma \in \mathcal{E}^{B},
\end{equation}
while the average reads
\begin{equation}
    \ave{\varphi} = \frac{1}{2}\left(\varphi^{+} + \varphi^{-}\right) \text{ if } \Gamma \in \mathcal{E}^{I} \qquad \ave{\varphi} = \varphi \text{ if } \Gamma \in \mathcal{E}^{B}.
\end{equation}
Analogous definitions apply for a vector function $\boldsymbol{\varphi}$. More specifically, we define
\begin{subequations}
    \begin{alignat}{2}
        \jump{\boldsymbol{\varphi}} &= \boldsymbol{\varphi}^{+} \cdot \mathbf{n}^{+} + \boldsymbol{\varphi}^{-} \cdot \mathbf{n}^{-} \text{ if } \Gamma \in \mathcal{E}^{I} \qquad && \jump{\boldsymbol{\varphi}} = \boldsymbol{\varphi} \cdot \mathbf{n} \text{ if } \Gamma \in \mathcal{E}^{B} \\
        \ave{\boldsymbol{\varphi}} &= \frac{1}{2}\left(\boldsymbol{\varphi}^{+} + \boldsymbol{\varphi}^{-}\right) \text{ if } \Gamma \in \mathcal{E}^{I} \qquad && \ave{\boldsymbol{\varphi}} = \boldsymbol{\varphi} \text{ if } \Gamma \in \mathcal{E}^{B}.
    \end{alignat}
\end{subequations}	
Finally, for vector functions, it is also useful to define a tensor jump as follows:
\begin{equation}
    \tjump{\boldsymbol{\varphi}} = \boldsymbol{\varphi}^{+} \otimes \mathbf{n}^{+} + \boldsymbol{\varphi}^{-} \otimes \mathbf{n}^{-} \text{ if } \Gamma \in \mathcal{E}^{I} \qquad \tjump{\boldsymbol{\varphi}} = \boldsymbol{\varphi} \otimes \mathbf{n} \text{ if } \Gamma \in \mathcal{E}^{B}.
\end{equation}
We also introduce the following finite element spaces
$$Q_{r} = \left\{v \in L^2(\Omega) : v\rvert_{K} \in \mathbb{Q}_{r} \quad \forall K \in \mathcal{T}_{\mathcal{H}}\right\} \qquad \mathbf{V}_{r} = \left[Q_{r}\right]^{d},$$
where $\mathbb{Q}_{r}$ is the space of polynomials of degree $r$ in each coordinate direction. We then denote by $\boldsymbol{\varphi}_{i}(\mathbf{x})$ the basis functions for the space $\mathbf{V}_{r}$ and by $\psi_i(\mathbf{x})$ the basis functions for the space $Q_{r}$, the finite element spaces chosen for the discretization of the velocity and of the pressure (as well as the density), respectively, so that
$$\mathbf{u} \approx \sum_{j = 1}^{\left|\mathcal{T}_{h}\right|(r+1)^{d}}u_{j}(t)\boldsymbol{\varphi}_j(\mathbf{x}) \qquad p \approx \sum_{j = 1}^{\left|\mathcal{T}_{h}\right|(r+1)^{d}}p_{j}(t)\psi_{j}(\mathbf{x}).$$
Here $\left|\mathcal{T}_{h}\right|$ denotes the number of elements of the computational mesh. Recall that $d$ denotes the space dimension. The number of degrees of freedom for scalar variable is indeed equal to $\left|\mathcal{T}_{h}\right|\left(r + 1\right)^{d}$ \cite{giraldo:2020}. The shape functions correspond to the products of Lagrange interpolation polynomials for the support points of $\left(r + 1\right)$-order Gauss-Lobatto quadrature rule in each coordinate direction (Figure \ref{fig:DG_mesh}). In particular, we have grid points at the boundaries of the elements, where the solution can be discontinuous and this simplifies the evaluation of the integrals at the boundary itself \cite{kronbichler:2019}. Hence, for any given edge, the only shape functions with non-zero values are exactly those whose node points are located on that edge \cite{kronbichler:2019}.

Given these definitions, the weak formulation for the momentum equation at each stage \eqref{eq:stage_euler_IMEX} reads as follows \cite{orlando:2022, orlando:2025}:
\begin{equation}
    \mathbf{A}^{(l)}\mathbf{U}^{(l)} + \mathbf{B}^{(l)}\mathbf{P}^{(l)} = \mathbf{F}^{(n,l)},
\end{equation}
with $\mathbf{U}^{(l)}$ denoting the vector of the degrees of freedom associated to the velocity field and $\mathbf{P}^{(l)}$ denoting the vector of the degrees of freedom associated to the pressure. Here we have set
\begin{eqnarray}
    A_{ij}^{(l)} &=& \sum_{K \in \mathcal{T}_{h}} \int_{K} \rho^{(l)}\boldsymbol{\varphi}_{j} \cdot \boldsymbol{\varphi}_{i}d\Omega \\
    B_{ij}^{(l)} &=& \sum_{K \in \mathcal{T}_{h}}\int_{K} -a_{ll}\frac{\Delta t}{M^{2}}\dive\boldsymbol{\varphi}_{i}\Psi_{j}d\Omega + \sum_{\Gamma \in \mathcal{E}}\int_{\Gamma}a_{ll}\frac{\Delta t}{M^{2}}\ave{\Psi_{j}}\jump{\boldsymbol{\varphi}_{i}}d\Sigma \\
    F_{i}^{(l)} &=& \sum_{K \in \mathcal{T}_{h}} \int_{K} \rho^{n}\mathbf{u}^{n} \cdot \boldsymbol{\varphi}_{i}d\Omega \nonumber \\
    &+& \sum_{m = 1}^{l-1}\sum_{K \in \mathcal{T}_{h}}\int_{K}\tilde{a}_{lm}\Delta t\left(\rho^{(m)}\mathbf{u}^{(m)} \otimes \mathbf{u}^{(m)}\right) : \grad\boldsymbol{\varphi}_{i}d\Omega + \sum_{m=1}^{l-1}\sum_{K \in \mathcal{T}_{h}}\int_{K}a_{lm}\frac{\Delta t}{M^{2}}p^{(m)}\dive\boldsymbol{\varphi}_{i}d\Omega \nonumber \\
    &-& \sum_{m=1}^{l-1}\sum_{\Gamma \in \mathcal{E}}\int_{\Gamma}\tilde{a}_{lm}\Delta t \ave{\rho^{(m)}\mathbf{u}^{(m)} \otimes \mathbf{u}^{(m)}} : \tjump{\boldsymbol{\varphi}_{i}}d\Sigma \nonumber \\
    &-& \sum_{m=1}^{l-1}\sum_{\Gamma \in \mathcal{E}}\int_{\Gamma}\tilde{a}_{lm}\Delta t\frac{\lambda^{(m)}}{2}\tjump{\rho^{(m)}\mathbf{u}^{(m)}} : \tjump{\boldsymbol{\varphi}_{i}}d\Sigma \nonumber \\
    &-& \sum_{m=1}^{l-1}\sum_{\Gamma \in \mathcal{E}}\int_{\Gamma}a_{lm}\frac{\Delta t}{M^{2}}\ave{p^{(m)}}\jump{\boldsymbol{\varphi}_{i}}d\Sigma.
\end{eqnarray}
Following the discussion in \cite{orlando:2022, orlando:2025}, one can notice that a centered flux is employed for the quantities defined implicitly, while an upwind-biased flux is adopted for the quantities computed explicitly. Moreover, following \cite{abbate:2019, orlando:2025}, in order to obtain a numerical method effective for a wider range of Mach numbers, we take
\begin{equation}\label{eq:lambda_convex}
    \lambda^{(m)} = \max\left[f\left(M_{loc}^{+,(m)}\right) \left(\left|\mathbf{u}^{+,(n,m)}\right| + \frac{1}{M}c^{+,(m)}\right), f\left(M_{loc}^{-,(m)}\right)\left(\left|\mathbf{u}^{-,(m)}\right| + \frac{1}{M}c^{-,(m)}\right)\right],
\end{equation}
with $M_{loc}^{\pm,(m)} = M\frac{\left|\mathbf{u}\right|^{\pm,(m)}}{c^{\pm,(m)}}$ and $f\left(M_{loc}\right) = \min\left(1, M_{loc}\right)$. This choice corresponds to the convex combination between a centered flux and a Rusanov flux \cite{rusanov:1962} as proposed in \cite{abbate:2019}, so that, for $M_{loc} \ge 1$, we resort to a Rusanov flux, whereas for $M_{loc} \ll 1$, only the local fluid velocity is relevant for the numerical dissipation. Further considerations on the numerical flux will be discussed at the end of the Section. The numerical integration is based on the so-called over-integration or consistent integration, so as to guarantee exact integration. In particular, we employ $2r + 1$ Gauss-Legendre quadrature points along each coordinate direction \cite{orlando:2024b}. Analogously, the energy equation in \eqref{eq:stage_euler_IMEX} can be expressed as
\begin{equation}
    \mathbf{C}^{(l)}\mathbf{U}^{(l)} + \mathbf{D}^{(l)}\mathbf{P}^{(l)} = \mathbf{G}^{(l)}.
\end{equation}
For the sake of completeness, as well as to point out the contribution due to the novel strategy presented in Section \ref{ssec:num_eos} to handle a generic EOS, we report the expression of the components $\mathbf{C}^{(l)}$ and $\mathbf{D}^{(l)}$. The expression of $\mathbf{G}^{(l)}$ can be easily inferred from \eqref{eq:stage_euler_IMEX} and its definition entails that centered fluxes are employed for the quantities defined implicitly, while an upwind-biased flux is used for the quantities computed explicitly (see also \cite{orlando:2025}). Hence, we obtain
\begin{eqnarray}
    C_{ij}^{(l)} &=& \sum_{K \in \mathcal{T}_{h}}\int_{K} -a_{ll}\Delta t \hspace{0.05cm} h^{(l)}\rho^{(l)}\boldsymbol{\varphi_{j}} \cdot \grad\Psi_{i}d\Omega + \sum_{\Gamma \in \mathcal{E}}\int_{\Gamma}a_{ll}\Delta t\ave{h^{(l)}\rho^{(l)}\boldsymbol{\varphi_{j}}} \cdot \jump{\Psi_{i}}d\Sigma \\
    D_{ij}^{(l)} &=& \sum_{K \in \mathcal{T}_{h}}\int_{K}\frac{(\rho e)^{(l)}}{p^{(l)}}\Psi_{j}\Psi_{i}d\Omega 
\end{eqnarray}
Formally, we can derive
\begin{equation}\label{eq:Schur_vel}
    \mathbf{U}^{(l)} = \left(\mathbf{A}^{(l)}\right)^{-1}\left(\mathbf{F}^{(l)} - \mathbf{B}^{(l)}\mathbf{P}^{(l)}\right), 
\end{equation}
so as to obtain
\begin{equation}\label{eq:Schur}
    \mathbf{D}^{(l)}\mathbf{P}^{(l)} + \mathbf{C}^{(l)}\left(\mathbf{A}^{(l)}\right)^{-1}\left(\mathbf{F}^{(l)} - \mathbf{B}^{(l)}\mathbf{P}^{(l)}\right) = \mathbf{G}^{(l)}.
\end{equation}
The above system is then solved following the fixed point procedure described in \cite{brugnano:2008, dumbser:2016a, orlando:2022}. More specifically, setting $\mathbf{U}^{(l,0)} = \mathbf{U}^{(l-1)}$, $\mathbf{P}^{(l,0)} = \mathbf{P}^{(l-1)}$, one solves for $\tilde{k} = 0, \dots, L$
\begin{equation}
    \left(\mathbf{D}^{(l,\tilde{k})} - \mathbf{C}^{(l,\tilde{k})}\left(\mathbf{A}^{(l)}\right)^{-1}\mathbf{B}^{(l)}\right)\mathbf{P}^{(l,\tilde{k}+1)} = \mathbf{G}^{(l,\tilde{k})} - \mathbf{C}^{(l,\tilde{k})}\left(\mathbf{A}^{(l)}\right)^{-1}\mathbf{F}^{(l)}
\end{equation}
and then updates the velocity solving
\begin{equation}
    \mathbf{A}^{(l)}\mathbf{U}^{(l,\tilde{k})} = \mathbf{F}^{(l)} - \mathbf{B}^{(l)}\mathbf{P}^{(l,\tilde{k}+1)}.
\end{equation}

The algebraic formulation associated to \eqref{eq:pressure_elliptic_SI_IMEX_mod} is obtained substituting the degrees of freedom of the velocity into the algebraic formulation of the energy equation. Relations \eqref{eq:Schur_vel} and \eqref{eq:Schur} can be therefore employed to achieve this goal. For the sake of completeness, we report the new definitions of $\mathbf{C}^{(l)}$ and $\mathbf{D}^{(l)}$, while analogous modifications apply to the other variables. Hence, $\mathbf{C}^{(l)}$ and $\mathbf{D}^{(l)}$ now read as follows:
\begin{eqnarray}
    C_{ij}^{(l)} &=& \sum_{K \in \mathcal{T}_{h}}\int_{K} -a_{ll}\Delta t \hspace{0.05cm} h_{E}^{(l)}\rho_{I}^{(l)}\boldsymbol{\varphi_{j}} \cdot \grad\Psi_{i}d\Omega \nonumber \\
    &+& \sum_{\Gamma \in \mathcal{E}}\int_{\Gamma}a_{ll}\Delta t\ave{h_{E}^{(l)}\rho_{I}^{(l)}\boldsymbol{\varphi_{j}}} \cdot \grad\Psi_{i}d\Sigma + \sum_{K \in \mathcal{T}_{h}}\int_{K} \frac{M^{2}}{2}\mathbf{u}_{E}^{(l)} \cdot \rho_{I}^{(l)}\boldsymbol{\varphi_{j}}\Psi_{i}d\Omega \\
    D_{ij}^{(l)} &=& \sum_{K \in \mathcal{T}_{h}}\int_{K}\frac{(\rho e)_{E}^{(l)}}{p_{E}^{(l)}}\Psi_{j}\Psi_{i}d\Omega.
\end{eqnarray}

A matrix-free approach is employed \cite{arndt:2023}, meaning that no global sparse matrix is built and only the action of the linear operators on a vector is actually implemented. Matrices $\mathbf{A}^{(l)}$ and $\mathbf{D}^{(l)}$ are symmetric and positive definite, while the matrix $\mathbf{C}^{(l)}\left(\mathbf{A}^{(l)}\right)^{-1}\mathbf{B}^{(l)}$ is not symmetric. We point out that if one directly discretizes \eqref{eq:pressure_elliptic_IMEX} and \eqref{eq:pressure_elliptic_SI_IMEX}, as done, e.g., in \cite{boscheri:2021a, reddy:2023}, a symmetric positive definite linear system can be obtained. However, this approach implies the direct numerical solution of an elliptic equation and the discretization of a second order operator that, in the framework of a DG method, would require, e.g., the use of the Symmetric Interior Penalty method \cite{arnold:1982}. The use of a Schur complement type technique, as the one described in this work, allows one to employ only the standard numerical fluxes of hyperbolic problems (Rusanov and upwind-biased in this work), without defining and setting penalization constants typical of the aforementioned numerical strategy for elliptic equations. A comparison between the approach employed in this work and the direct solution of the Helmholtz-type equations \eqref{eq:pressure_elliptic_IMEX} and \eqref{eq:pressure_elliptic_SI_IMEX} will be matter of future work. In view of these considerations, a preconditioned conjugate gradient method with a geometric multigrid preconditioner is employed to solve the symmetric positive definite linear systems. The GMRES solver with a Jacobi preconditioner is employed for the solution of the non-symmetric linear systems. In future developments, we aim to implement and employ multigrid preconditioners also for the non-symmetric linear systems in the matrix-free framework so as to further improve the performance of the solver.

The DG method naturally allows for high-order accuracy without the use of reconstructions which involve large stencils. However, as discussed in \cite{jung:2024b}, its accuracy in the very low Mach regime depends on the numerical flux and on the shape of the elements. More specifically, a simplicial mesh is needed to establish low Mach accuracy. The lack of accuracy in the very low Mach limit can lead to a numerical scheme which is not convergent with a finite volume scheme, while an order reduction is observed in the case of the Discontinuous Galerkin method \cite{jung:2024b}. A low Mach fix for the Euler equations resolved with the finite volume method on Cartesian grids was proposed in \cite{barsukow:2021}. Moreover, it is known that an upwind scheme fails to solve very subsonic flows \cite{guillard:1999}. However, for moderate low Mach numbers, i.e. $M > 10^{-4}$, the convex combination \eqref{eq:lambda_convex} leads to a correct scaling of the pressure fluctuations \cite{abbate:2019, orlando:2025} and the low Mach number inaccuracy is counterbalanced by the high-order nature of the numerical scheme \cite{jung:2024b, orlando:2025}.

A detailed analysis of the spatial discretization and the development of possible remedies for very subsonic flows goes beyond the scope of the present work and will subject of future developments (see also Sections \ref{ssec:Taylor_Green_vortex} and \ref{ssec:traveling_vortex}). Since the main focus of this work is the comparison of the time discretization methods presented in Section \ref{ssec:IMEX} and \ref{ssec:SI_IMEX}, we believe that considering a minimum Mach number around $10^{-3}-10^{-4}$ allows us to perform this analysis in the low Mach limit and compares well with the minimum Mach number chosen to analyze asymptotic-preserving schemes in the literature, see, e.g., \cite{abbate:2019, boscheri:2020, boscheri:2025, cordier:2012, noelle:2014, thomann:2019}. 

\section{Numerical results}
\label{sec:num_res}

The numerical methods outlined in Section \ref{sec:num_disc} are now validated in a number of relevant benchmarks. The implementation is carried out in the framework of the \texttt{deal.II} library \cite{arndt:2023, bangerth:2007}, that is a C++ open-source software supporting the creation of finite element codes. Several libraries based on \texttt{deal.II} have been developed in the last years \cite{africa:2022, guermond:2022, kronbichler:2012}. All the simulations are performed in double precision. The employed time discretization schemes are reported in Appendix \ref{app:IMEX_coeffs}. Discrete parameter choices are associated to two Courant numbers, one based on the speed of sound denoted by $C$, the so-called acoustic Courant number, and one based on the local velocity of the flow, the so-called advective Courant number, denoted by $C_{u}$:
\begin{equation}
    C = \frac{1}{M}rc\frac{\Delta t}{\mathcal{H}}\sqrt{d} \qquad C_{u} = ru\frac{\Delta t}{\mathcal{H}}\sqrt{d}.
\end{equation}
Here, $\mathcal{H} = \min\left\{\text{diam}(K) | K \in \mathcal{T}_{h}\right\}$, $r$ is the polynomial degree employed for the spatial discretization, $c$ is the speed of sound, and $u$ is the magnitude of the flow velocity. Recall that $d$ denotes the space dimension. For what concerns the tests with the ideal gas law, $\gamma = 1.4$ is employed in \eqref{eq:ideal_gas}. Finally, following \cite{orlando:2022}, the fixed point loop \eqref{eq:pressure_elliptic_IMEX_fixed} is stopped at the iteration $\tilde{k}$ for which the maximum relative difference for the pressure is below a tolerance $\eta$, namely
\begin{equation}\label{eq:fixed_point_tol}
    \frac{\left\|\xi^{(\tilde{k})} - \xi^{(\tilde{k}-1)}\right\|_{\infty}}{\left\|\xi^{(\tilde{k})}\right\|_{\infty}} < \eta.
\end{equation}

\subsection{Taylor-Green vortex}
\label{ssec:Taylor_Green_vortex}

As a first benchmark to verify the scaling properties of the numerical methods with respect to the Mach number $M$, we consider the Taylor-Green vortex \cite{chorin:1968, zampa:2025}, that represents an exact steady solution of the incompressible Euler equations. The initial condition in non-dimensional variables reads as follows:
\begin{equation}
    \rho(\mathbf{x}, 0) = 1 \qquad \mathbf{u}(\mathbf{x}, 0) = \begin{pmatrix}
	\sin(x)\cos(y) \\
	-\cos(x)\sin(y)
    \end{pmatrix} \qquad p = 1 + \frac{1}{4}M^{2}\left(\cos(2x) + \cos(2y)\right).
\end{equation}
The computational domain is $\Omega = \left(0, 2\pi\right)^{2}$ endowed with periodic boundary conditions. The time step is such that the maximum advective Courant number is $C_{u} \approx 0.095$, while the maximum acoustic Courant number is $C \approx 0.11/M$. We employ $\eta = 10^{-10}$.

\begin{figure}[h!]
    \centering
    \includegraphics[width = 0.8\textwidth]{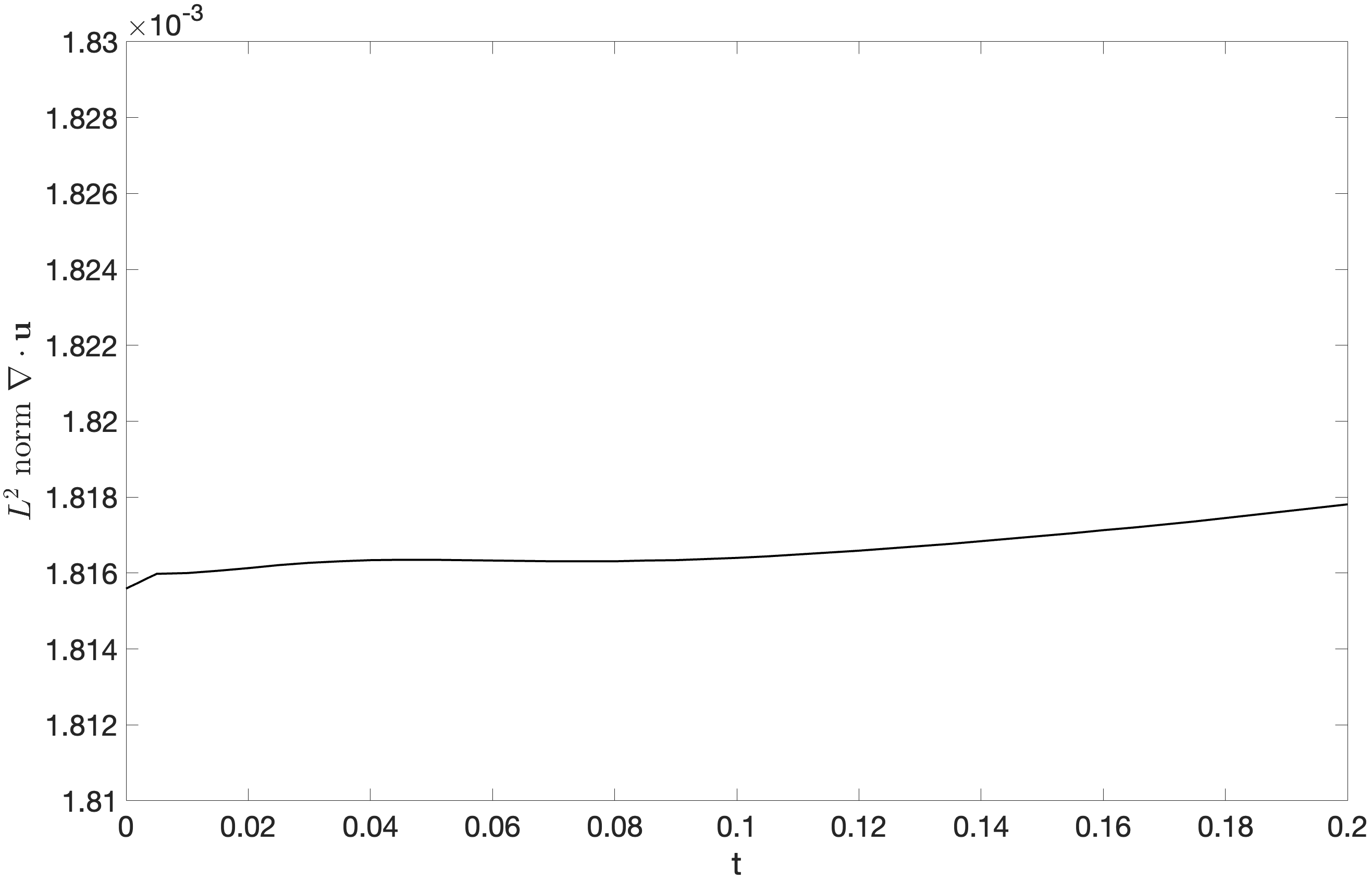}
    \caption{Taylor-Green vortex test case, time evolution of the divergence of the velocity at $M = 10^{-3}$. The results are obtained using the \textbf{IMEX} method with the IMEX-RK(3,3,3) scheme of type II in Table \ref{tab:rk3_butch_type_II} together with polynomial degree $r = 2$ and $N_{el} = 60$.}
    \label{fig:TG_vortex_Mach_div_u_IMEX_type_II}
\end{figure}

First, we consider the IMEX-RK(3,3,3) scheme of type II (Table \ref{tab:rk3_butch_type_II}) in combination with polynomial degree $r = 2$ and $N_{el} = 60$ elements along each coordinate direction. We employ the IMEX-DG method. One can easily notice that pressure fluctuations scale as $\mathcal{O}(M^{2})$ up to $M = 10^{-4}$, as expected, whereas the density fluctuations scale as $\mathcal{O}(M^{2})$ up to $M = 10^{-3}$ and a degradation is experienced at $M = 10^{-4}$. This degradation is likely related to the low Mach inaccuracy of the DG method on quadrilateral cells discussed in Section \ref{ssec:space_disc}. For what concerns the divergence of the velocity field, it scales as $\mathcal{O}\left(\mathcal{H}^{2}\right)$ (see the convergence analysis in Section \ref{ssec:traveling_vortex}), but it does not vanish as $M \to 0$ (Table \ref{tab:TG_vortex_Mach_IMEX_type_II}). As discussed in \cite{orlando:2025}, since the initial velocity field is solenoidal and the vortex is stationary, a quadratic convergence with respect to $M$ is expected for the divergence of the velocity field. This result is dependent on the spatial discretization. Indeed, since our method employs a standard nodal DG method, the divergence-free property is not imposed pointwise and the error associated to $\dive\mathbf{u}$ is basically constant in time and it is therefore related to the interpolation of the initial datum into the employed finite element space (Figure \ref{fig:TG_vortex_Mach_div_u_IMEX_type_II}). A quadratic convergence with respect to $M$ was recently obtained in \cite{zampa:2025} and preliminary results in our framework suggest that the use of Raviart-Thomas finite elements \cite{arnold:2005} for the velocity field improves the scaling properties as $M \to 0$. As already discussed at the end of Section \ref{ssec:space_disc}, a more detailed analysis of the spatial discretization is currently under investigation, while the primary goal of the present work is to perform a quantitative comparison between two different IMEX-RK approaches for the Euler equations. However, further considerations about the spatial discretization will be added at the end of the upcoming Section \ref{ssec:traveling_vortex}.

Next, we consider the SI-IMEX-DG method. A stable solution is obtained up to $M = 10^{-3}$ and the density fluctuations scale as $\mathcal{O}(M^{2})$ only up to $M = 10^{-2}$ (Table \ref{tab:TG_vortex_Mach_SI_IMEX_type_II}). This degradation is again likely mainly related to an early manifestation of the low Mach inaccuracy of the DG method on quadrilateral cells. However, we point out that some issues in the low Mach regime employing schemes of type II for the semi-implicit time discretization were already experienced in \cite{boscarino:2022}

\begin{table}[pos=H]
    \centering
    \footnotesize
    \begin{tabularx}{0.95\linewidth}{XXXXXXX}
	\toprule
        $M$ & $L^{2}$ norm $\dive\mathbf{u}$ & Rate $\dive\mathbf{u}$ & $L^{2}$ norm $\delta\rho$ & Rate $\delta\rho$ & $L^{2}$ norm $\delta p$ & Rate $\delta p$ \\
	\midrule
        $10^{-1}$ & $\num{5.76e-3}$ & & $\num{1.38e-3}$ & & $\num{1.57e-2}$ & \\
	\midrule
        $10^{-2}$ & $\num{1.82e-3}$ & $0.5$ & $\num{1.57e-5}$ & $1.9$ & $\num{1.57e-4}$ & $2.0$ \\
	\midrule
        $10^{-3}$ & $\num{1.82e-3}$ & $-$ & $\num{1.58e-7}$ & $2.0$ & $\num{1.57e-6}$ & $2.0$ \\
	\midrule
        $10^{-4}$ & $\num{1.82e-3}$ & $-$ & $\num{1.24e-8}$ & $1.1$ & $\num{1.62e-8}$ & $2.0$ \\
	\bottomrule
    \end{tabularx}
    \caption{Mach number scaling of the density and pressure fluctuations and of the divergence of the velocity field for the Taylor-Green vortex test case. Here, and in the following Tables, $\delta\rho = \rho - 1$ and $\delta p = p - 1$. The results are obtained using the \textbf{IMEX} method with the IMEX-RK(3,3,3) scheme of type II in Table \ref{tab:rk3_butch_type_II} together with polynomial degree $r = 2$ and $N_{el} = 60$.}
    \label{tab:TG_vortex_Mach_IMEX_type_II}
\end{table}

\begin{table}[pos=H]
    \centering
    \footnotesize 
    \begin{tabularx}{0.95\linewidth}{XXXXXXX}
	\toprule
        $M$ & $L^{2}$ norm $\dive\mathbf{u}$ & Rate $\dive\mathbf{u}$ & $L^{2}$ norm $\delta\rho$ & Rate $\delta\rho$ & $L^{2}$ norm $\delta p$ & Rate $\delta p$ \\
	\midrule
        $10^{-1}$ & $\num{5.76e-3}$ & & $\num{1.38e-3}$ & & $\num{1.57e-2}$ & \\
	\midrule
        $10^{-2}$ & $\num{1.82e-3}$ & $0.5$ & $\num{1.61e-5}$ & $1.9$ & $\num{1.57e-4}$ & $2.0$ \\
	\midrule
        $10^{-3}$ & $\num{1.82e-3}$ & $-$ & $\num{4.28e-6}$ & $0.6$ & $\num{1.57e-6}$ & $2.0$ \\
	\bottomrule
    \end{tabularx}
    \caption{Mach number scaling of the density and pressure fluctuations and of the divergence of the velocity field for the Taylor-Green vortex test case. The results are obtained using the \textbf{SI-IMEX} method with the IMEX-RK(3,3,3) scheme of type II in Table \ref{tab:rk3_butch_type_II} together with polynomial degree $r = 2$ and $N_{el} = 60$.}
    \label{tab:TG_vortex_Mach_SI_IMEX_type_II}
\end{table}

Next, we consider the IMEX-RK(4,3,3) scheme of type I (Table \ref{tab:rk3_butch_type_I}). For what concerns the IMEX method, analogous results with respect to the scheme of type II are obtained up to $M = 10^{-3}$ with an improvement of the scaling of the density fluctuations at $M = 10^{-4}$ (Table \ref{tab:TG_vortex_Mach_IMEX_type_I}). For what concerns the SI-IMEX method, a stable solution is established at $M = 10^{-4}$ with a stagnation of the density fluctuations (Table \ref{tab:TG_vortex_Mach_SI_IMEX_type_I}). These results confirm the superior stability of schemes of type I with respect to schemes of type II for low Mach numbers flows already experienced in \cite{boscarino:2022}. Moreover, we can infer that the low Mach accuracy is influenced also by the time discretization strategy and by the time discretization scheme.

\begin{table}[pos=H]
    \centering
    \footnotesize  
    \begin{tabularx}{0.95\linewidth}{XXXXXXX}
	\toprule
        $M$ & $L^{2}$ norm $\dive\mathbf{u}$ & Rate $\dive\mathbf{u}$ & $L^{2}$ norm $\delta\rho$ & Rate $\delta\rho$ & $L^{2}$ norm $\delta p$ & Rate $\delta p$ \\
	\midrule
        $10^{-1}$ & $\num{5.76e-3}$ & & $\num{1.38e-3}$ & & $\num{1.57e-2}$ & \\
	\midrule
        $10^{-2}$ & $\num{1.82e-3}$ & $0.5$ & $\num{1.58e-5}$ & $1.9$ & $\num{1.57e-4}$ & $2.0$ \\
	\midrule
        $10^{-3}$ & $\num{1.82e-3}$ & $-$ & $\num{1.58e-7}$ & $2.0$ & $\num{1.57e-6}$ & $2.0$ \\
	\midrule
        $10^{-4}$ & $\num{1.82e-3}$ & $-$ & $\num{4.66e-9}$ & $1.5$ & $\num{1.63e-8}$ & $2.0$ \\
	\bottomrule
    \end{tabularx}
    \caption{Mach number scaling of the density and pressure fluctuations and of the divergence of the velocity field for the Taylor-Green vortex test case. The results are obtained using the \textbf{IMEX} method with the IMEX-RK(4,3,3) scheme of type I in Table \ref{tab:rk3_butch_type_I} together with polynomial degree $r = 2$ and $N_{el} = 60$.}
    \label{tab:TG_vortex_Mach_IMEX_type_I}
\end{table}

\begin{table}[h!]
    \centering
    \footnotesize  
    \begin{tabularx}{0.95\linewidth}{XXXXXXX}
	\toprule
        $M$ & $L^{2}$ norm $\dive\mathbf{u}$ & Rate $\dive\mathbf{u}$ & $L^{2}$ norm $\delta\rho$ & Rate $\delta\rho$ & $L^{2}$ norm $\delta p$ & Rate $\delta p$ \\
	\midrule
        $10^{-1}$ & $\num{5.76e-3}$ & & $\num{1.38e-3}$ & & $\num{1.57e-2}$ & \\
	\midrule
        $10^{-2}$ & $\num{1.82e-3}$ & $0.5$ & $\num{1.58e-5}$ & $1.9$ & $\num{1.57e-4}$ & $2.0$ \\
	\midrule
        $10^{-3}$ & $\num{1.82e-3}$ & $-$ & $\num{1.69e-6}$ & $1.0$ & $\num{1.57e-6}$ & $2.0$ \\
        \midrule
        $10^{-4}$ & $\num{1.82e-3}$ & $-$ & $\num{1.68e-6}$ & $-$ & $\num{1.63e-8}$ & $2.0$ \\
	\bottomrule
    \end{tabularx}
    \caption{Mach number scaling of the density and pressure fluctuations and of the divergence of the velocity field for the Taylor-Green vortex test case. The results are obtained using the \textbf{SI-IMEX} method with the IMEX-RK(4,3,3) scheme of type I in Table \ref{tab:rk3_butch_type_I} together with polynomial degree $r = 2$ and $N_{el} = 60$.}
    \label{tab:TG_vortex_Mach_SI_IMEX_type_I}
\end{table}

\subsection{Traveling vortex at low Mach}
\label{ssec:traveling_vortex}

Next, we consider for an ideal gas a two-dimensional traveling vortex inspired by the inviscid isentropic vortex studied, e.g., in \cite{orlando:2022, tavelli:2017, zeifang:2019}. For this test, a time-dependent analytic solution is available and the convergence properties of a numerical scheme can be therefore assessed. The exact solution is indeed a propagation of the initial condition at the background velocity
$$\rho(\mathbf{x},t) = \rho(\mathbf{x} - \mathbf{u}_{\infty}t, 0) \qquad \mathbf{u}(\mathbf{x}, t) = \mathbf{u}(\mathbf{x} - \mathbf{u}_{\infty}t, 0) \qquad p(\mathbf{x},t) = p(\mathbf{x} - \mathbf{u}_{\infty}t,0).$$
Notice that a different version of this test case, which allows for a steady solution, was employed in \cite{boscheri:2025, boscheri:2021a, orlando:2025}. Following \cite{zeifang:2019}, in order to emphasize the role of the Mach number $M$, we define the perturbation as
\begin{equation}
    \delta T = \frac{1 - \gamma}{8\gamma\pi^{2}}M^{2}\beta^{2}e^{1 - \tilde{r}^2},
\end{equation}
with
$$\tilde{r}^{2} = \frac{(x - x_{0})^{2} + (y - y_{0})^{2}}{r_{0}^{2}}$$
denoting the scaled radial coordinate and $\beta$ being the vortex strength. We set $\beta = 10$, $r_{0} = 1$, and
\begin{equation}
    \rho(\mathbf{x},0) = \left(1 + \delta T\right)^{\frac{1}{\gamma-1}} \quad
    p(\mathbf{x},0) = 1 + M^{2}\left(1 + \delta T\right)^{\frac{\gamma}{\gamma-1}} = 1 + M^{2}\rho^{\gamma}.
\end{equation}
For what concerns the velocity, we define its perturbation as
\begin{equation}
    \delta\mathbf{u} = \beta M \begin{pmatrix}
	-\left(y - y_{0}\right) \\
	\left(x - x_{0}\right)
    \end{pmatrix}\frac{e^{\frac{1}{2}\left(1 - \tilde{r}^{2}\right)}}{2\pi}
\end{equation}
Finally, we set $\mathbf{u}_{\infty} = [1, 1]^{\top}$, whereas the final time is $T_{f} = 3$. To avoid problems related to the definition of boundary conditions, we choose a sufficiently large domain $\Omega = \left(-10, 10\right)^2$ and periodic boundary conditions. The behavior of the numerical methods is investigated at different Mach numbers. More specifically, we consider $M \in \left[10^{-4}, 10^{-1}\right]$. The tolerance in \eqref{eq:fixed_point_tol} is set to $\eta = 10^{-10}$.

\begin{figure}[h!]
    \begin{subfigure}{0.475\textwidth}
	\centering
        \includegraphics[width = 0.95\textwidth]{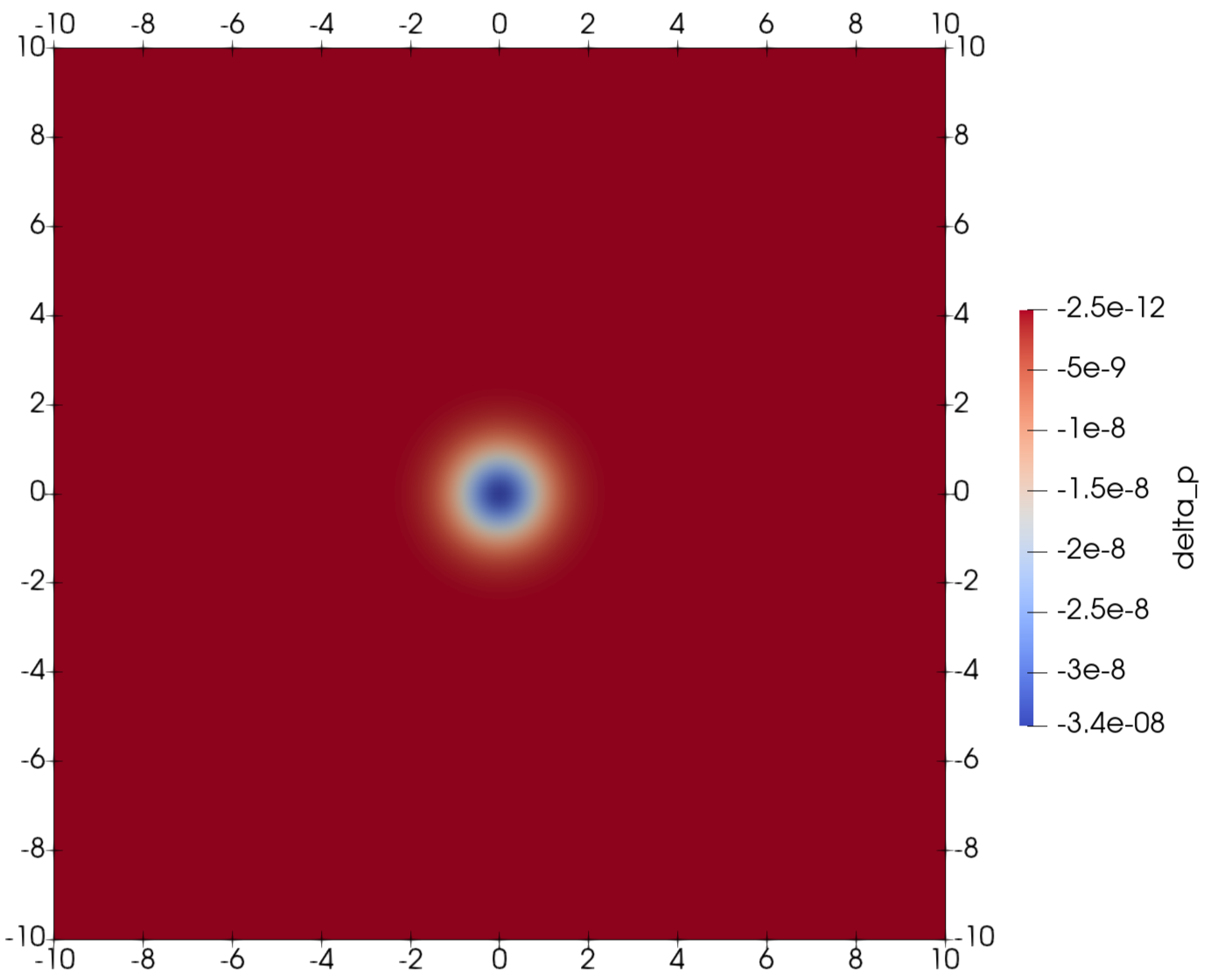}
    \end{subfigure}
    \begin{subfigure}{0.475\textwidth}
	\centering
        \includegraphics[width = 0.95\textwidth]{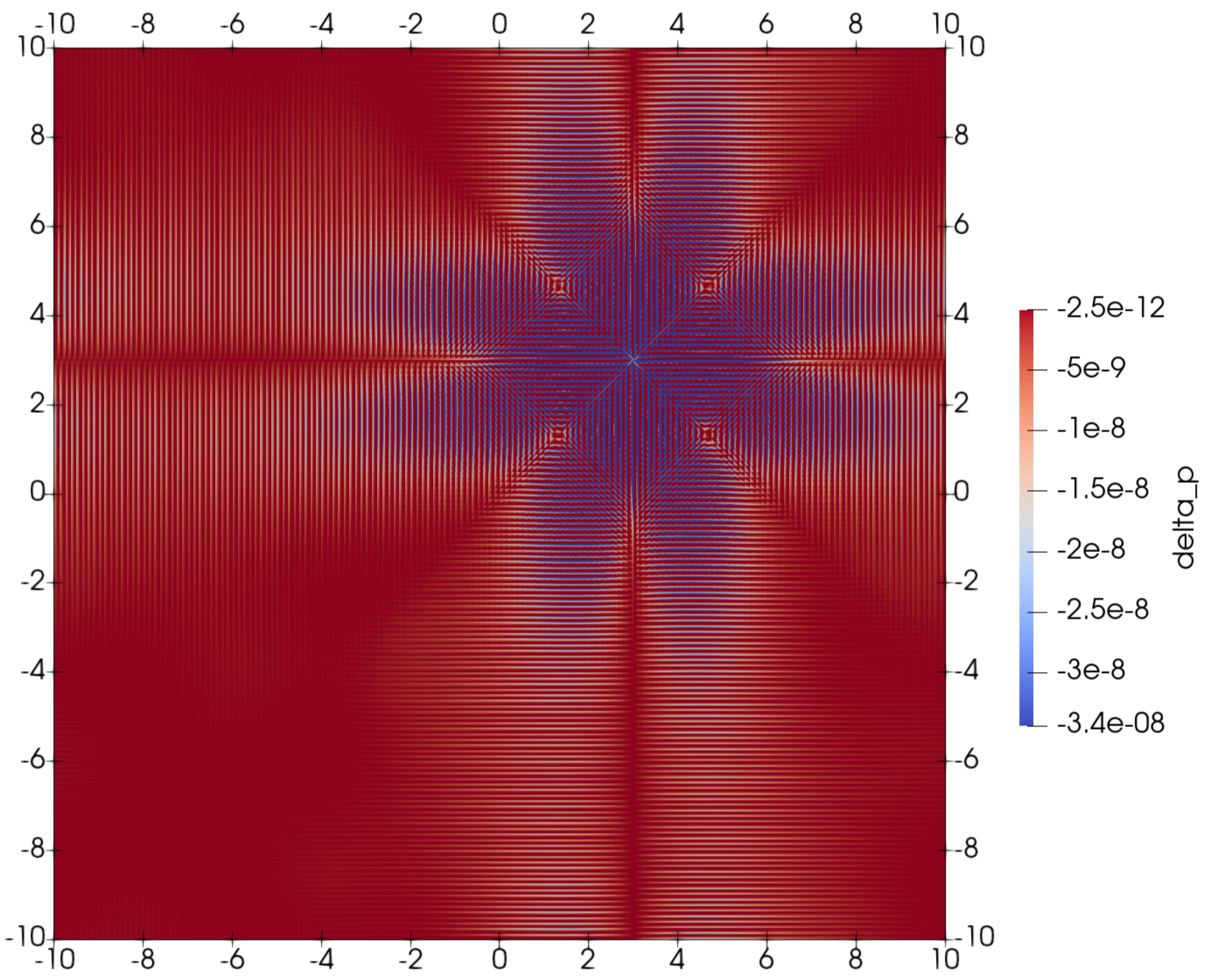}
    \end{subfigure}
    \caption{Traveling vortex test case with the IMEX-RK(2,2,2) scheme (Table \ref{tab:rk2_butch}) and polynomial degree $r = 1$, contour plots of the pressure perturbation $p - (1 - M^{2})$ at $M = 10^{-2}$ with $N_{el} = 160$. Left: initial field. Right: field at final time $t = T_{f} = 3$.}
    \label{fig:vortex_second_order_pressure}
\end{figure}

\begin{figure}[h!]
    \begin{subfigure}{0.475\textwidth}
	\centering
        \includegraphics[width = 0.95\textwidth]{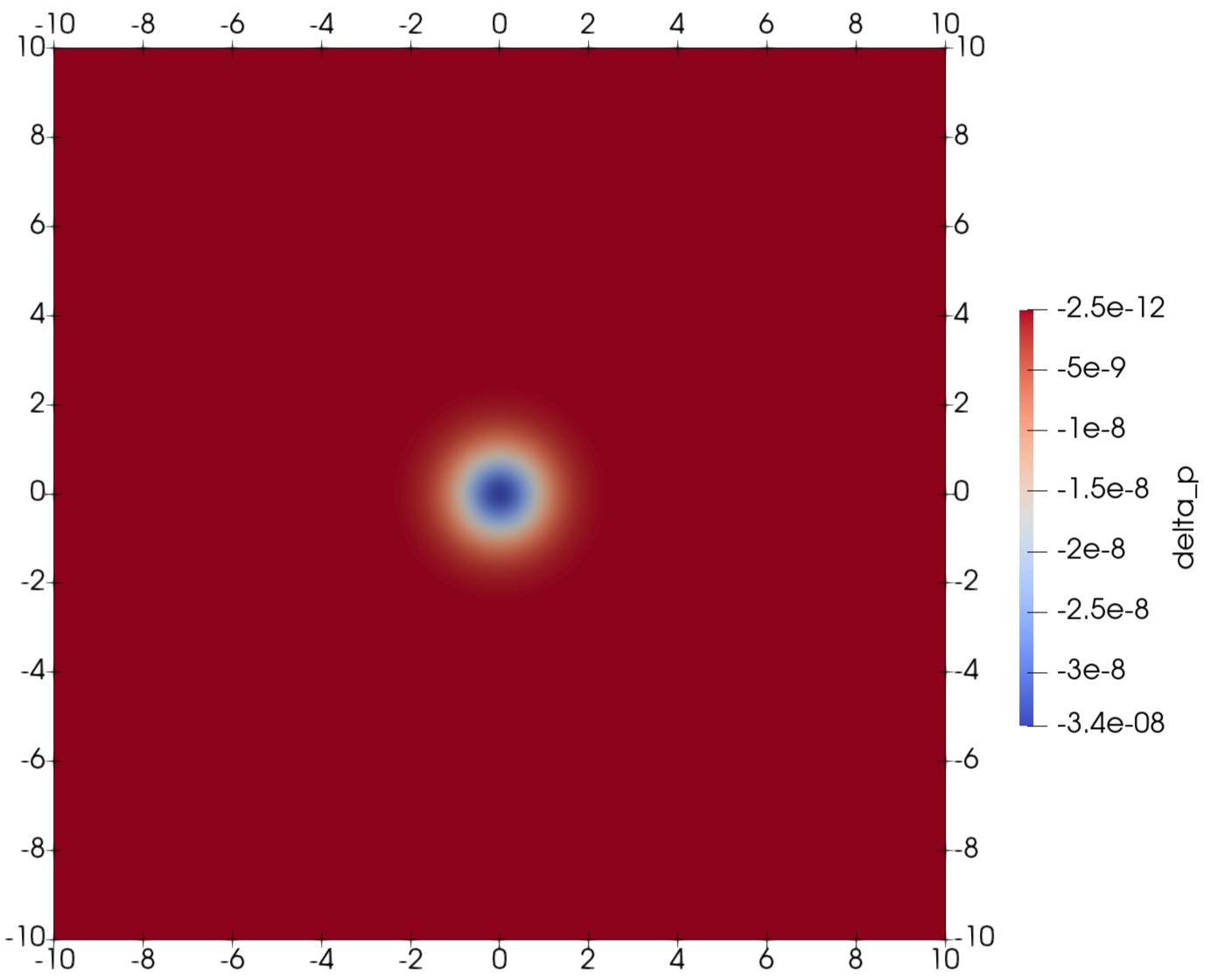}
    \end{subfigure}
    \begin{subfigure}{0.475\textwidth}
	\centering
        \includegraphics[width = 0.95\textwidth]{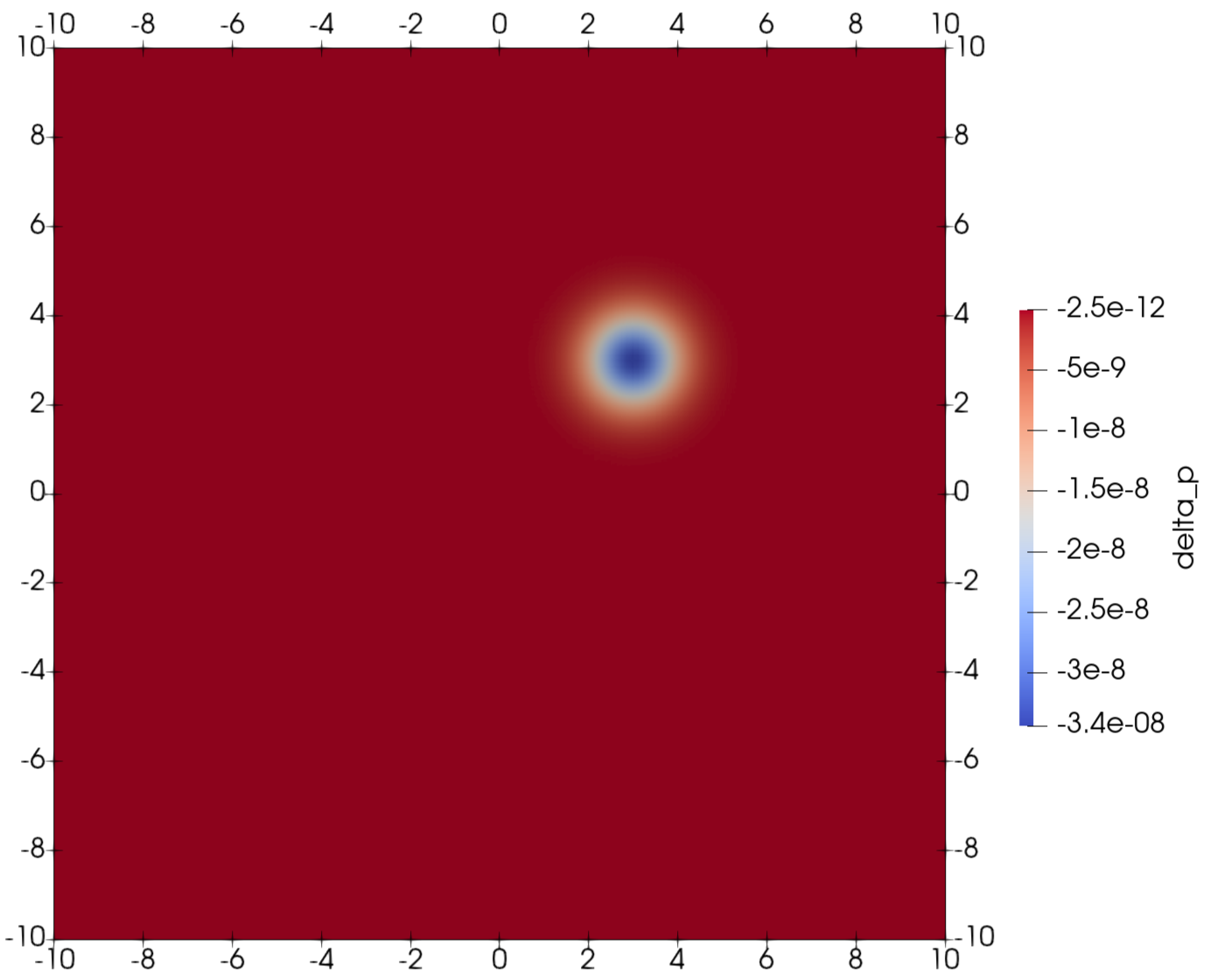}
    \end{subfigure}
    \caption{Traveling vortex test case with the IMEX-RK(2,2,2) scheme (Table \ref{tab:rk2_butch}), $r = 2$ for the velocity field and $r = 1$ for the remaining variables, contour plots of the pressure perturbation $p - (1 + M^{2})$ at $M = 10^{-2}$ with $N_{el} = 160$. Left: initial field. Right: field at final time $t = T_{f} = 3$.}
    \label{fig:vortex_second_order_pressure_corrected}
\end{figure}

First, we analyze the results obtained employing the IMEX-RK(2,2,2) scheme (Table \ref{tab:rk2_butch}) with polynomial degree $r = 1$ for the space discretization. The time step is chosen in such a way that the maximum advective Courant number is $C_{u} \approx 0.16$ and the maximum acoustic Courant number is $C \approx 0.12/M$. One can easily notice that IMEX and SI-IMEX time discretization methods yield the same level of accuracy (Tables \ref{tab:vortex_second_order_IMEX}-\ref{tab:vortex_second_order_SI_IMEX}). A degradation is experienced at $M = 10^{-2}$ and analogous results are obtained for lower values of $M$. Contour plots of the pressure perturbation show that the shape of the vortex is not preserved (Figure \ref{fig:vortex_second_order_pressure}). This is likely related to well known issues of collocated finite element type discretization on quadrilateral meshes in the low Mach regime \cite{jung:2024b, rieper:2009} (see also the discussion at the end of Section \ref{ssec:space_disc}). A simple workaround consists of increasing the polynomial degree of the finite element space employed for the discretization of the velocity field. We also refer to \cite{jung:2024a}, where a space enrichment for the velocity has been used in the framework of a finite volume scheme. One can easily notice that, if we take $r = 2$ for the finite element space of the velocity, the shape of the vortex is preserved (Figure \ref{fig:vortex_second_order_pressure_corrected}) and the expected second order convergence is established (Table \ref{tab:vortex_second_order_IMEX_M0,01}).

\begin{table}[h!]
    \centering
    \footnotesize
    \begin{tabularx}{0.5\columnwidth}{rXccXc}
	\toprule
        \multirow{1}{*}{$N_{el}$} & \multicolumn{2}{c}{$M = 10^{-1}$} & & \multicolumn{2}{c}{$M = 10^{-2}$} \\
	\cmidrule(l){2-3}\cmidrule(l){5-6}
        & $L^{2}$ relative error $\delta\mathbf{u}$ & $L^{2}$ rate $\delta\mathbf{u}$ & & $L^{2}$ relative error $\delta\mathbf{u}$ & $L^{2}$ rate $\delta\mathbf{u}$ \\
	\midrule
	$40$ & $\num{1.28e-1}$ & & & $\num{1.33e-1}$ & \\
	\midrule
        $80$ & $\num{3.99e-2}$ & $1.7$ & & $\num{5.18e-2}$ & $1.3$ \\
	\midrule
        $160$ & $\num{1.02e-2}$ & $2.1$ & & $\num{2.27e-2}$ & $\mathbf{1.2}$ \\
	\midrule
        $320$ & $\num{2.03e-3}$ & $\mathbf{2.1}$ & & \\
	\bottomrule
    \end{tabularx}
    \caption{Convergence analysis for the traveling vortex test case using the \textbf{IMEX} method with the IMEX-RK(2,2,2) scheme (Table \ref{tab:rk2_butch}) and polynomial degree $r = 1$. Here, and in the following Tables, $N_{el}$ denotes the number of elements along each direction.}
    \label{tab:vortex_second_order_IMEX}
\end{table}

\begin{table}[h!]
    \centering
    \footnotesize
    \begin{tabularx}{0.5\columnwidth}{rXccXc}
	\toprule
        \multirow{1}{*}{$N_{el}$} & \multicolumn{2}{c}{$M = 10^{-1}$} & & \multicolumn{2}{c}{$M = 10^{-2}$} \\
	\cmidrule(l){2-3}\cmidrule(l){5-6}
        & $L^{2}$ relative error $\delta\mathbf{u}$ & $L^{2}$ rate $\delta\mathbf{u}$ & & $L^{2}$ relative error $\delta\mathbf{u}$ & $L^{2}$ rate $\delta\mathbf{u}$ \\
	\midrule
	$40$ & $\num{1.28e-1}$ & & & $\num{1.32e-1}$ & \\
	\midrule
        $80$ & $\num{3.99e-2}$ & $1.7$ & & $\num{5.18e-2}$ & $1.3$ \\
	\midrule
        $160$ & $\num{1.02e-2}$ & $2.0$ & & $\num{2.27e-2}$ & $\mathbf{1.2}$ \\
	\midrule
        $320$ & $\num{2.03e-3}$ & $\mathbf{2.3}$ & & \\
	\bottomrule
    \end{tabularx}
    \caption{Convergence analysis for the traveling vortex test case using the \textbf{SI-IMEX} method with the IMEX-RK(2,2,2) scheme (Table \ref{tab:rk2_butch}) and polynomial degree $r = 1$.}
    \label{tab:vortex_second_order_SI_IMEX}
\end{table}

\begin{table}[h!]
    \centering
    \footnotesize
    \begin{tabularx}{0.35\columnwidth}{rcc}
	\toprule
	\multirow{1}{*}{$N_{el}$} & \multicolumn{2}{c}{$M = 10^{-2}$} \\
	\cmidrule(l){2-3}
        & $L^{2}$ relative error $\delta\mathbf{u}$ & $L^{2}$ rate $\delta\mathbf{u}$ \\
	\midrule
	$20$ & $\num{2.48e-2}$ & \\
	\midrule
	$40$ & $\num{3.08e-3}$ & $3.0$ \\
	\midrule
	$80$ & $\num{4.86e-4}$ & $2.7$ \\
	\midrule
	$160$ & $\num{9.30e-5}$ & $\mathbf{2.4}$ \\
	\bottomrule
    \end{tabularx}
    \caption{Convergence analysis for the traveling vortex test case using the \textbf{IMEX} method with the IMEX-RK(2,2,2) scheme (Table \ref{tab:rk2_butch}), polynomial degree $r = 2$ for the velocity and $r = 1$ for the remaining variables.}
    \label{tab:vortex_second_order_IMEX_M0,01}
\end{table}

Next, we focus on third order time discretization schemes (Tables \ref{tab:rk3_butch_type_I}-\ref{tab:rk3_butch_type_II}) in combination with polynomial degree $r = 2$. The time step is such that the maximum advective Courant number is $C_{u} \approx 0.094$, while the maximum acoustic Courant number is $C \approx 0.07/M$. All the schemes provide a similar level of accuracy and the expected third order convergence rate is established (Tables \ref{tab:vortex_third_order_IMEX_type_II}-\ref{tab:vortex_third_order_SI_IMEX_type_I}) up to $M = 10^{-3}$, except for the SI-IMEX method in combination with the scheme of type II at $M = 10^{-3}$. As already remarked in Section \ref{ssec:Taylor_Green_vortex}, some issues in the low Mach regime employing schemes of type II for the semi-implicit time discretization were already experienced in \cite{boscarino:2022}. Moreover, for $M = 10^{-4}$, a stable solution for the SI-IMEX method is obtained only employing the scheme of type I. This further confirms the superior stability of schemes of type I with respect to schemes of type II for low Mach numbers when using the SI-IMEX method. 

\begin{figure}[h!]
    \begin{subfigure}{0.475\textwidth}
	\centering
        \includegraphics[width = 0.95\textwidth]{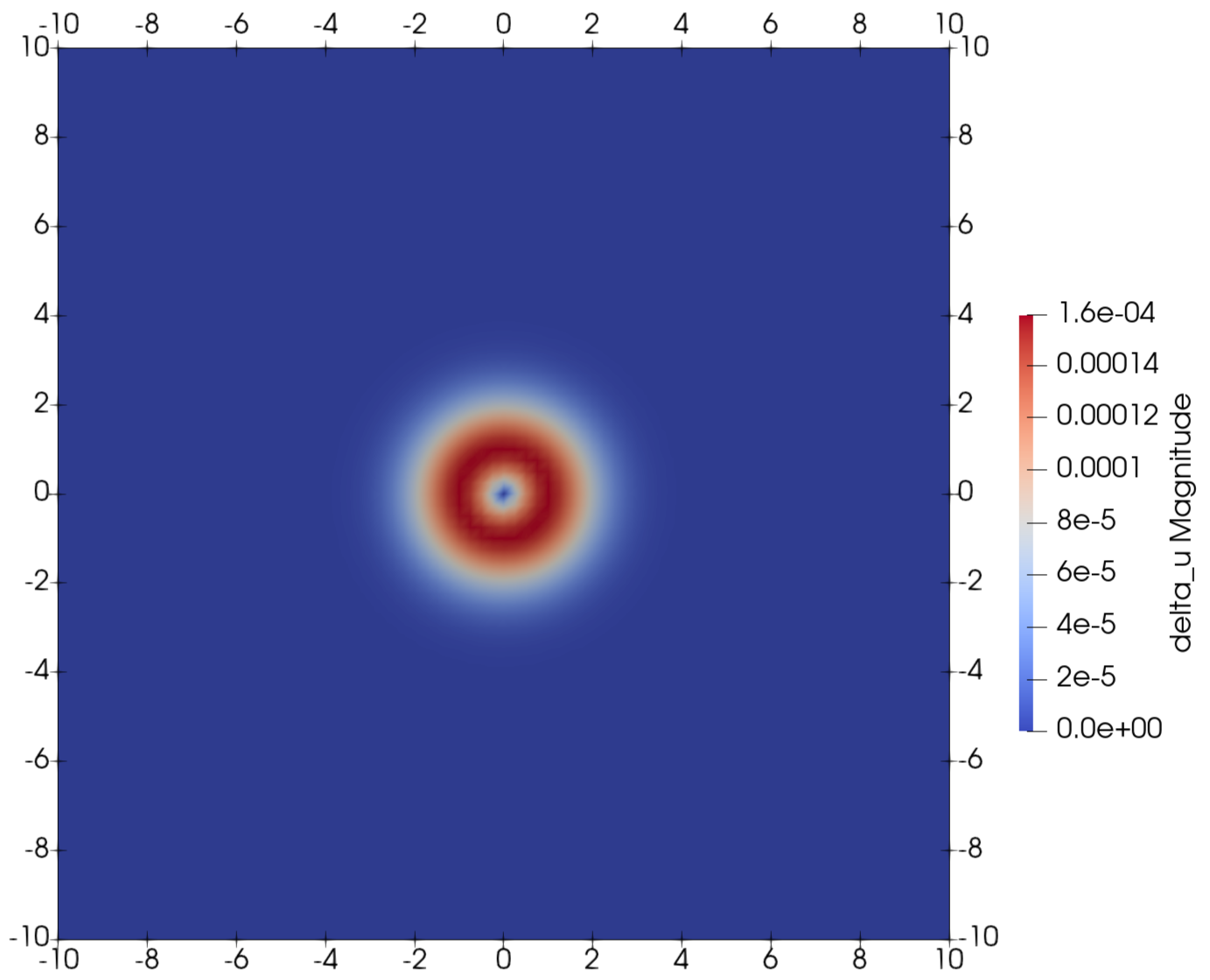}
    \end{subfigure}
    \begin{subfigure}{0.475\textwidth}
	\centering
        \includegraphics[width = 0.95\textwidth]{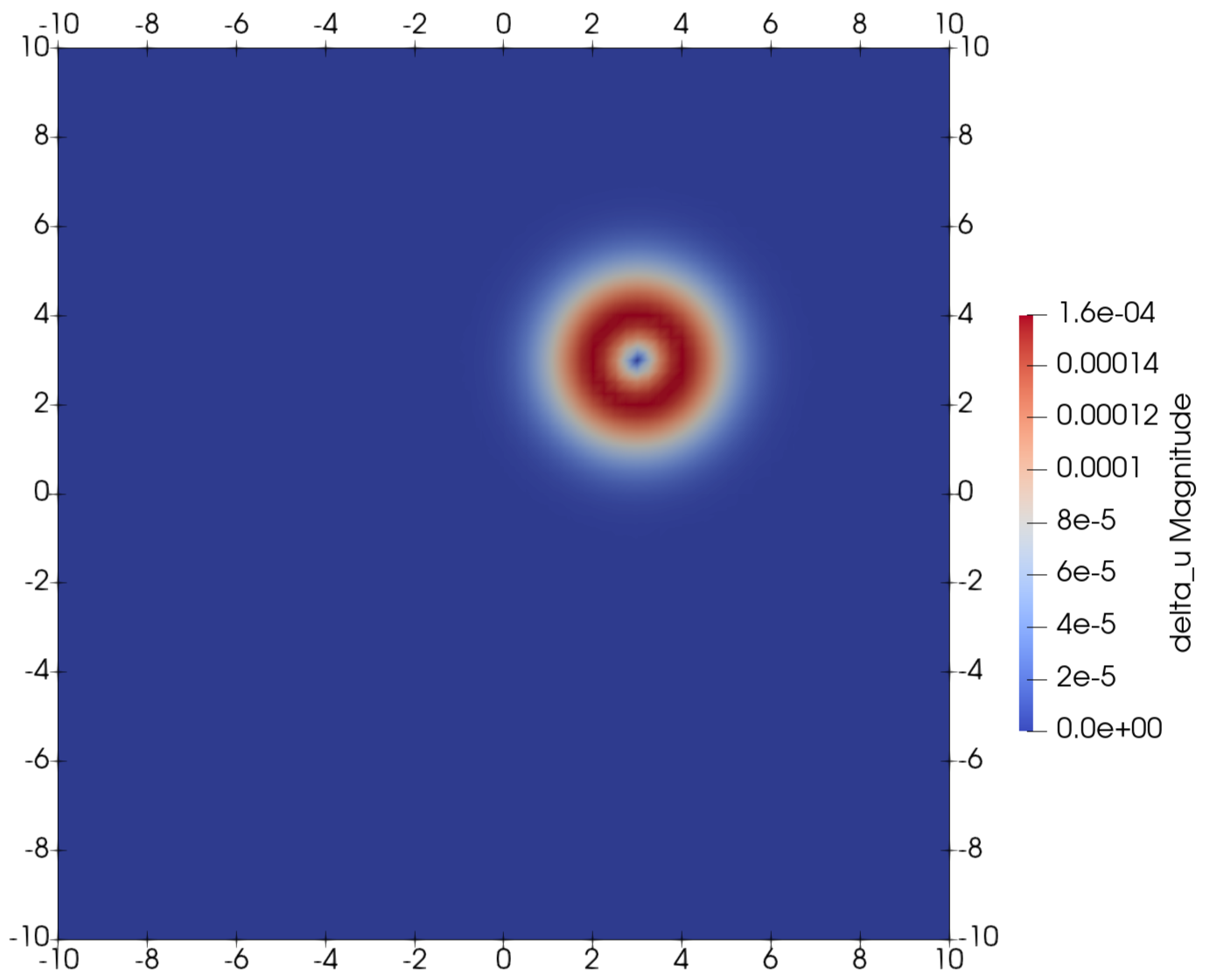}
    \end{subfigure}
    \caption{Traveling vortex test case, contour plot of the velocity perturbation at $M = 10^{-4}$. Left: initial field. Right: field at final time $t = T_{f} = 3$. The results are obtained using the IMEX method with the IMEX-RK(3,3,3) scheme (Table \ref{tab:rk3_butch_type_II}) and polynomial degree $r = 2$ with $N_{el} = 40$.}
    \label{fig:vortex_third_order_delta_u}
\end{figure}

A saturation of the error is experienced at $M = 10^{-4}$. Analogous results with a saturation of the error at $M = 10^{-4}$ were obtained in \cite{zeifang:2019}. Since in most low Mach number applications the main interest lies on the velocity field rather than on the acoustics, we notice, similarly to \cite{zeifang:2019}, that no visible spurious effect arises in the velocity field in spite of the order reduction (Figure \ref{fig:vortex_third_order_delta_u}). As already discussed in Section \ref{ssec:space_disc}, the DG method on quadrilateral cells becomes low Mach inaccurate for $M < 10^{-4}$ and therefore this order reduction is likely related to a manifestation of this phenomenon. Moreover, round-off errors are dominant in this configuration \cite{zeifang:2019} and the use of quadruple precision is crucial to reach really very low Mach numbers \cite{boscheri:2021b}. However, it is worth to notice that in the case of the steady vortex, a correct scaling was established at $M = 10^{-4}$ using the third order method of type II in \cite{orlando:2025}. Hence, the saturation of the error could be also related to the well-known order reduction phenomenon experienced for stiff ODE's problems \cite{wanner:1996}. In future work, as already mentioned, we aim to analyze the behavior employing a spatial discretization based on compatible finite elements, as recently done in \cite{zampa:2025}, so as to guarantee that the initial condition of the velocity field is setup as the discrete derivative of a vector potential, and to use simplicial meshes \cite{jung:2024b, zampa:2025} or Voronoi meshes \cite{boscheri:2025}, which have been shown to be low Mach accurate for a steady vortex.

\begin{table}[h!]
    \centering
    \footnotesize
    \begin{tabularx}{0.97\columnwidth}{rXXcXXcXXcXX}
	\toprule
        \multirow{1}{*}{$N_{el}$} & \multicolumn{2}{c}{$M = 10^{-1}$} & & \multicolumn{2}{c}{$M = 10^{-2}$} & & \multicolumn{2}{c}{$M = 10^{-3}$} & & \multicolumn{2}{c}{$M = 10^{-4}$} \\
        \cmidrule(l){2-3}\cmidrule(l){5-6}\cmidrule(l){8-9}\cmidrule(l){11-12}
        & $L^{2}$ relative error $\delta\mathbf{u}$ & $L^{2}$ rate $\delta\mathbf{u}$ & & $L^{2}$ relative error $\delta\mathbf{u}$ & $L^{2}$ rate $\delta\mathbf{u}$ & & $L^{2}$ relative error $\delta\mathbf{u}$ & $L^{2}$ rate $\delta\mathbf{u}$ & & $L^{2}$ relative error $\delta\mathbf{u}$ & $L^{2}$ rate $\delta\mathbf{u}$ \\
	\midrule
        $20$ & $\num{3.12e-2}$ & & & $\num{3.55e-2}$ & & & $\num{3.55e-2}$ & & & $\num{3.59e-2}$ & \\
	\midrule
        $40$ & $\num{2.40e-3}$ & $3.7$ & & $\num{2.32e-3}$ & $3.9$ & & $\num{2.32e-3}$ & $3.9$ & & $\num{5.51e-3}$ & $2.7$ \\
	\midrule
        $80$ & $\num{3.17e-4}$ & $2.9$ & & $\num{2.92e-4}$ & $3.0$ & & $\num{2.90e-4}$ & $3.0$ & & $\num{5.71e-3}$ & $-$ \\
	\midrule
        $160$ & $\num{4.13e-5}$ & $\mathbf{2.9}$ & & $\num{3.73e-5}$ & $\mathbf{3.0}$ & & $\num{4.61e-5}$ & $\mathbf{2.7}$ & & & \\
	\bottomrule
    \end{tabularx}
    \caption{Convergence analysis for the traveling vortex test case using the \textbf{IMEX} method with the IMEX-RK(3,3,3) scheme of type II (Table \ref{tab:rk3_butch_type_II}) and polynomial degree $r = 2$.}
    \label{tab:vortex_third_order_IMEX_type_II}
\end{table}

\begin{table}[h!]
    \centering
    \footnotesize
    \begin{tabularx}{0.8\columnwidth}{rXXcXXcXX}
	\toprule
        \multirow{1}{*}{$N_{el}$} & \multicolumn{2}{c}{$M = 10^{-1}$} & & \multicolumn{2}{c}{$M = 10^{-2}$} & & \multicolumn{2}{c}{$M = 10^{-3}$} \\
        \cmidrule(l){2-3}\cmidrule(l){5-6}\cmidrule(l){8-9}
        & $L^{2}$ relative error $\delta\mathbf{u}$ & $L^{2}$ rate $\delta\mathbf{u}$ & & $L^{2}$ relative error $\delta\mathbf{u}$ & $L^{2}$ rate $\delta\mathbf{u}$ & & $L^{2}$ relative error $\delta\mathbf{u}$ & $L^{2}$ rate $\delta\mathbf{u}$ \\
	\midrule
        $20$ & $\num{3.12e-2}$ & & & $\num{3.56e-2}$ & & & $\num{3.57e-2}$ \\
	\midrule
        $40$ & $\num{2.40e-3}$ & $3.7$ & & $\num{2.33e-3}$ & $3.9$ & & $\num{2.33e-3}$ & $3.9$ \\
	\midrule
        $80$ & $\num{3.17e-4}$ & $2.9$ & & $\num{2.92e-4}$ & $3.0$ & & $\num{3.85e-4}$ & $2.6$ \\
	\midrule
        $160$ & $\num{4.13e-5}$ & $\mathbf{2.9}$ & & $\num{3.73e-5}$ & $\mathbf{3.0}$ & & $\num{2.42e-4}$ & $0.7$ \\
	\bottomrule
    \end{tabularx}
    \caption{Convergence analysis for the traveling vortex test case using the \textbf{SI-IMEX} method with the IMEX-RK(3,3,3) scheme of type II (Table \ref{tab:rk3_butch_type_II}) and polynomial degree $r = 2$.}
    \label{tab:vortex_third_order_SI_IMEX_type_II}
\end{table}

\begin{table}[h!]
    \centering
    \footnotesize
    \begin{tabularx}{0.97\columnwidth}{rXXcXXcXXcXX}
	\toprule
        \multirow{1}{*}{$N_{el}$} & \multicolumn{2}{c}{$M = 10^{-1}$} & & \multicolumn{2}{c}{$M = 10^{-2}$} & & \multicolumn{2}{c}{$M = 10^{-3}$} & & \multicolumn{2}{c}{$M = 10^{-4}$} \\
        \cmidrule(l){2-3}\cmidrule(l){5-6}\cmidrule(l){8-9}\cmidrule(l){11-12}
        & $L^{2}$ relative error $\delta\mathbf{u}$ & $L^{2}$ rate $\delta\mathbf{u}$ & & $L^{2}$ relative error $\mathbf{u}$ & $L^{2}$ rate $\delta\mathbf{u}$ & & $L^{2}$ relative error $\delta\mathbf{u}$ & $L^{2}$ rate $\delta\mathbf{u}$ & & $L^{2}$ relative error $\delta\mathbf{u}$ & $L^{2}$ rate $\delta\mathbf{u}$ \\
	\midrule
        $20$ & $\num{3.12e-2}$ & & & $\num{3.55e-2}$ & & & $\num{3.55e-2}$ & & & $\num{3.63e-2}$ & \\
	\midrule
        $40$ & $\num{2.40e-3}$ & $3.7$ & & $\num{2.32e-3}$ & $3.9$ & & $\num{2.31e-3}$ & $3.9$ & & $\num{3.57e-3}$ & $3.3$ \\
	\midrule
        $80$ & $\num{3.17e-4}$ & $2.9$ & & $\num{2.92e-4}$ & $3.0$ & & $\num{2.90e-4}$ & $3.0$ & & $\num{4.50e-3}$ & $-$ \\
	\midrule
        $160$ & $\num{4.13e-5}$ & $\mathbf{2.9}$ & & $\num{3.73e-5}$ & $\mathbf{3.0}$ & & $\num{4.61e-5}$ & $\mathbf{2.7}$ & & & \\
	\bottomrule
    \end{tabularx}
    \caption{Convergence analysis for the traveling vortex test case using the \textbf{IMEX} method with the IMEX-RK(4,3,3) scheme of type I (Table \ref{tab:rk3_butch_type_I}) and polynomial degree $r = 2$.}
    \label{tab:vortex_third_order_IMEX_type_I}
\end{table}

\begin{table}[h!]
    \centering
    \footnotesize
    \begin{tabularx}{0.97\columnwidth}{rXXcXXcXXcXX}
	\toprule
        \multirow{1}{*}{$N_{el}$} & \multicolumn{2}{c}{$M = 10^{-1}$} & & \multicolumn{2}{c}{$M = 10^{-2}$} & & \multicolumn{2}{c}{$M = 10^{-3}$} & & \multicolumn{2}{c}{$M = 10^{-4}$} \\
        \cmidrule(l){2-3}\cmidrule(l){5-6}\cmidrule(l){8-9}\cmidrule(l){11-12}
        & $L^{2}$ relative error $\delta\mathbf{u}$ & $L^{2}$ rate $\delta\mathbf{u}$ & & $L^{2}$ relative error $\delta\mathbf{u}$ & $L^{2}$ rate $\delta\mathbf{u}$ & & $L^{2}$ relative error $\delta\mathbf{u}$ & $L^{2}$ rate $\delta\mathbf{u}$ & & $L^{2}$ relative error $\delta\mathbf{u}$ & $L^{2}$ rate $\delta\mathbf{u}$ \\
	\midrule
        $20$ & $\num{3.12e-2}$ & & & $\num{3.55e-2}$ & & & $\num{3.56e-2}$ & & & $\num{4.70e-2}$ & \\
	\midrule
        $40$ & $\num{2.40e-3}$ & $3.7$ & & $\num{2.32e-3}$ & $3.9$ & & $\num{2.32e-3}$ & $3.9$ & & $\num{2.12e-2}$ & $1.1$ \\
	\midrule
        $80$ & $\num{3.17e-4}$ & $2.9$ & & $\num{2.92e-4}$ & $3.0$ & & $\num{2.90e-4}$ & $3.0$ & & $\num{1.01e-2}$ & $1.1$ \\
	\midrule
        $160$ & $\num{4.13e-5}$ & $\mathbf{2.9}$ & & $\num{3.73e-5}$ & $\mathbf{2.9}$ & & $\num{4.69e-5}$ & $\mathbf{2.7}$ & & & \\
	\bottomrule
    \end{tabularx}
    \caption{Convergence analysis for the traveling vortex test case using the \textbf{SI-IMEX} method with IMEX-RK(4,3,3) scheme of type I (Table \ref{tab:rk3_butch_type_I}) and polynomial degree $r = 2$}
    \label{tab:vortex_third_order_SI_IMEX_type_I}
\end{table}

\begin{figure}[h!]
    \centering
    \begin{subfigure}{0.475\textwidth}
	\centering
        \includegraphics[width = 0.95\textwidth]{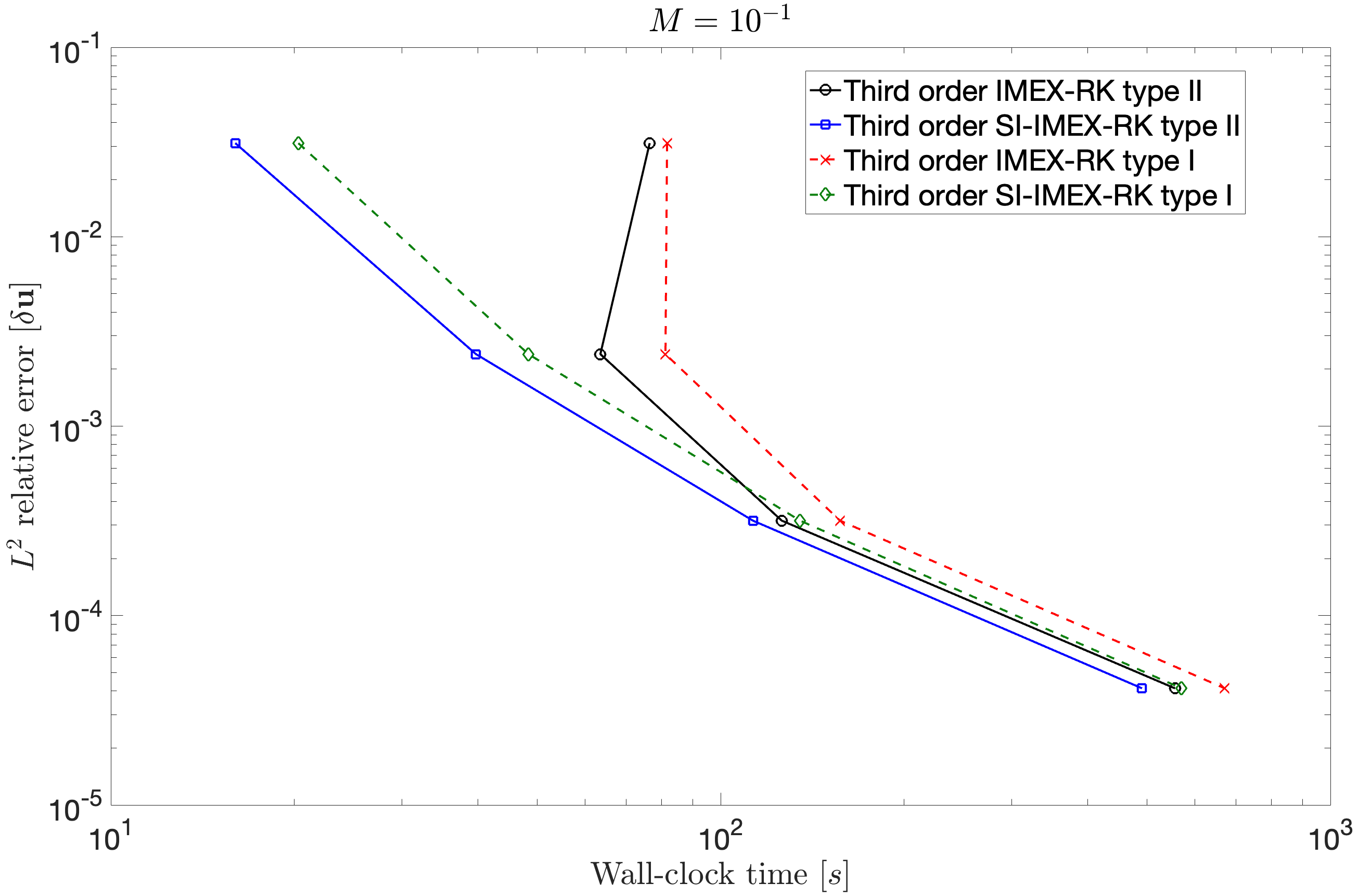}
    \end{subfigure}
    \begin{subfigure}{0.475\textwidth}
	\centering
        \includegraphics[width = 0.95\textwidth]{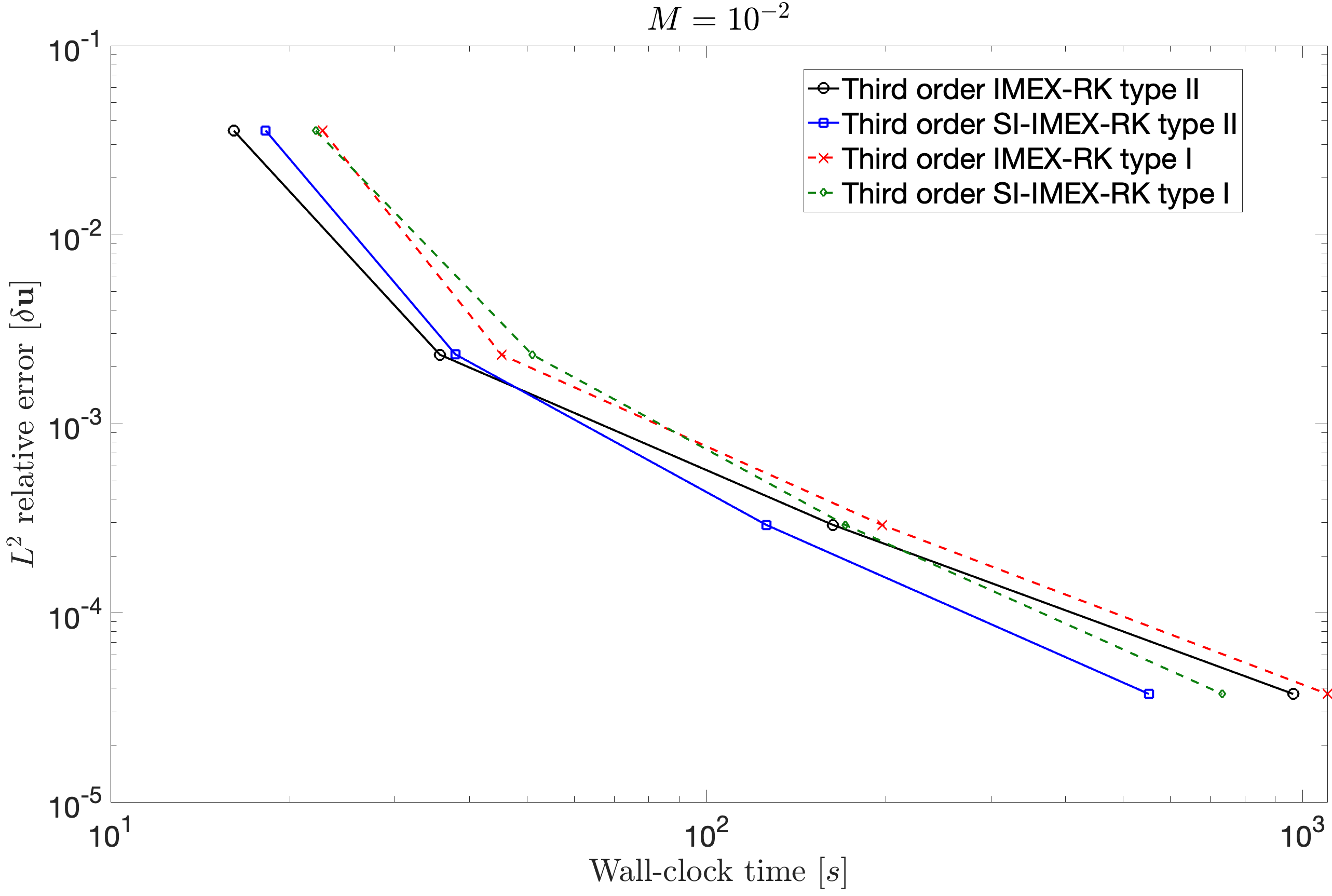}
    \end{subfigure}
    \begin{subfigure}{0.475\textwidth}
	\centering
        \includegraphics[width = 0.95\textwidth]{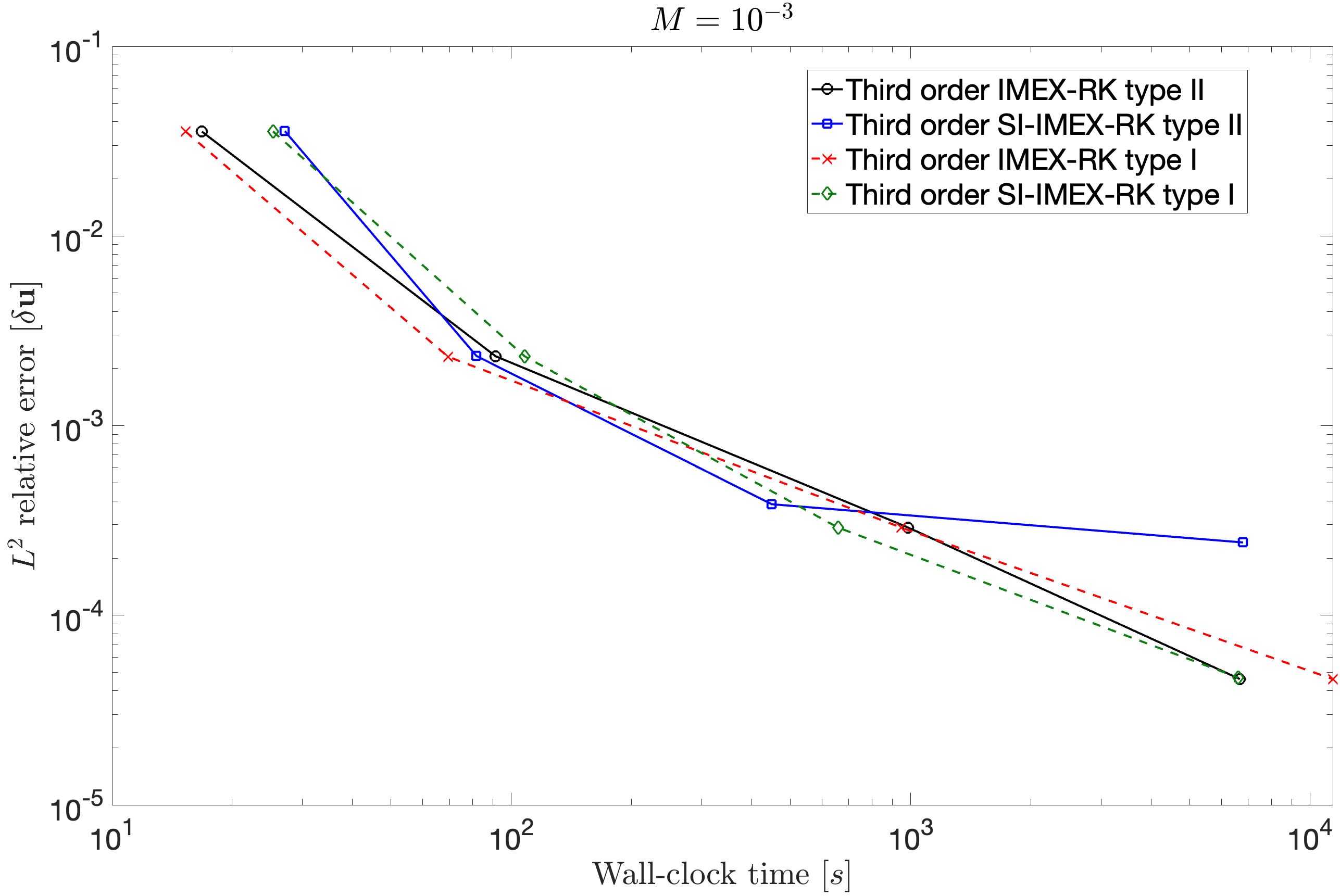}
    \end{subfigure}
    \caption{Traveling vortex test case, comparison of wall-clock times employing third order time discretization schemes and polynomial degree $r = 2$. Top-left: $M = 10^{-1}$. Top-right: $M = 10^{-2}$. Bottom: $M = 10^{-3}$. An advective Courant number $C_{u} \approx 0.19$ is employed for the IMEX method. The solid black lines denote the results obtained with the IMEX method with the scheme of type II (Table \ref{tab:rk3_butch_type_II}), the solid blue lines show the results with the SI-IMEX method with the scheme of type II (Table \ref{tab:rk3_butch_type_II}), the dashed red lines report the results achieved with the IMEX method with the scheme of type I (Table \ref{tab:rk3_butch_type_I}), whereas the dashed green lines represent the results established with the SI-IMEX method with the scheme of type I (Table \ref{tab:rk3_butch_type_I}).}
    \label{fig:vortex_third_order_WT}
\end{figure}

The IMEX method described in Section \ref{ssec:IMEX} allows the use of a larger time step ensuring the same level of accuracy, because of the superior stability provided by the fixed point loop \cite{dumbser:2016a}. In particular, the time step can be doubled, yielding a maximum advective Courant number $C_{u} \approx 0.19$. In spite of the use of a smaller time step, the SI-IMEX method and, more specifically, the SI-IMEX time discretization using the scheme of type II, provides in general better computational performance (Figure \ref{fig:vortex_third_order_WT}). As already discussed, the scheme of type I is more robust for $M = 10^{-3}$ and it becomes also more efficient as the spatial resolution increases. In the low Mach regime, differences in terms of computational cost between the IMEX method and the SI-IMEX method are significantly smaller: the efficiency gain reduces to around $5\%$ for $M = 10^{-3}$, meaning that the two methods are essentially equivalent in this respect. We will further discuss this point in Section \ref{ssec:fixed_point}.

Next, we employ the fourth order time discretization schemes of type ARS (Table \ref{tab:rk4_butch_ARS}) and of type II (Table \ref{tab:rk4_butch_type_II}) in combination with polynomial degree $r = 3$. The time step is such that the maximum advective Courant number is $C_{u} \approx 0.05$ and the maximum acoustic Courant number is $C \approx 0.035/M$. The SI-IMEX method requires a smaller time step for $M \le 10^{-2}$. More specifically, an advective Courant number $C_{u} \approx 0.025$ is required to achieve a stable solution. One can easily notice that the expected convergence rates are established up to $M = 10^{-3}$ (Tables \ref{tab:vortex_fourth_order_IMEX_ARS}-\ref{tab:vortex_fourth_order_SI_IMEX_ARS}), except for the IMEX method at $M = 10^{-2}$, for which an order reduction is experienced. Analogous considerations to those reported for the third order discretization schemes for $M \le 10^{-4}$ are valid. In particular, the correct scaling at $M = 10^{-4}$ for the steady vortex using the IMEX method with the scheme of type ARS was obtained in \cite{orlando:2025}, so that the experienced order reduction is likely dependent on both the low Mach inaccuracy of the DG method on quadrilateral cells and the order reduction phenomenon typical of stiff problems \cite{wanner:1996}.

For what concerns the time discretization method of type II (Table \ref{tab:rk4_butch_type_II}), the expected convergence rates are established for the IMEX method and the accuracy is preserved up to $M = 10^{-3}$ (Table \ref{tab:vortex_fourth_order_IMEX_type_II}). On the contrary, the SI-IMEX method shows the expected convergence rates for $M \ge 10^{-2}$, whereas severe issues start from $M = 10^{-3}$ (Table \ref{tab:vortex_fourth_order_SI_IMEX_type_II}), unlike the scheme of type ARS. As already mentioned, schemes of type II show some issues in the low Mach regime for the semi-implicit time discretization \cite{boscarino:2022}. Moreover, following also the results in Section \ref{ssec:Taylor_Green_vortex}, we can further infer that the activation of spurious modes (hence the low Mach accuracy) is influenced by the time discretization strategy and by the time discretization scheme. 

\begin{table}[h!]
    \centering
    \footnotesize
    \begin{tabularx}{0.97\columnwidth}{rXXcXXcXXcXX}
	\toprule
        \multirow{1}{*}{$N_{el}$} & \multicolumn{2}{c}{$M = 10^{-1}$} & & \multicolumn{2}{c}{$M = 10^{-2}$} & & \multicolumn{2}{c}{$M = 10^{-3}$} & & \multicolumn{2}{c}{$M = 10^{-4}$} \\
        \cmidrule(l){2-3}\cmidrule(l){5-6}\cmidrule(l){8-9}\cmidrule(l){11-12}
        & $L^{2}$ relative error $\delta\mathbf{u}$ & $L^{2}$ rate $\delta\mathbf{u}$ & & $L^{2}$ relative error $\delta\mathbf{u}$ & $L^{2}$ rate $\delta\mathbf{u}$ & & $L^{2}$ relative error $\delta\mathbf{u}$ & $L^{2}$ rate $\delta\mathbf{u}$ & & $L^{2}$ relative error $\delta\mathbf{u}$ & $L^{2}$ rate $\delta\mathbf{u}$ \\
	\midrule
        $10$ & $\num{4.85e-2}$ & & & $\num{4.88e-2}$ & & & $\num{4.88e-2}$ & & & $\num{9.15e-2}$ & \\
	\midrule
        $20$ & $\num{2.88e-3}$ & $4.1$ & & $\num{2.73e-3}$ & $4.2$ & & $\num{2.71e-3}$ & $4.2$ & & $\num{6.37e-2}$ & $0.5$ \\
	\midrule
        $40$ & $\num{1.23e-4}$ & $4.5$ & & $\num{1.30e-4}$ & $\mathbf{4.4}$ & & $\num{1.31e-4}$ & $4.4$ & & $\num{2.23e-2}$ & $1.5$ \\
	\midrule
        $80$ & $\num{8.47e-6}$ & $\mathbf{3.9}$ & & $\num{6.41e-5}$ & $1.0$ & & $\num{1.22e-5}$ & $\mathbf{3.4}$ & & & \\
        \bottomrule
    \end{tabularx}
    \caption{Convergence analysis for the traveling vortex test case using the \textbf{IMEX} method with the IMEX-RK(5,5,4) scheme of type ARS (Table \ref{tab:rk4_butch_ARS}) and polynomial degree $r = 3$.}
    \label{tab:vortex_fourth_order_IMEX_ARS}
\end{table}

\begin{table}[h!]
    \centering
    \footnotesize
    \begin{tabularx}{0.97\columnwidth}{rXXcXXcXXcXX}
	\toprule
        \multirow{1}{*}{$N_{el}$} & \multicolumn{2}{c}{$M = 10^{-1}$} & & \multicolumn{2}{c}{$M = 10^{-2}$} & & \multicolumn{2}{c}{$M = 10^{-3}$} & & \multicolumn{2}{c}{$M = 10^{-4}$} \\
        \cmidrule(l){2-3}\cmidrule(l){5-6}\cmidrule(l){8-9}\cmidrule(l){11-12}
        & $L^{2}$ relative error $\delta\mathbf{u}$ & $L^{2}$ rate $\mathbf{u}$ & & $L^{2}$ relative error $\delta\mathbf{u}$ & $L^{2}$ rate $\delta\mathbf{u}$ & & $L^{2}$ relative error $\delta\mathbf{u}$ & $L^{2}$ rate $\delta\mathbf{u}$ & & $L^{2}$ relative error $\delta\mathbf{u}$ & $L^{2}$ rate $\delta\mathbf{u}$ \\
	\midrule
        $10$ & $\num{4.85e-2}$ & & & $\num{4.89e-2}$ & & & $\num{4.88e-2}$ & & & $\num{7.38e-2}$ & \\
	\midrule
        $20$ & $\num{2.88e-3}$ & $4.1$ & & $\num{2.74e-3}$ & $4.2$ & & $\num{2.71e-3}$ & $4.2$ & & $\num{6.52e-2}$ & $0.2$ \\
	\midrule
        $40$ & $\num{1.23e-4}$ & $4.5$ & & $\num{1.30e-4}$ & $4.4$ & & $\num{1.31e-4}$ & $4.4$ & & $\num{6.27e-2}$ & $-$ \\
	\midrule
        $80$ & $\num{8.47e-6}$ & $\mathbf{3.9}$ & & $\num{1.10e-5}$ & $\mathbf{3.6}$ & & $\num{1.24e-5}$ & $\mathbf{3.4}$ & & & \\
	\bottomrule
    \end{tabularx}
    \caption{Convergence analysis for the traveling vortex test case using the \textbf{SI-IMEX} method with the IMEX-RK(5,5,4) scheme of type ARS (Table \ref{tab:rk4_butch_ARS}) and polynomial degree $r = 3$. Note that for $M = 10^{-1}$ the advective Courant number is $C_{u} \approx 0.05$, while for $M \le 10^{-2}$ the advective Courant number is $C_{u} \approx 0.025$ to achieve a stable solution.}
    \label{tab:vortex_fourth_order_SI_IMEX_ARS}
\end{table}

\begin{table}[h!]
    \centering
    \footnotesize
    \begin{tabularx}{0.97\columnwidth}{rXXcXXcXXcXX}
        \toprule
        \multirow{1}{*}{$N_{el}$} & \multicolumn{2}{c}{$M = 10^{-1}$} & & \multicolumn{2}{c}{$M = 10^{-2}$} & & \multicolumn{2}{c}{$M = 10^{-3}$} & & \multicolumn{2}{c}{$M = 10^{-4}$} \\
        \cmidrule(l){2-3}\cmidrule(l){5-6}\cmidrule(l){8-9}\cmidrule(l){11-12}
        & $L^{2}$ relative error $\delta\mathbf{u}$ & $L^{2}$ rate $\delta\mathbf{u}$ & & $L^{2}$ relative error $\delta\mathbf{u}$ & $L^{2}$ rate $\delta\mathbf{u}$ & & $L^{2}$ relative error $\delta\mathbf{u}$ & $L^{2}$ rate $\delta\mathbf{u}$ & & $L^{2}$ relative error $\delta\mathbf{u}$ & $L^{2}$ rate $\delta\mathbf{u}$ \\
	\midrule
        $10$ & $\num{4.85e-2}$ & & & $\num{4.88e-2}$ & & & $\num{4.88e-2}$ & & & $\num{5.90e-2}$ & \\
	\midrule
        $20$ & $\num{2.88e-3}$ & $4.1$ & & $\num{2.73e-3}$ & $4.2$ & & $\num{2.71e-3}$ & $4.2$ & & $\num{1.45e-2}$ & $2.0$ \\
	\midrule
        $40$ & $\num{1.23e-4}$ & $4.5$ & & $\num{1.30e-4}$ & $4.4$ & & $\num{1.30e-4}$ & $4.4$ & & $\num{5.56e-3}$ & $1.4$ \\
	\midrule
        $80$ & $\num{8.47e-6}$ & $\mathbf{3.9}$ & & $\num{1.10e-5}$ & $\mathbf{3.6}$ & & $\num{1.26e-5}$ & $\mathbf{3.4}$ & & & \\
	\bottomrule
    \end{tabularx}
    \caption{Convergence analysis for the traveling vortex test case using the \textbf{IMEX} method with the IMEX-RK(6,6,4) scheme of type II (Table \ref{tab:rk4_butch_type_II}) and polynomial degree $r = 3$.}
    \label{tab:vortex_fourth_order_IMEX_type_II}
\end{table}

\begin{table}[h!]
    \centering
    \footnotesize
    \begin{tabularx}{0.55\columnwidth}{rXXcXX}
	\toprule
        \multirow{1}{*}{$N_{el}$} & \multicolumn{2}{c}{$M = 10^{-1}$} & & \multicolumn{2}{c}{$M = 10^{-2}$} \\
	\cmidrule(l){2-3}\cmidrule(l){5-6}
        & $L^{2}$ relative error $\delta\mathbf{u}$ & $L^{2}$ rate $\delta\mathbf{u}$ & & $L^{2}$ relative error $\delta\mathbf{u}$ & $L^{2}$ rate $\delta\mathbf{u}$ \\
	\midrule
        $10$ & $\num{4.85e-2}$ & & & $\num{4.88e-2}$ & \\
	\midrule
        $20$ & $\num{2.88e-3}$ & $4.1$ & & $\num{2.74e-3}$ & $4.2$ \\
	\midrule
        $40$ & $\num{1.23e-4}$ & $4.5$ & & $\num{1.30e-4}$ & $4.4$ \\
	\midrule
        $80$ & $\num{8.47e-6}$ & $\mathbf{3.9}$ & & $\num{1.10e-5}$ & $\mathbf{3.6}$ \\
	\bottomrule
    \end{tabularx}
    \caption{Convergence analysis for the traveling vortex test case using the \textbf{SI-IMEX} method with the IMEX-RK(6,6,4) scheme of type II (Table \ref{tab:rk4_butch_type_II}) and polynomial degree $r = 3$. Note that for $M = 10^{-1}$ the advective Courant number is $C_{u} \approx 0.05$, while for $M = 10^{-2}$ the advective Courant number is $C_{u} \approx 0.025$ to achieve a stable solution.}
    \label{tab:vortex_fourth_order_SI_IMEX_type_II}
\end{table}

\subsubsection{Investigation of fixed point iterations}
\label{ssec:fixed_point}

In this Section, we analyze the impact of the value of the tolerance $\eta$ in the fixed point loop on the overall performance of the IMEX method depicted in Section \ref{ssec:IMEX}. For the sake of brevity, we focus on $M = 10^{-1}$ and on $M = 10^{-3}$. We consider the third order scheme of type II (Table \ref{tab:rk3_butch_type_II}), but analogous considerations hold for the other schemes. For moderate values of the Mach number, a sufficiently low value of the tolerance $\eta$ has to be chosen in order to achieve full convergence and the number of fixed point iterations depends on the value of $\eta$ (Table \ref{tab:vortex_third_order_IMEX_eta}). For what concerns low values of the Mach number, one fixed point iteration is sufficient to achieve full convergence, independently of $\eta$ (Table \ref{tab:vortex_third_order_IMEX_eta}). This explains why the IMEX method and the SI-IMEX method behave similarly in this regime in terms of computational cost. These considerations are further confirmed by the experimental contraction rate (ECR), as defined in \cite{noelle:2014}, i.e.
\begin{equation}\label{eq:ECR}
   \text{ECR}(\tilde{k}) = \frac{\left\|p_{\tilde{k} + 1} - p_{\tilde{k}}\right\|_{\infty}}{\left\|p_{\tilde{k}} - p_{\tilde{k} - 1}\right\|_{\infty}}, \tilde{k} > 0,
\end{equation}
which is reported in Table \ref{tab:vortex_third_order_IMEX_ECR_EOC}. One can easily notice that the error is significantly reduced after three fixed point iterations as evident by the small value of the contraction rate. This results confirms that very few iterations of the fixed point loop are sufficient to obtain a satisfactory solution. We also monitor the experimental order of convergence of the fixed point method (EOC), defined as
\begin{equation}\label{eq:EOC}
    \text{EOC}(\tilde{k}) = \frac{\log\left(\frac{\left\|p_{\tilde{k} + 1} - p_{\tilde{k}}\right\|_{\infty}}{\left\|p_{\tilde{k}} - p_{\tilde{k} - 1}\right\|_{\infty}}\right)}{\log\left(\frac{\left\|p_{\tilde{k}} - p_{\tilde{k} - 1}\right\|_{\infty}}{\left\|p_{\tilde{k} - 1} - p_{\tilde{k} - 2}\right\|_{\infty}}\right)}, \tilde{k} > 1
\end{equation}
and we notice that, as expected, the fixed point method is of order 1 (Table \ref{tab:vortex_third_order_IMEX_ECR_EOC}). From now on, we set $\eta = 10^{-6}$ in the following test cases.

\begin{table}[h!]
    \centering
    \footnotesize
    \begin{tabularx}{0.6\columnwidth}{rXXcXX}
	\toprule
        \multirow{1}{*}{$\eta$} & \multicolumn{2}{c}{$M = 10^{-1}$} & & \multicolumn{2}{c}{$M = 10^{-4}$} \\
	\cmidrule(l){2-3}\cmidrule(l){5-6}
        & $L^{2}$ relative error $\delta\mathbf{u}$ & N. iters & & $L^{2}$ relative error $\delta\mathbf{u}$ & N. iters \\
	\midrule
        $10^{-4}$ & $\num{3.35e-4}$ & $1$ & & $\num{2.90e-4}$ & $1$ \\
	\midrule
        $10^{-6}$ & $\num{3.17e-4}$ & $3$ & & $\num{2.90e-4}$ & $1$ \\
	\midrule
        $10^{-10}$ & $\num{3.17e-4}$ & $5$ & & $\num{2.90e-4}$ & $1$ \\
	\bottomrule
    \end{tabularx}
    \caption{Traveling vortex test case, investigation of the fixed point loop using the third order scheme of type II (Table \ref{tab:rk3_butch_type_II}) and polynomial degree $r = 2$. Here, $\eta$ denotes the tolerance for the stopping criterion of the fixed point loop \eqref{eq:fixed_point_tol}, while N. iters represents the average number of iterations in the fixed point loop. The computational mesh of the reported results is composed by $80 \times 80 = 6400$ elements.}
    \label{tab:vortex_third_order_IMEX_eta}
\end{table}

\begin{table}[h!]
    \centering
    \footnotesize
    \begin{tabularx}{0.75\columnwidth}{rXXXcXXX}
	\toprule
        \multirow{1}{*}{$\tilde{k}$} & \multicolumn{2}{c}{1st time step} & & & & \multicolumn{2}{c}{200th time step} \\
	\cmidrule(l){2-4}\cmidrule(l){6-8}
        & $\left\|p_{\tilde{k} + 1} - p_{\tilde{k}}\right\|_{\infty}$ & ECR & EOC & & $\left\|p_{\tilde{k} + 1} - p_{\tilde{k}}\right\|_{\infty}$ & ECR & EOC \\
	\midrule
        $0$ & $\num{1.02e-5}$ & $-$ & $-$ & & $\num{5.84e-6}$ & $-$ & $-$ \\
	\midrule
        $1$ & $\num{5.42e-6}$ & $0.53$ & $-$ & & $\num{5.34e-7}$ & $0.91$ & $-$ \\
	\midrule
        $2$ & $\num{6.31e-8}$ & $0.012$ & $7.1$ & & $\num{4.25e-8}$ & $0.0080$ & $54$ \\
	\midrule
        $3$ & $\num{1.10e-9}$ & $0.017$ & $0.91$ & & $\num{3.03e-10}$ & $0.0071$ & $1.02$ \\
        \midrule
        $4$ & $\num{3.42e-11}$ & $0.031$ & $0.86$ & & $\num{4.68e-12}$ & $0.015$ & $0.84$ \\
        \midrule
        $5$ & $\num{2.75e-12}$ & $0.08$ & $0.73$ & & $\num{1.42e-13}$ & $0.03$ & $0.84$ \\
        \midrule
        $6$ & $\num{1.86e-13}$ & $0.07$ & $\mathbf{1.07}$ & & $\num{7.55e-15}$ & $0.05$ & $\mathbf{0.84}$ \\
	\bottomrule
    \end{tabularx}
    \caption{Traveling vortex test case at $M = 10^{-1}$, investigation of the fixed point loop using the third order scheme of type II (Table \ref{tab:rk3_butch_type_II}) and polynomial degree $r = 2$. Here $\tilde{k}$ is the index of the fixed point iteration, ECR represents the experimental convergence rate \eqref{eq:ECR}, whereas EOC denotes the experimental order of convergence \eqref{eq:EOC}. The computational mesh of the reported results is composed by $80 \times 80 = 6400$ elements.}
    \label{tab:vortex_third_order_IMEX_ECR_EOC}
\end{table}

\subsection{Sod shock tube for the Peng-Robinson EOS}
\label{ssec:Sod_PR_EOS}

\begin{figure}[h!]
    \centering
    \begin{subfigure}{0.475\textwidth}
        \centering
        \includegraphics[width = 0.95\textwidth]{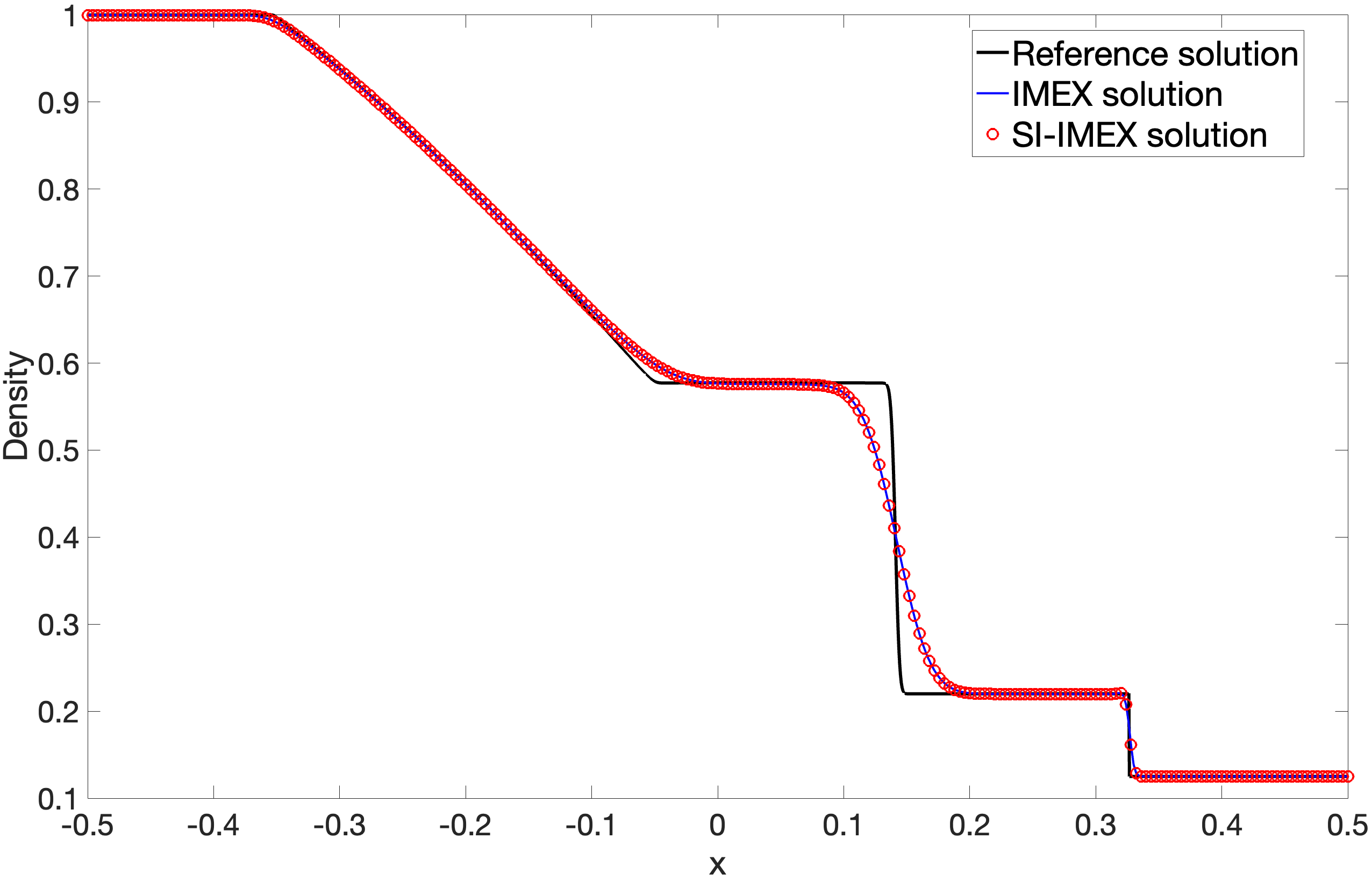} a)
    \end{subfigure}
    \begin{subfigure}{0.475\textwidth}
        \centering
        \includegraphics[width = 0.95\textwidth]{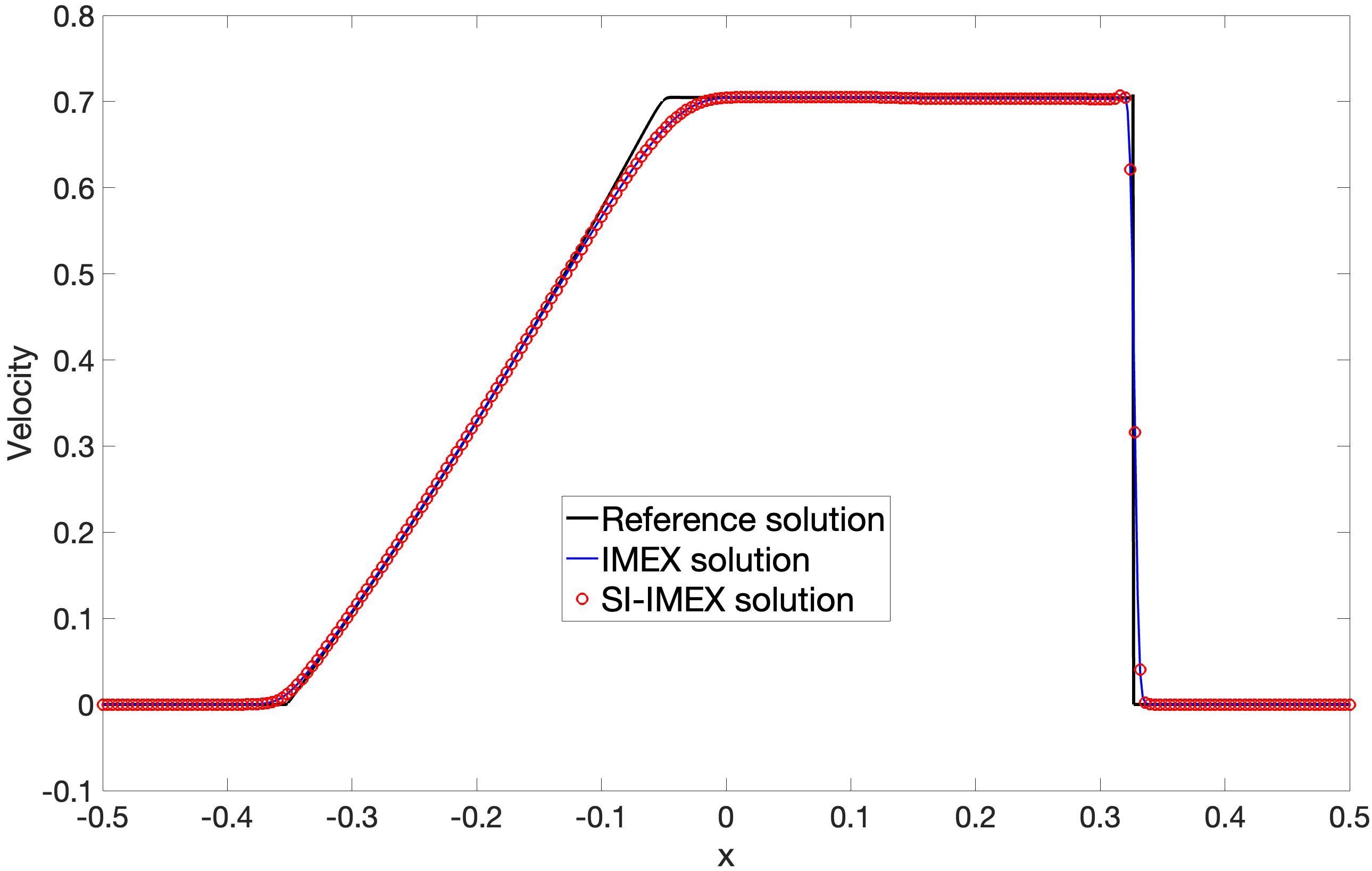} b)
    \end{subfigure}
    \begin{subfigure}{0.475\textwidth}
        \centering
        \includegraphics[width = 0.95\textwidth]{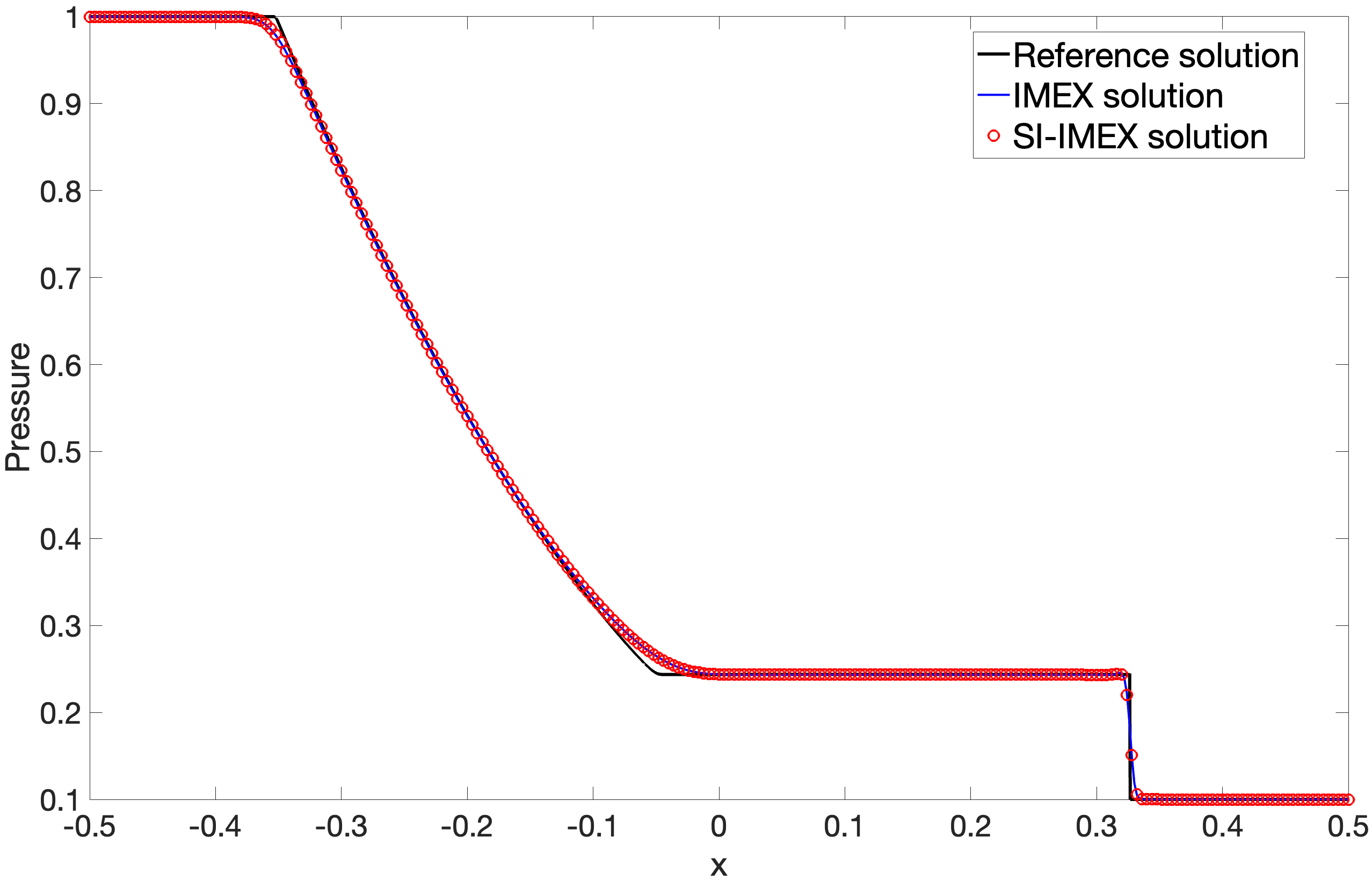} c)
    \end{subfigure}
    \caption{Sod shock tube for the Peng-Robinson EOS \eqref{eq:int_energy_general_cubic_eos_constant} with polynomial degree $r = 0$, a) density, b) velocity, c) pressure. The continuous black lines represent the reference solution computed using the third order optimal explicit strong stability-preserving scheme in \cite{gottlieb:2001}, the continuous blue show the results obtained with the IMEX method using the IMEX-RK(2,2,2) scheme in Table \ref{tab:rk2_butch}, whereas the red dots report the numerical results obtained using the SI-IMEX method.}
    \label{fig:Sod_PR_r_0}
\end{figure}

The main aim of the proposed time discretization methods is to use them for low Mach number flows so as to avoid the acoustic CFL restriction. In the case of high Mach number flows, acoustic waves are not negligible. Hence, one is interested in resolving both acoustic and material waves and explicit time discretization schemes are well suited to achieve this goal. However, we show that, when coupled with a monotonicity preserving spatial discretization, both IMEX and SI-IMEX methods can be effective also for high Mach number flows. We consider the Sod shock tube problem \cite{sod:1978} for the Peng-Robinson EOS \eqref{eq:int_energy_general_cubic_eos_constant}. The computational domain is $\Omega = \left(-0.5, 0.5\right)$, whereas the final time is $T_{f} = 0.2$. The initial condition reads as follows: 
\begin{equation}
    \left(\rho, u, p\right)(x, 0) = 
    \begin{cases}
        \left(1,0,1\right) \qquad &\text{if } x < 0 \\
        \left(0.125,0,0.1\right) \qquad &\text{if } x > 0. \\
    \end{cases}
\end{equation}
Dirichlet boundary conditions are imposed. Following \cite{dumbser:2016a}, we take $c_{v} = 1$ and $R_{g} = 0.4$. Moreover, we take $b = 0.5$ and $a(T) = 0.5/\sqrt{T}$. We employ the IMEX-RK(2,2,2) scheme (Table \ref{tab:rk2_butch}) with polynomial degree $r = 0$ so as to avoid the oscillations that arise in the case of discontinuous solution when using high-order space discretization methods. The computational mesh is composed by $N_{el} = 500$ elements and the time step is $\Delta t \approx 1.33 \times 10^{-3}$, yielding a maximum acoustic Courant number $C \approx 1.23$ and $C_{u} \approx 0.47$. In particular, we point out that the acoustic Courant number is greater than one. A reference solution is computed using the optimal third order explicit strong stability-preserving scheme presented in \cite{gottlieb:2001} with $N_{el} = 32000$. An excellent agreement is established between the IMEX and the SI-IMEX method and a good agreement is established with the reference solution (Figure \ref{fig:Sod_PR_r_0}). Moreover, one can start appreciating the effectiveness of the linearization presented in Section \ref{ssec:num_eos}. We will further discuss this point in Section \ref{ssec:KH_instability}.

Finally, we point out that the employed spatial discretization is not TVD for $r > 0$. Hence, spurious oscillations arise in the case of discontinuous solutions. A detailed discussion of possible approaches to overcome this issue is not in the scope of the present work. However, a number of approaches have been proposed in the literature to obtain essentially monotone schemes using high order DG methods, see, e.g., \cite{dumbser:2016b, orlando:2023a}.

\subsection{Flow in an open tube}
\label{ssec:open_tube}

Next, we consider the test case III originally proposed in \cite{klein:1995} for an ideal gas and also employed in \cite{orlando:2025}. The domain is $\Omega = \left(0, 10\right)$ with a time-dependent density and velocity prescribed at the left-end, while a time-dependent outflow pressure with a large amplitude variation is imposed at the right-end. More specifically, the initial conditions read as follows:
\begin{equation}
    \left(\rho, u, p\right)\left(x, 0\right) = \left(1, 1, 1\right),
\end{equation}
while the boundary conditions are
\begin{subequations}
    \begin{eqnarray}
	\rho\left(0, t\right) &=& 1 + \frac{3}{10}\sin\left(4t\right) \\
	u\left(0, t\right) &=& 1 + \frac{1}{2}\sin\left(2t\right) \\
        p\left(L, t\right) &=& 1 + \frac{1}{4}\sin\left(3t\right), \label{eq:p_init_open_tube}
    \end{eqnarray}
\end{subequations}
with $L = 10$. The final time is $T_{f} = 7.47$, whereas the Mach number is set to $M = 10^{-4}$. Since \eqref{eq:dive_u} reduces to
\begin{equation}\label{eq:dive_u_IG_tube}
    \dive\bar{\mathbf{u}} = -\frac{1}{\gamma\bar{p}}\frac{d\bar{p}}{dt},
\end{equation}
the velocity field is not solenoidal as $M \to 0$, and, in one space dimension, it is a linear function of the space with a given time-dependent slope and boundary value at $x = 0$ \cite{klein:1995, orlando:2025}. In particular, the leading order solution for the velocity reads as follows:
\begin{equation}
    \bar{u}(x,t) = u(0,t) - \frac{1}{\gamma p(L, t)}\frac{dp(L,t)}{dt}x = 1 + \frac{1}{2}\sin\left(2t\right) - \frac{3\cos\left(3t\right)}{4\gamma\left(1 + \frac{1}{4}\sin\left(3t\right)\right)}x.
\end{equation}
Moreover, the leading order relation \eqref{eq:continuity_limit} reduces to
\begin{equation}
    \frac{D\log\bar{\rho}}{Dt} = -\frac{\partial\bar{u}}{\partial x} = \frac{1}{\gamma p(L, t)}\frac{dp(L,t)}{dt} = \frac{3\cos\left(3t\right)}{4\gamma\left(1 + \frac{1}{4}\sin\left(3t\right)\right)},
\end{equation}
where $\frac{D\log\bar{\rho}}{Dt} = \frac{\partial\log\bar{\rho}}{\partial t} + \bar{u}\frac{\partial\log\bar{\rho}}{\partial x}$.

\begin{figure}[h!]
    \centering
    \begin{subfigure}{0.475\textwidth}
	\centering
        \includegraphics[width = 0.95\textwidth]{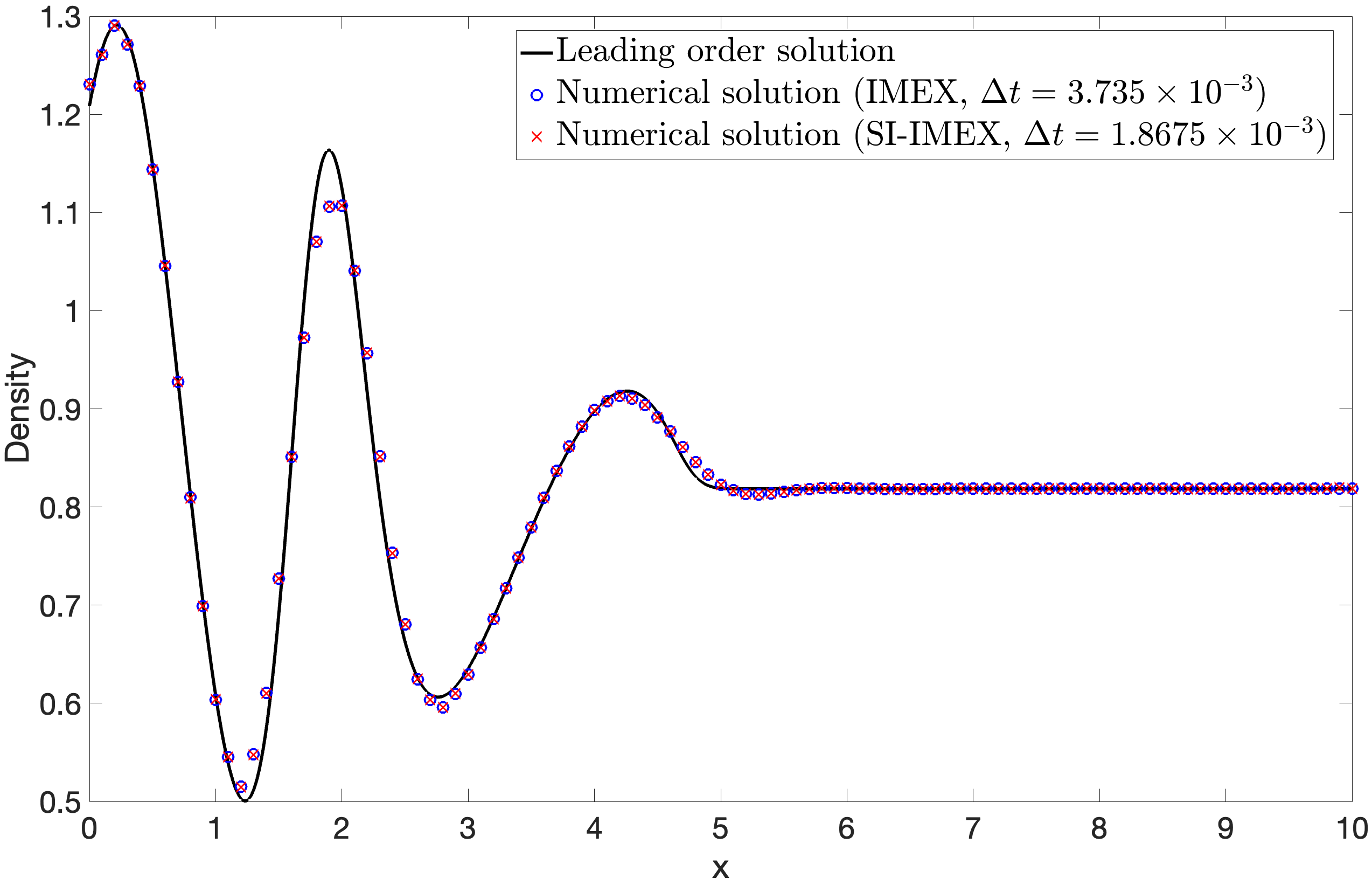}
    \end{subfigure}
    \begin{subfigure}{0.475\textwidth}
	\centering
        \includegraphics[width = 0.95\textwidth]{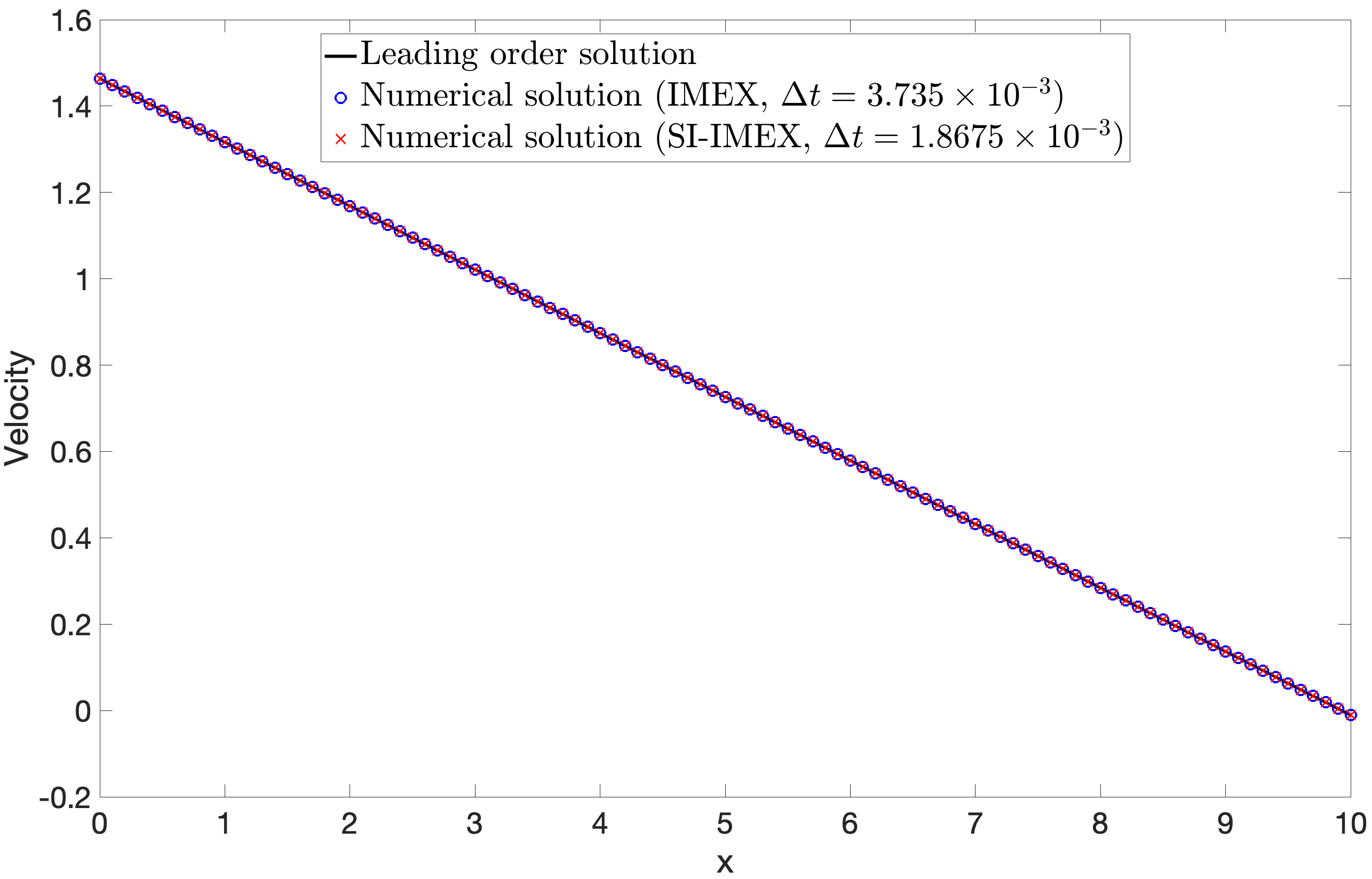}
    \end{subfigure}
    \begin{subfigure}{0.475\textwidth}
	\centering
        \includegraphics[width = 0.95\textwidth]{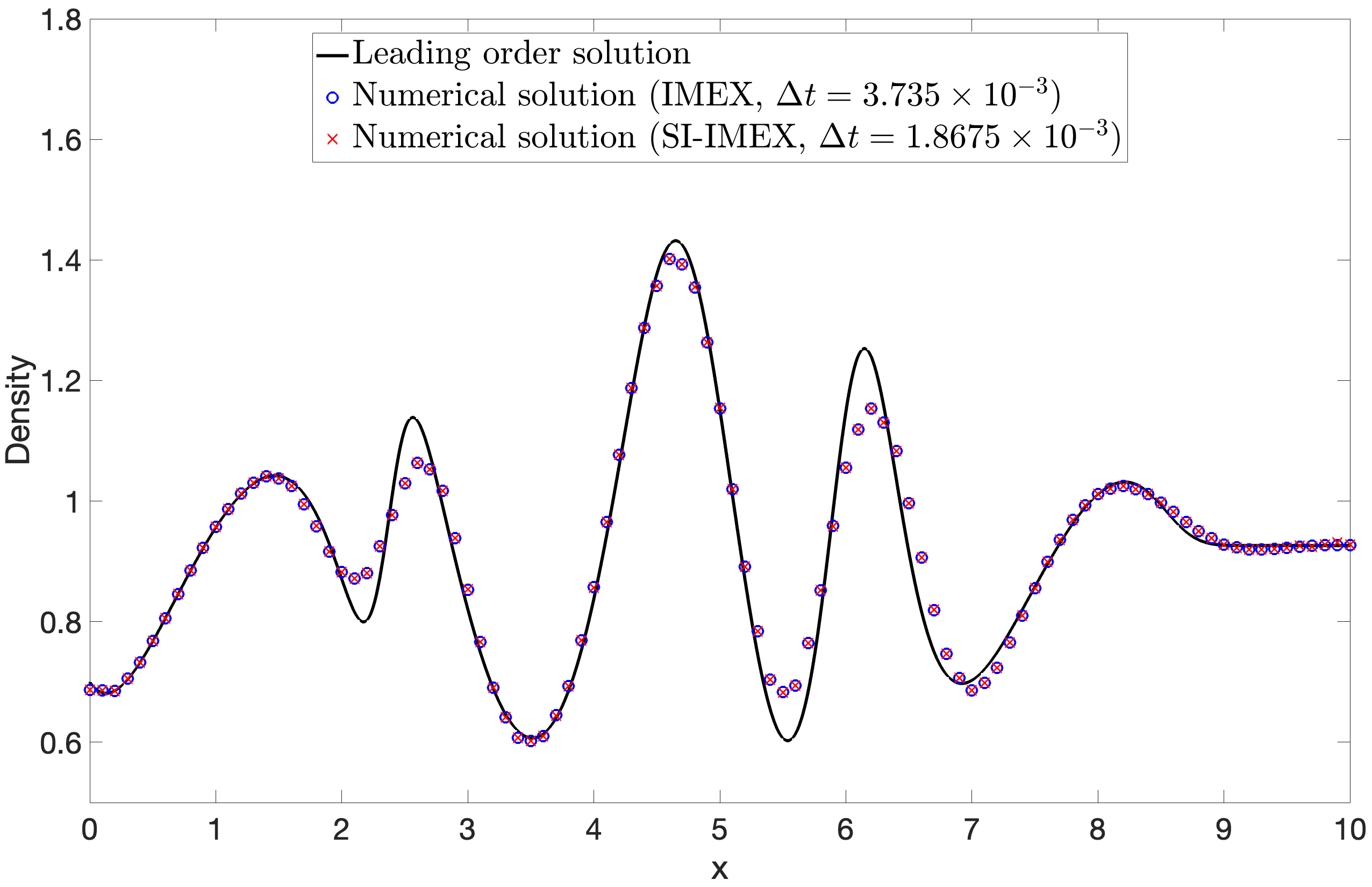}
    \end{subfigure}
    \begin{subfigure}{0.475\textwidth}
	\centering
        \includegraphics[width = 0.95\textwidth]{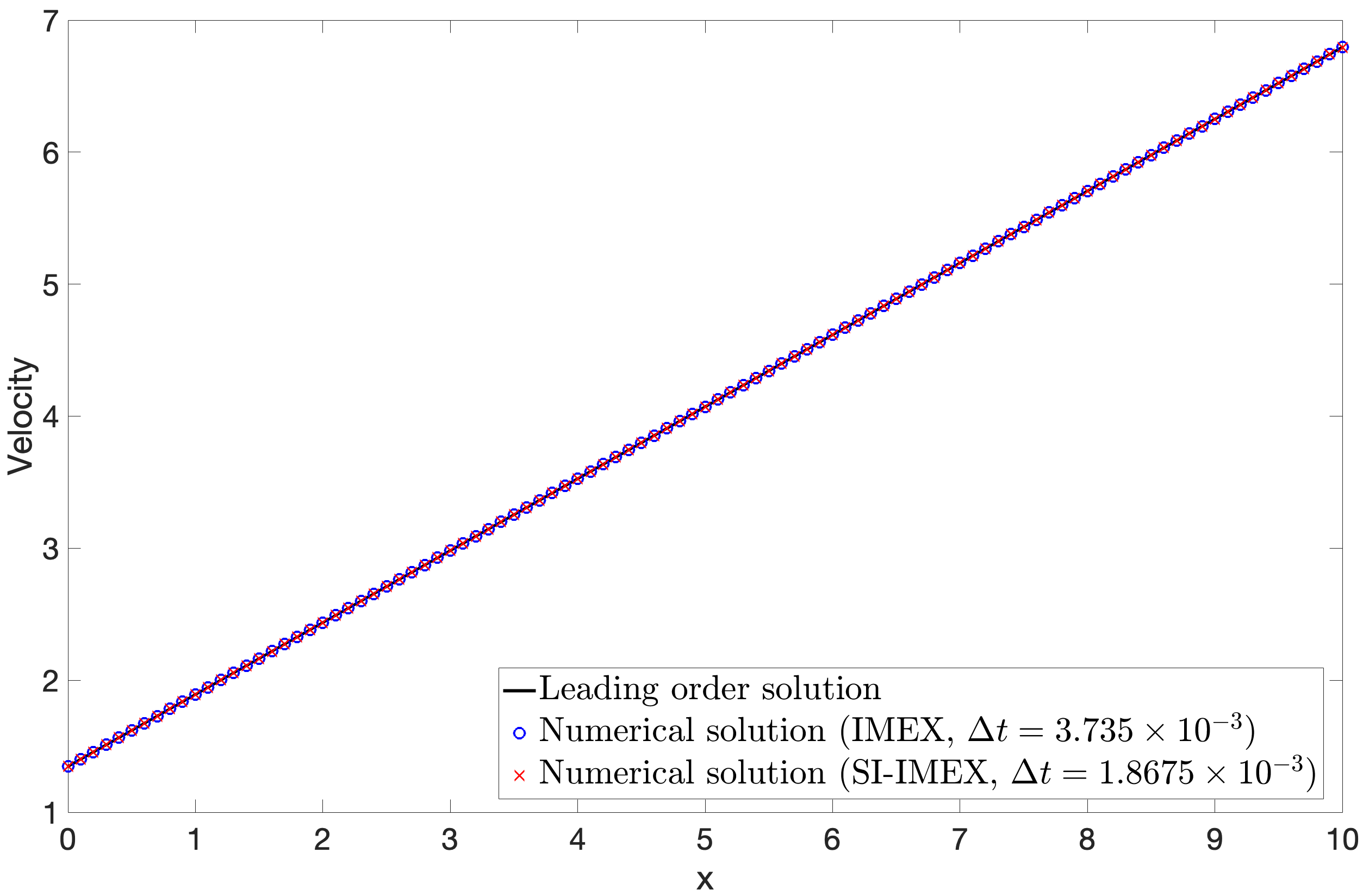}
    \end{subfigure}
    \caption{Open tube test case with the ideal gas law \eqref{eq:ideal_gas}, comparison between the IMEX method and the SI-IMEX method using the second order scheme (Table \ref{tab:rk2_butch}) and polynomial degree $r = 1$ with $N_{el} = 50$. Top: results at $t = \frac{T_{f}}{2} = 3.735$. Bottom: results at $t = T_{f} = 7.47$. Left: density. Right: velocity. The continuous black lines show the leading order solution as $M \to 0$, the blue dots report the numerical results with the IMEX method and $\Delta t = 3.735 \times 10^{-3}$, whereas the red crosses represent the results with the SI-IMEX method and $\Delta t = 1.8675 \times 10^{-3}$.}
    \label{fig:open_tube_IG_comparison} 
\end{figure}

First, we employ the second order IMER-RK(2,2,2) (Table \ref{tab:rk2_butch}) with polynomial degree $r = 1$. We consider a number of elements $N_{el} = 50$, whereas the time step is $\Delta t = 3.735 \times 10^{-3}$, leading to a maximum advective Courant number $C_{u} \approx 0.13$ and a maximum acoustic Courant number $C \approx 310$. The results at $t = \frac{T_{f}}{2}$ and at $t = T_{f}$ obtained with the IMEX method are those expected by the asymptotic analysis \cite{klein:1995, orlando:2025} for both the density and velocity profiles (Figure \ref{fig:open_tube_IG_comparison}), while the SI-IMEX method does not achieve a stable solution employing this time step. We also notice that no significant issue seems to arise in the low Mach regime.

The SI-IMEX method requires a smaller time step to achieve a stable solution. More specifically, we take $\Delta t = 1.8675 \times 10^{-3}$, which leads to a maximum acoustic Courant number $C \approx 155$ and a maximum advective Courant number $C_{u} \approx 0.065$. The use of the SI-IMEX method yields a computational time saving of around $37\%$ at fixed time step. One can easily notice that a good agreement with the leading order solution is established for both the density and the velocity (Figure \ref{fig:open_tube_IG_comparison}).

Next, we consider the third order scheme IMEX-RK(4,3,3) of type I (Table \ref{tab:rk3_butch_type_I}) and the third order scheme IMEX-RK(3,3,3) of type II (Table \ref{tab:rk3_butch_type_II}) using polynomial degree $r = 2$ with $N_{el} = 50$. We take $\Delta t =  1.8675 \times 10^{-3}$, yielding a maximum advective Courant number $C_{u} \approx 0.14$ and a maximum acoustic Courant number $C \approx 292$. Similarly to the second order scheme, a good agreement with the leading order solution is established for both the time discretization schemes, in particular for the scheme of type I, while the SI-IMEX method does not achieve a stable solution. We take $\Delta t = 4.66875 \times 10^{-4}$, for which a stable solution for the scheme of type II is obtained. On the contrary, a stable solution is not achieved for the scheme of type I. This is likely related to the employed boundary conditions. The AP property and the AA property of the SI-IMEX method indeed has been proven considering periodic or no-flux boundary conditions \cite{arun:2021, boscarino:2022, huang:2022}, while we are considering time-dependent Dirichlet boundary conditions. Hence, the system is not autonomous and therefore it is outside of the theoretical framework of the SI-IMEX method depicted in Section \ref{ssec:SI_IMEX}. Since $c \neq \tilde{c}$ for schemes of type I, $\mathbf{U}_{E}$ and $\mathbf{U}_{I}$ are computed at different time instants, which seems to cause issues for time-dependent Dirichlet boundary conditions. An excellent agreement with the leading order solution is established for both the density and velocity field for $\Delta t = 4.66875 \times 10^{-4}$. Hence, for configurations which involve large variations of density and pressure, and time-dependent boundary conditions, the IMEX method is globally more robust and allows for sizeable larger time steps.

\begin{figure}[h!]
    \centering
    \begin{subfigure}{0.475\textwidth}
	\centering
        \includegraphics[width = 0.95\textwidth]{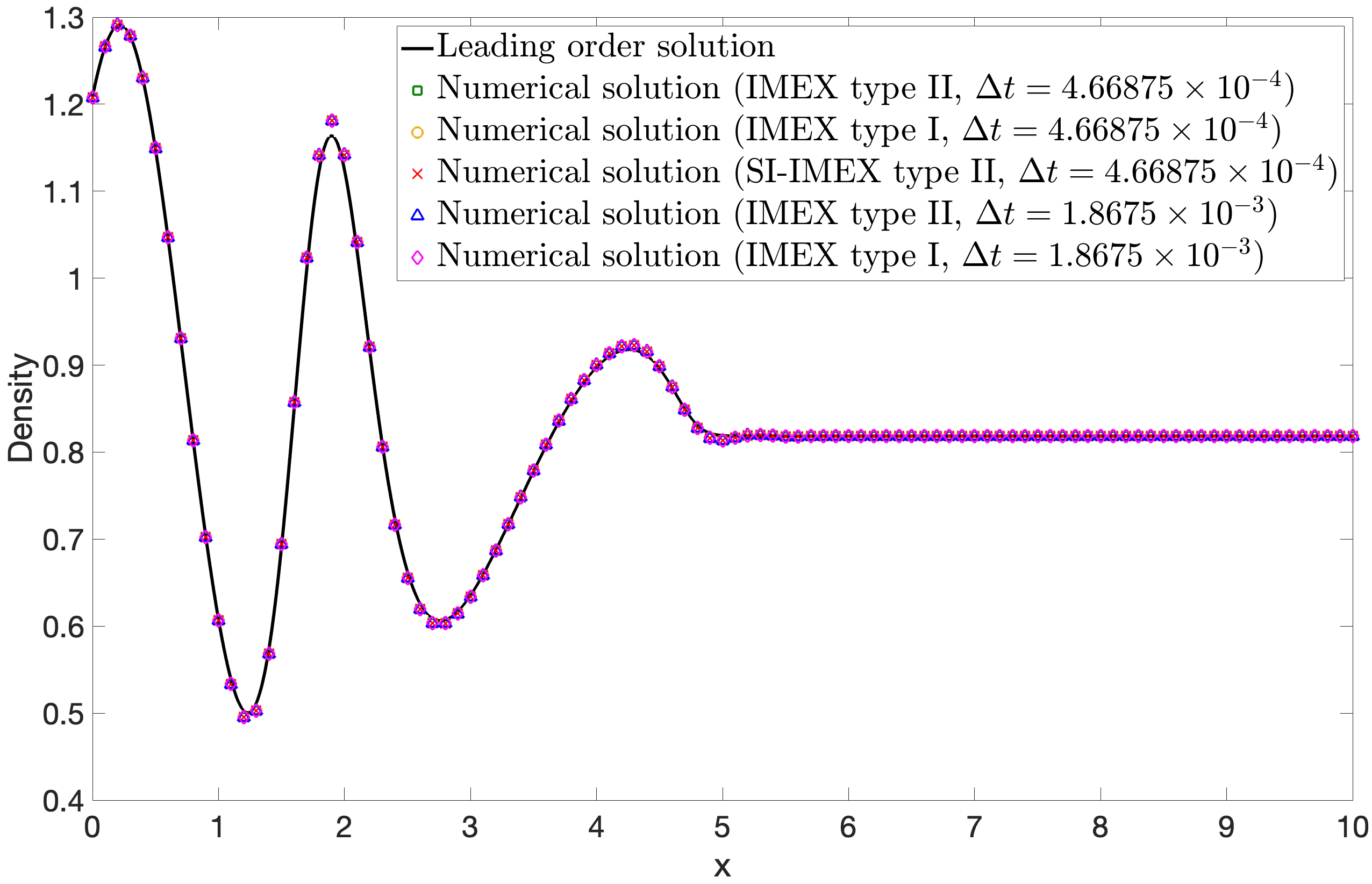}
    \end{subfigure}
    \begin{subfigure}{0.475\textwidth}
	\centering
        \includegraphics[width = 0.95\textwidth]{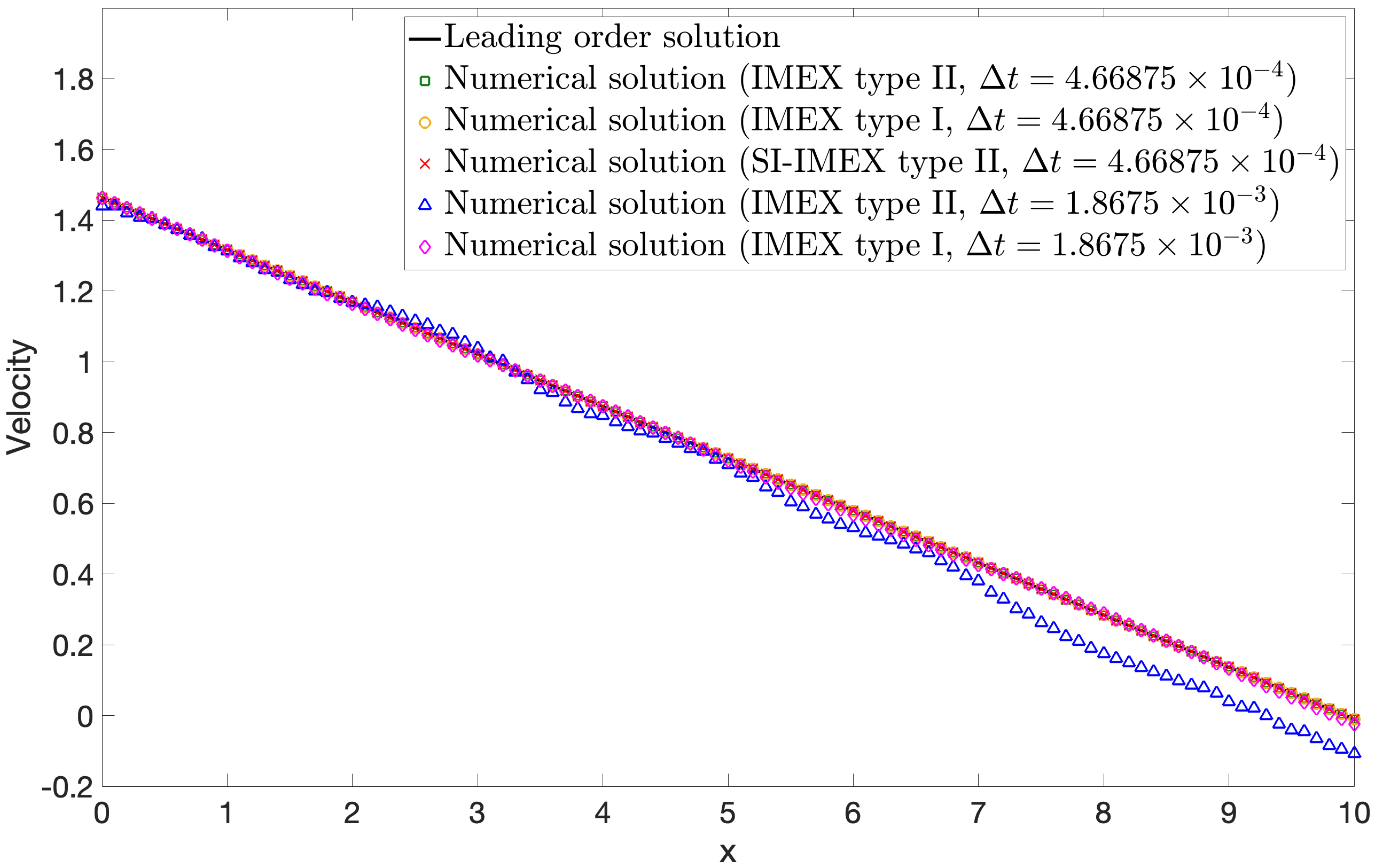}
    \end{subfigure}
    \begin{subfigure}{0.475\textwidth}
	\centering
        \includegraphics[width = 0.95\textwidth]{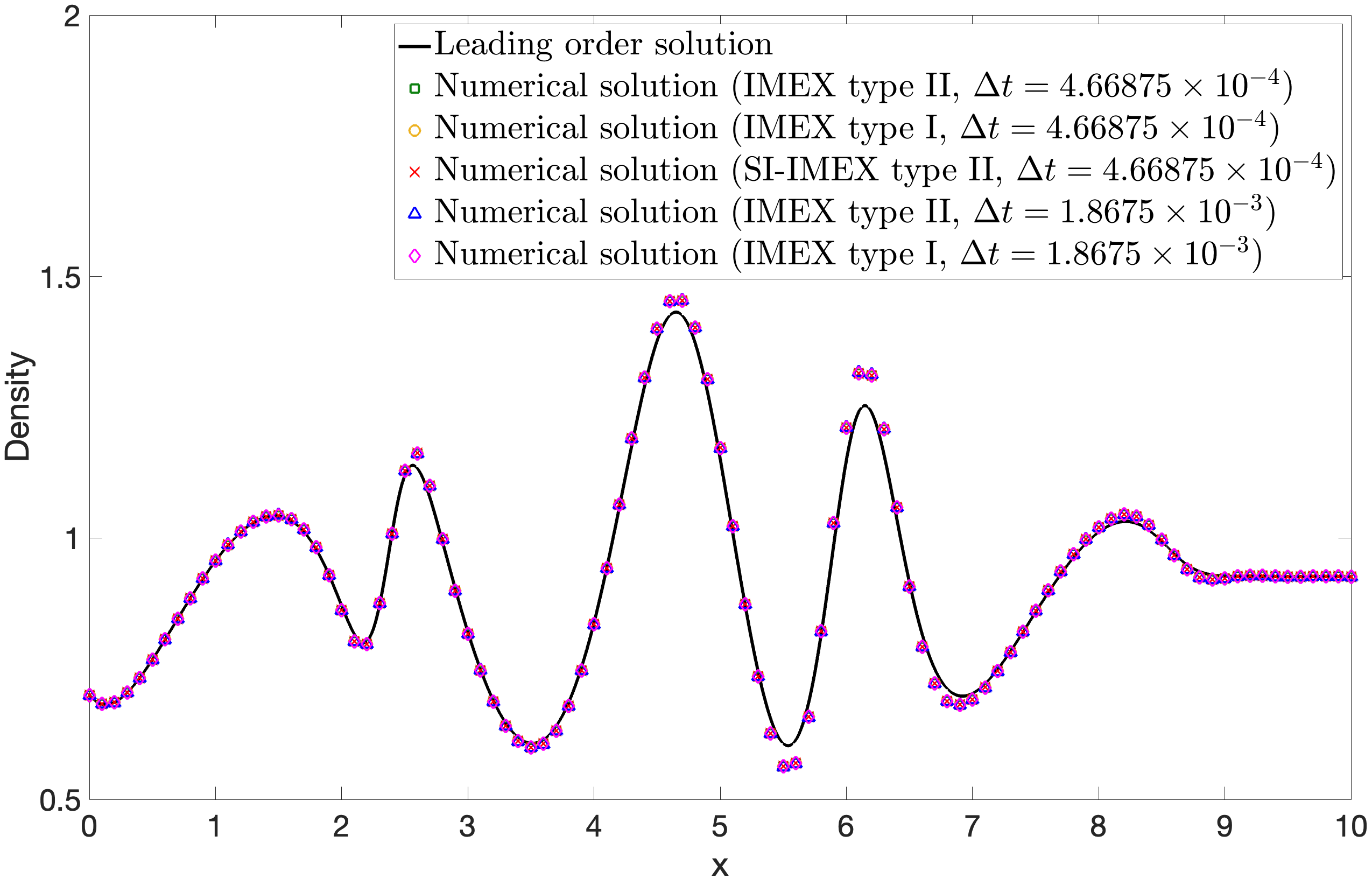}
    \end{subfigure}
    \begin{subfigure}{0.475\textwidth}
	\centering
        \includegraphics[width = 0.95\textwidth]{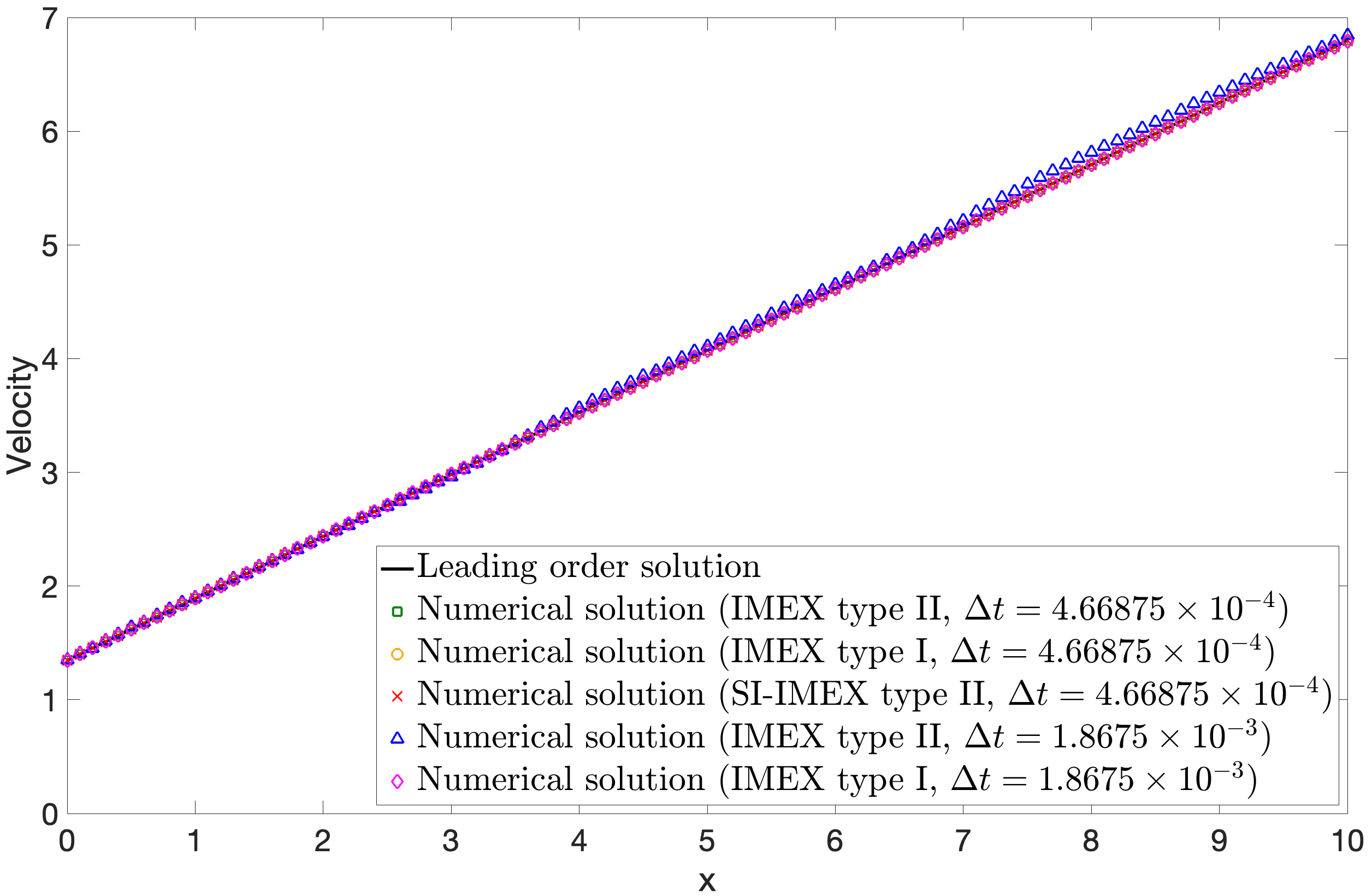}
    \end{subfigure}
    \caption{Open tube test case with the ideal gas law \eqref{eq:ideal_gas}, comparison between the IMEX method and the SI-IMEX method using the third order schemes (Tables \ref{tab:rk3_butch_type_I}-\ref{tab:rk3_butch_type_II}) in combination with polynomial degree $r = 2$ with $N_{el} = 50$. Top: results at $t = \frac{T_{f}}{2} = 3.735$. Bottom: results at $t = T_{f} = 7.47$. Left: density. Right: velocity. The continuous black lines show the leading order solution as $M \to 0$, the green squares report the numerical results with the IMEX method of type II and $\Delta t = 2.334375 \times 10^{-4}$, the orange circles show the numerical results with the IMEX method of type I and $\Delta t = 2.334375 \times 10^{-4}$, the red crosses represent the results with the SI-IMEX method of type II and $\Delta t = 2.334375 \times 10^{-4}$, the blue triangles show the results with the IMEX method of type II and $\Delta t = 1.245 \times 10^{-3}$, while the magenta diamonds  report the results with the IMEX method of type I and $\Delta t = 1.245 \times 10^{-3}$.}
    \label{fig:open_tube_IG_third_order} 
\end{figure}

\begin{figure}[h!]
    \centering
    \begin{subfigure}{0.475\textwidth}
	\centering
        \includegraphics[width = 0.95\textwidth]{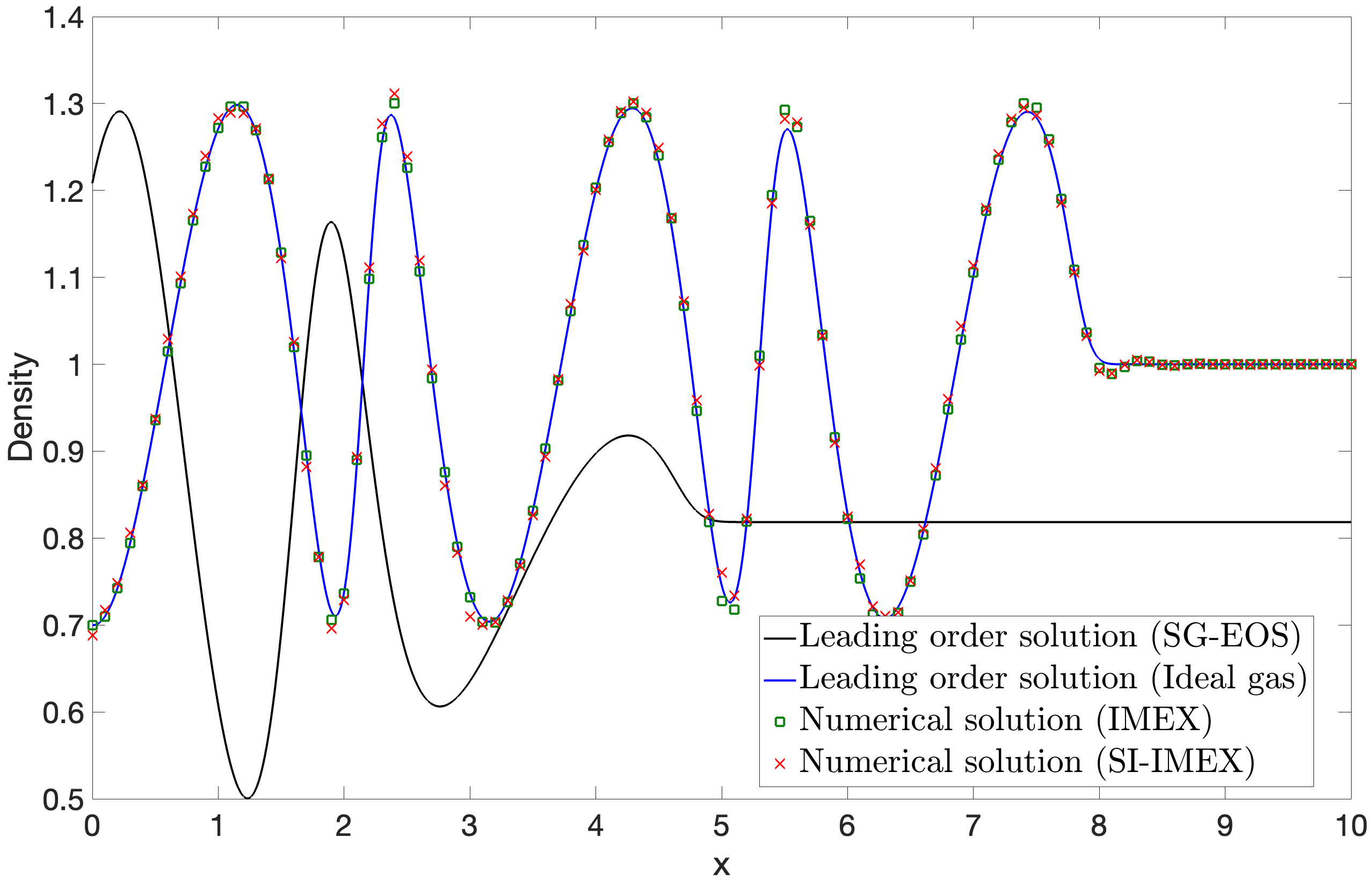}
    \end{subfigure}
    \begin{subfigure}{0.475\textwidth}
	\centering
        \includegraphics[width = 0.95\textwidth]{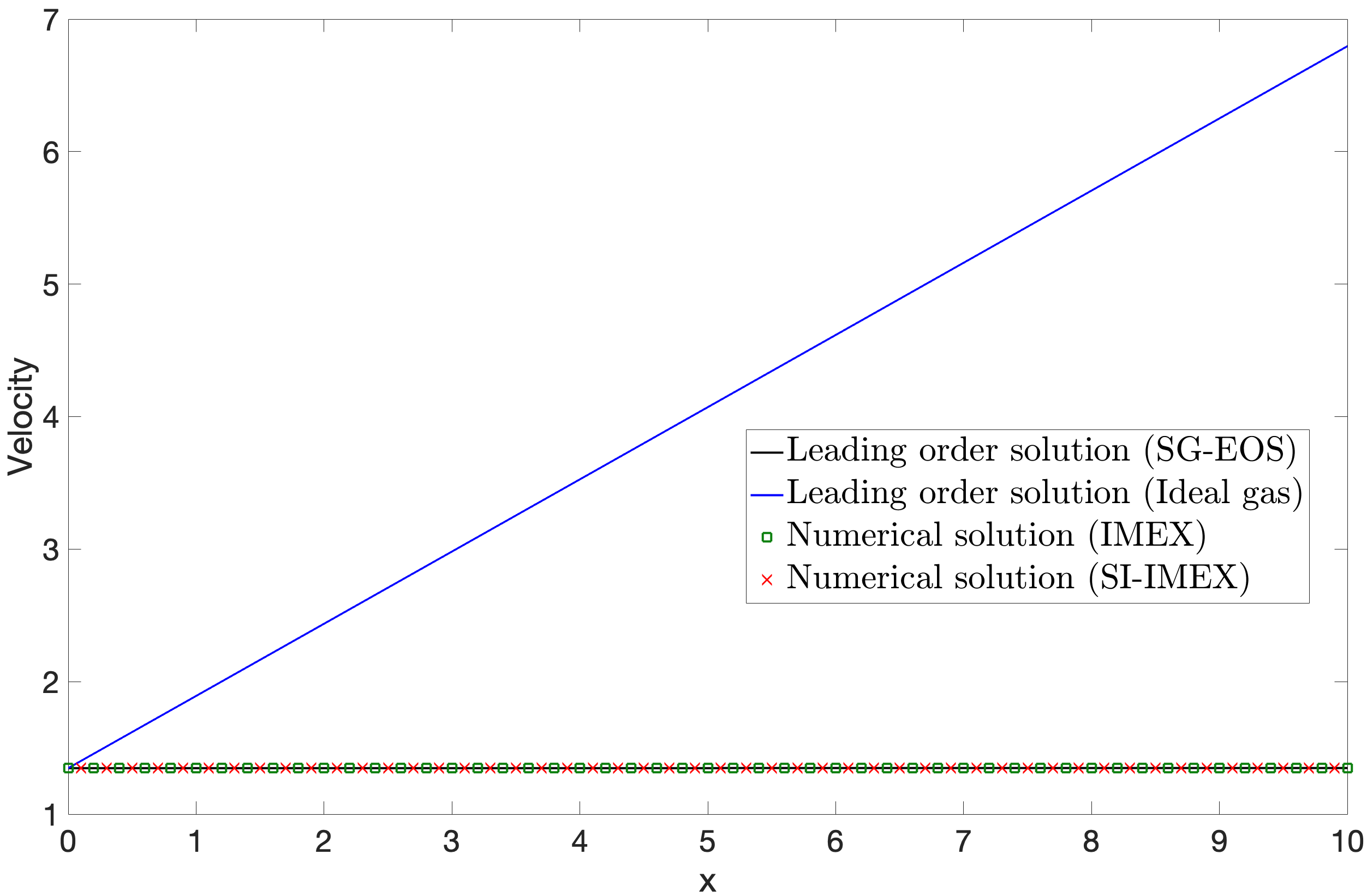}
    \end{subfigure}
    \caption{Open tube test case with the SG-EOS \eqref{eq:sg_eos}, comparison between the IMEX method and the SI-IMEX method using the third order schemes (Table \ref{tab:rk3_butch_type_II}) and polynomial degree $r = 2$ with $N_{el} = 50$. Results at $t = T_{f} = 7.47$. Left: density. Right: velocity. The continuous black lines show the leading order solution as $M \to 0$, the continuous blue lines report the leading order solution as $M \to 0$ for the ideal gas with $\gamma = 1.4$, the green squares report the numerical results with the IMEX method, while the red crosses represent the results with the SI-IMEX method.}
    \label{fig:open_tube_SG} 
\end{figure}

Next, we consider an extension of this test case using the stiffened gas equation of state (SG-EOS) \eqref{eq:sg_eos}. We take $\gamma = 4.4, \pi_{\infty} = 6.8 \times 10^{2}$, and $q_{\infty} = 0$ in \eqref{eq:sg_eos}. We employ the third-order scheme of type II (Table \ref{tab:rk3_butch_type_II}) in combination with $r = 2$. We take $\Delta t = 2.334375 \times 10^{-4}$, yielding a maximum acoustic Courant number $C \approx 1530$. An excellent agreement is established between the IMEX method and the SI-IMEX method (Figure \ref{fig:open_tube_SG}). Moreover, one can easily notice that the leading order solution changes with the equation of state and modifying its parameters \cite{orlando:2025} (Figure \ref{fig:open_tube_SG}). In particular, relation \eqref{eq:dive_u} for the SG-EOS reduces to
\begin{equation}
    \dive\bar{\mathbf{u}} = -\frac{1}{\gamma\left(\bar{p} + \pi_{\infty}\right)}\frac{d\bar{p}}{dt}.
\end{equation}
and therefore the leading order solution for the velocity reads as follows:
\begin{equation}
    \bar{u}(x,t) = 1 + \frac{1}{2}\sin\left(2t\right) - \frac{3\cos\left(3t\right)}{4\left(\gamma + \pi_{\infty}\right)\left(1 + \frac{1}{4}\sin\left(3t\right)\right)}x.
\end{equation}
Hence, since from \eqref{eq:p_init_open_tube} $\bar{p} \ge \frac{3}{4}$ and $\left|\frac{d\bar{p}}{dt}\right| \le \frac{3}{4}$, we obtain
\begin{equation}
    \left|\dive\bar{\mathbf{u}}\right| \le \frac{1}{\gamma\left(\frac{3}{4} + \pi_{\infty}\right)}\frac{3}{4} \approx 2.5 \times 10^{-4},
\end{equation}
meaning that the velocity field is almost constant. 

\subsection{Kelvin-Helmholtz instability at low Mach}
\label{ssec:KH_instability}

\begin{figure}[h!]
    \centering
    \begin{subfigure}{0.475\textwidth}
	\centering
        \includegraphics[width = 0.85\textwidth]{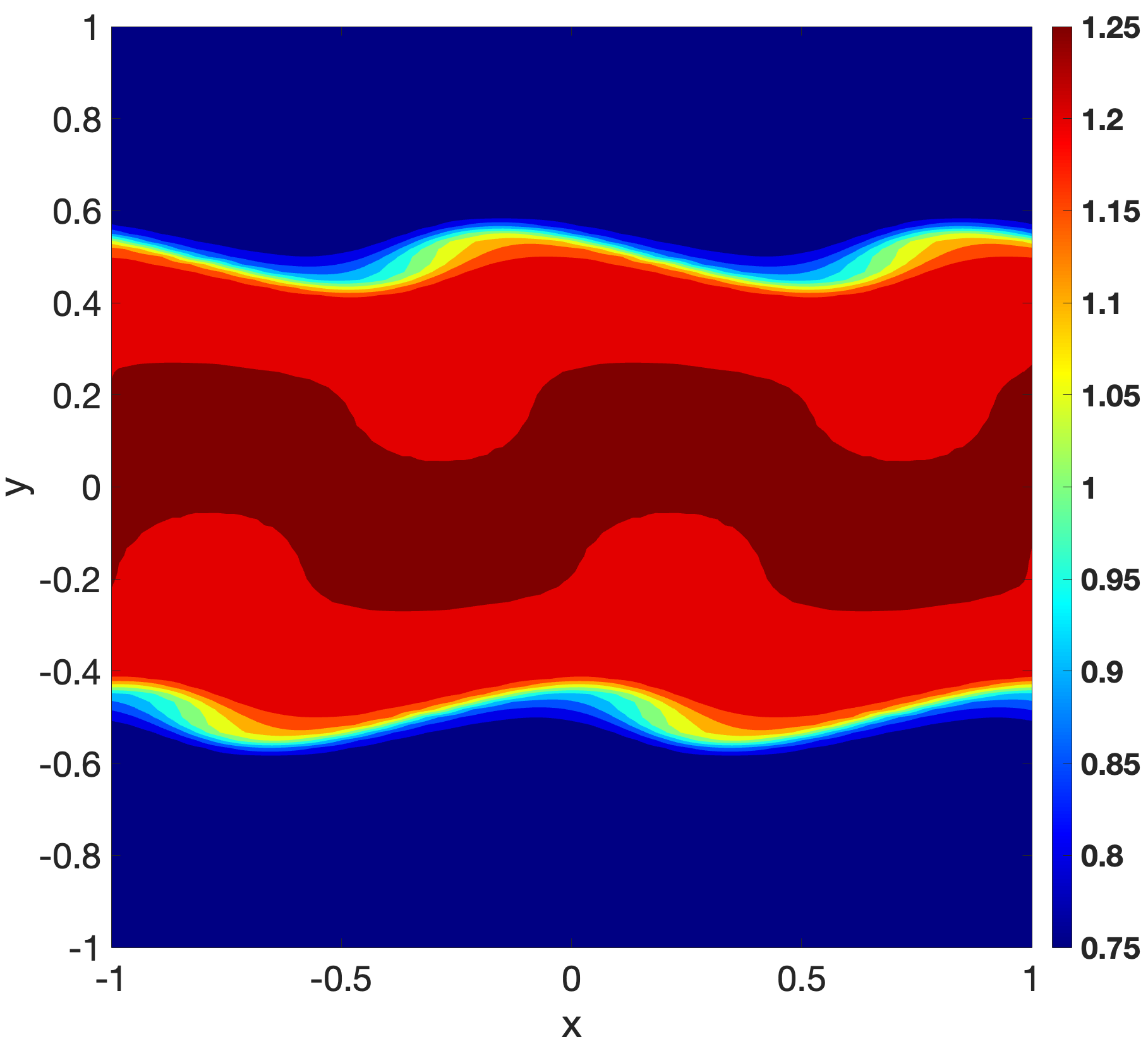}
    \end{subfigure}
    \begin{subfigure}{0.475\textwidth}
	\centering
        \includegraphics[width = 0.8\textwidth]{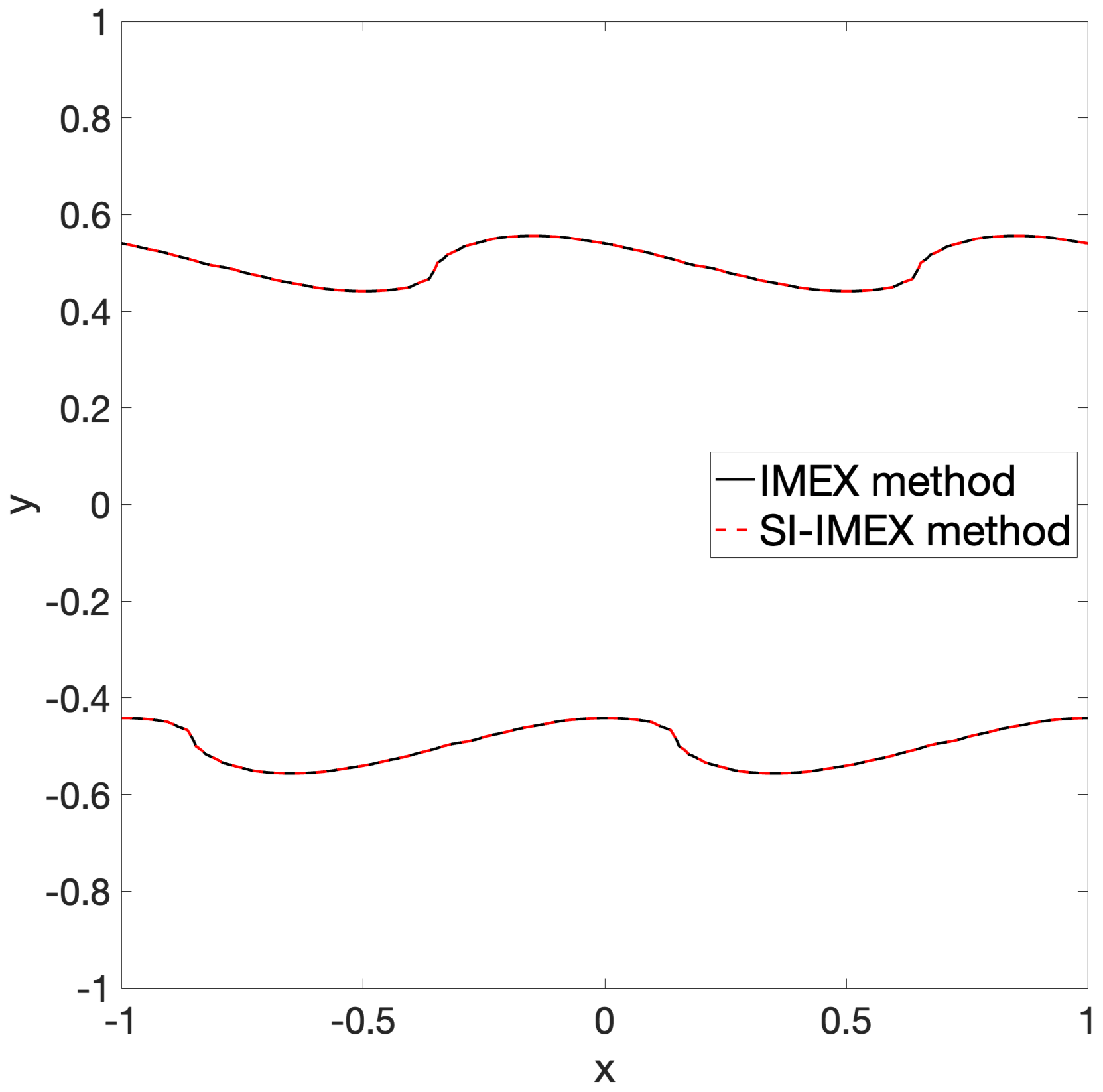}
    \end{subfigure}
    \begin{subfigure}{0.475\textwidth}
	\centering
        \includegraphics[width = 0.85\textwidth]{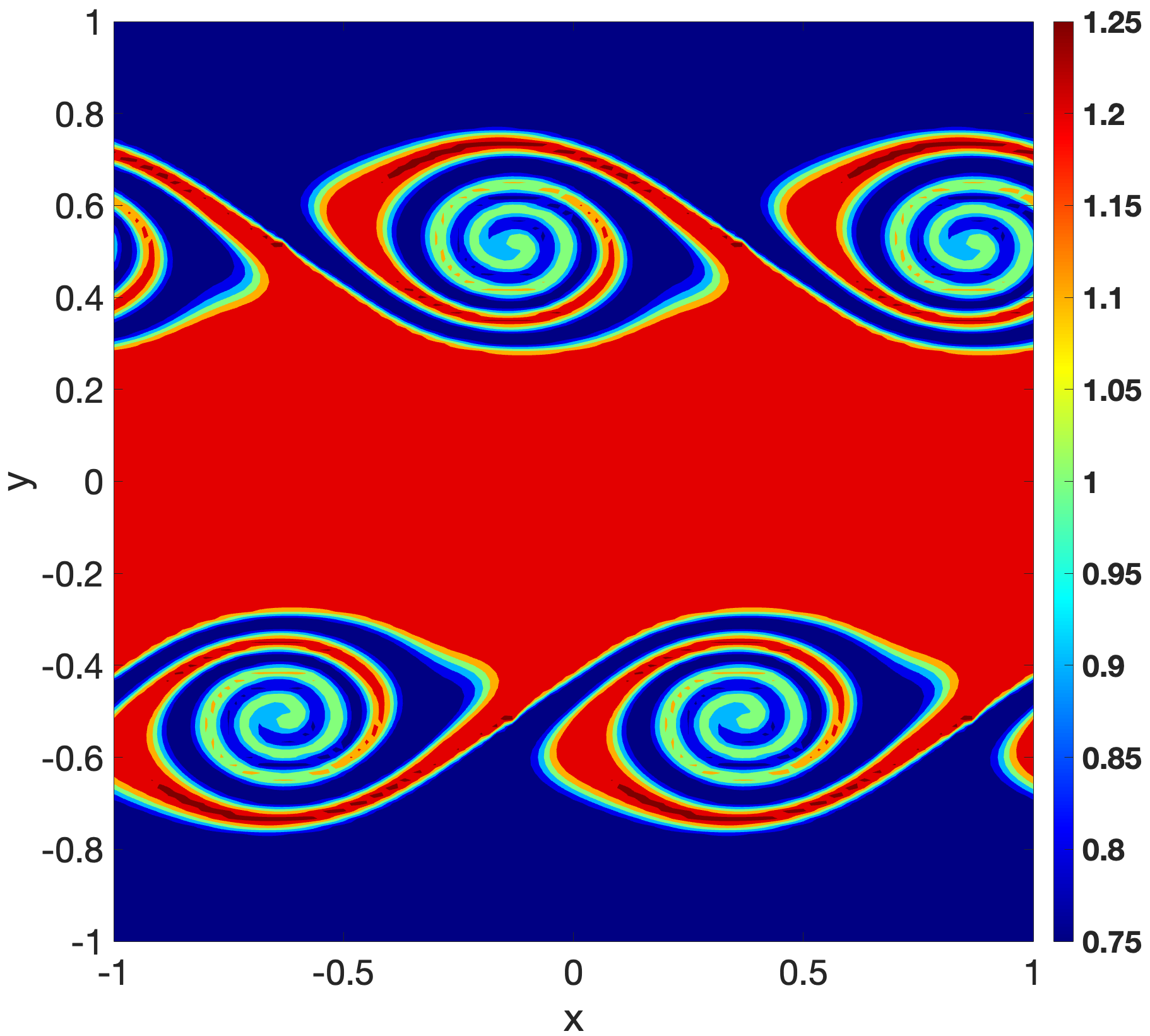}
    \end{subfigure}
    \begin{subfigure}{0.475\textwidth}
	\centering
        \includegraphics[width = 0.8\textwidth]{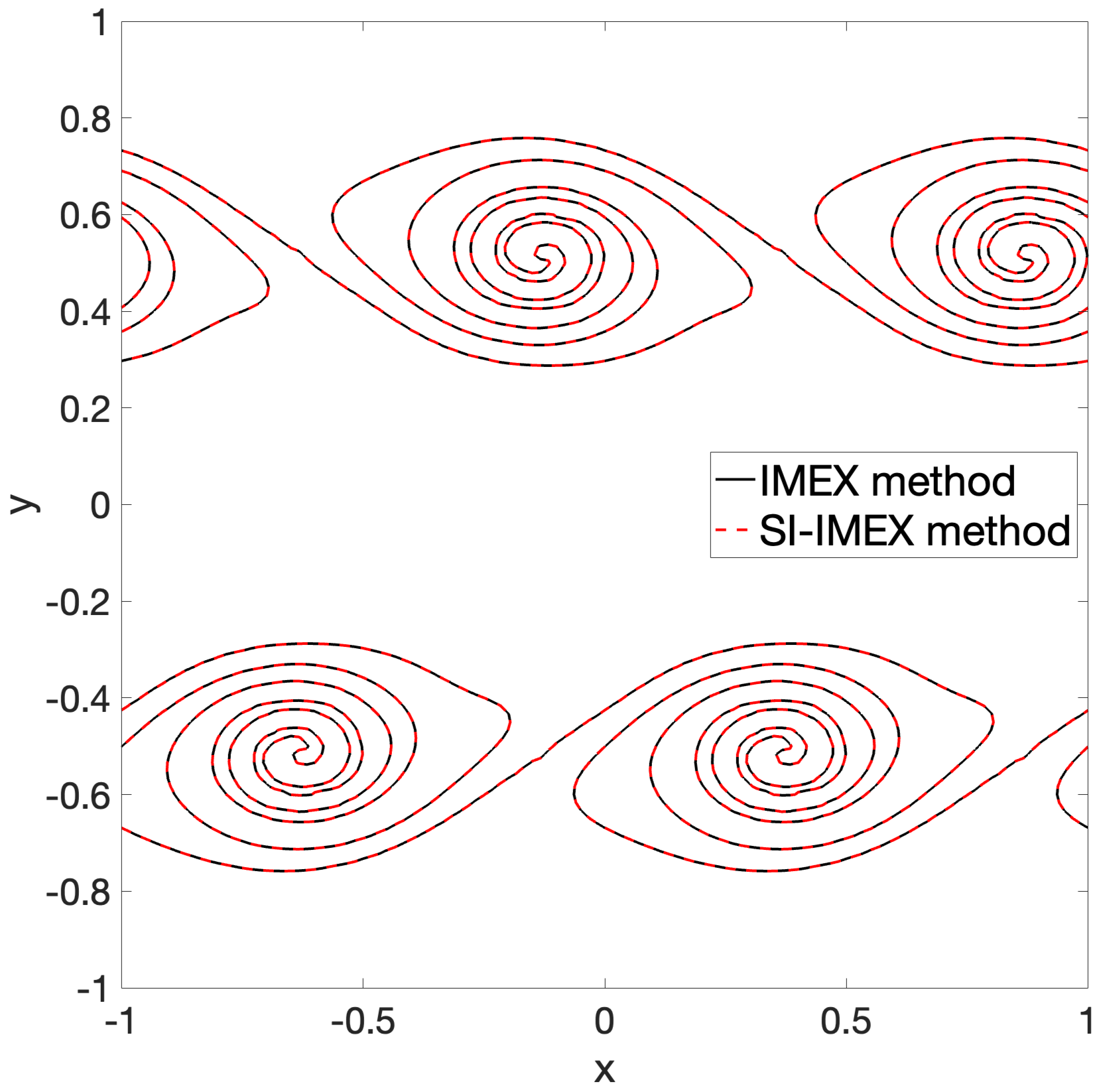}
    \end{subfigure}
    \caption{Kelvin-Helmholtz instability with the ideal gas law \eqref{eq:ideal_gas}, results using the second order scheme (Table \ref{tab:rk2_butch}) and polynomial degree $r = 1$ with $N_{el} = 120$. Top: results at $t = 2$. Bottom: results at $t = T_{f} = 5$. Left: contour plots of the density field obtained with the IMEX method. Right: comparison between the IMEX method and the SI-IMEX method for the isoline equal to $1$. The continuous black lines show the results with IMEX method, while the dashed red lines represent the results with the SI-IMEX.}
    \label{fig:KH_IG_second_order} 
\end{figure}

In a final test, we consider the Kelvin-Helmholtz instability at low Mach studied, e.g., in \cite{zampa:2025}, which we briefly recall here for the convenience of the reader. The computational domain is the square $\Omega = \left(-1,1\right)^{2}$ endowed with periodic boundary conditions, while the final time is $T_{f} = 5$. The initial conditions are
\begin{subequations}
    \begin{eqnarray}
        \rho\left(\mathbf{x},0\right) &=& 1 - \frac{1}{4}\tanh\left[25\left(\left|y\right| - \frac{1}{2}\right)\right] \\
	p\left(\mathbf{x},0\right) &=& \frac{10^{4}}{\gamma} \\
        u\left(\mathbf{x},0\right) &=& -\frac{1}{2}\tanh\left[25\left(\left|y\right| - \frac{1}{2}\right)\right] \\
        v\left(\mathbf{x},0\right) &=& \frac{1}{100}\sin\left(2\pi x\right)\cos\left(2\pi y\right),
    \end{eqnarray}
\end{subequations}	
with $\gamma = 1.4$. Notice that, as discussed in \cite{zampa:2025}, this configuration is in the incompressible regime, but the density is not constant. This shows that assuming a constant background density in \eqref{eq:rho_expansion}, as done in some contributions, is not always a valid assumption, even in the incompressible regime. First, we employ the ideal gas law \eqref{eq:ideal_gas}. We consider the IMEX-RK(2,2,2) scheme (Table \ref{tab:rk2_butch}) with polynomial degree $r = 1$ for the space discretization. The computational mesh is composed by $N_{el} = 120 \times 120 = 14400$ elements along each direction, whereas the time step is $\Delta t = 2 \times 10^{-3}$, yielding a maximum acoustic Courant number $C \approx 15$ and a maximum advective Courant number $C_{u} \approx 0.10$. The contours of density at $t = 2$ and at $t = 5$ are in good agreement with the reference results reported in \cite{zampa:2025} (Figure \ref{fig:KH_IG_second_order}). Moreover, since we are analyzing a fluid mechanic instability, every small variation in the flow can lead to large variations \cite{orlando:2023b, orlando:2024a}, and the excellent agreement obtained between the IMEX method and the SI-IMEX further confirms the correctness of our implementation and the properties of both methods.

\begin{figure}[h!]
    \centering
    \begin{subfigure}{0.475\textwidth}
	\centering
        \includegraphics[width = 0.85\textwidth]{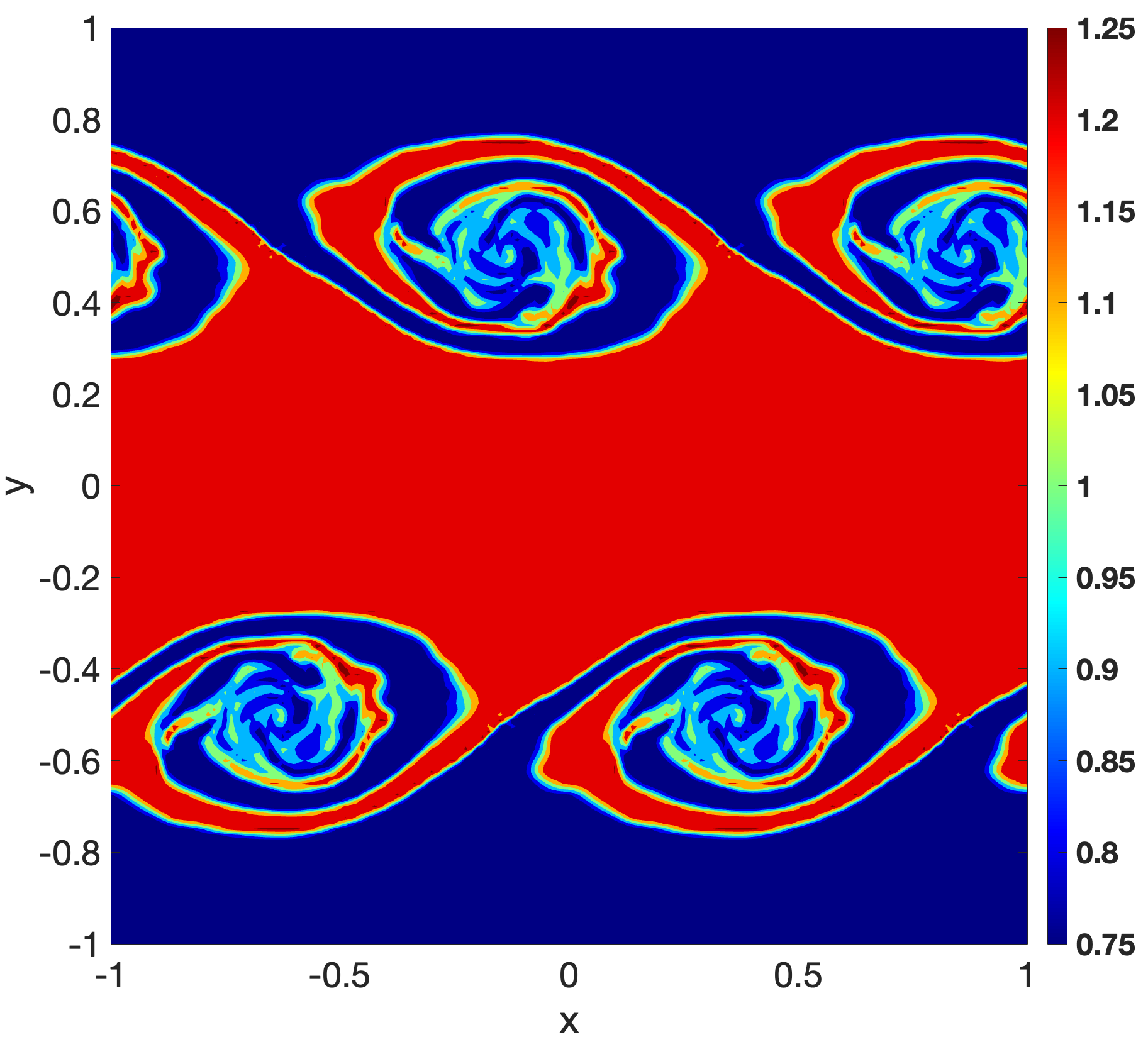}
    \end{subfigure}
    \begin{subfigure}{0.475\textwidth}
	\centering
        \includegraphics[width = 0.8\textwidth]{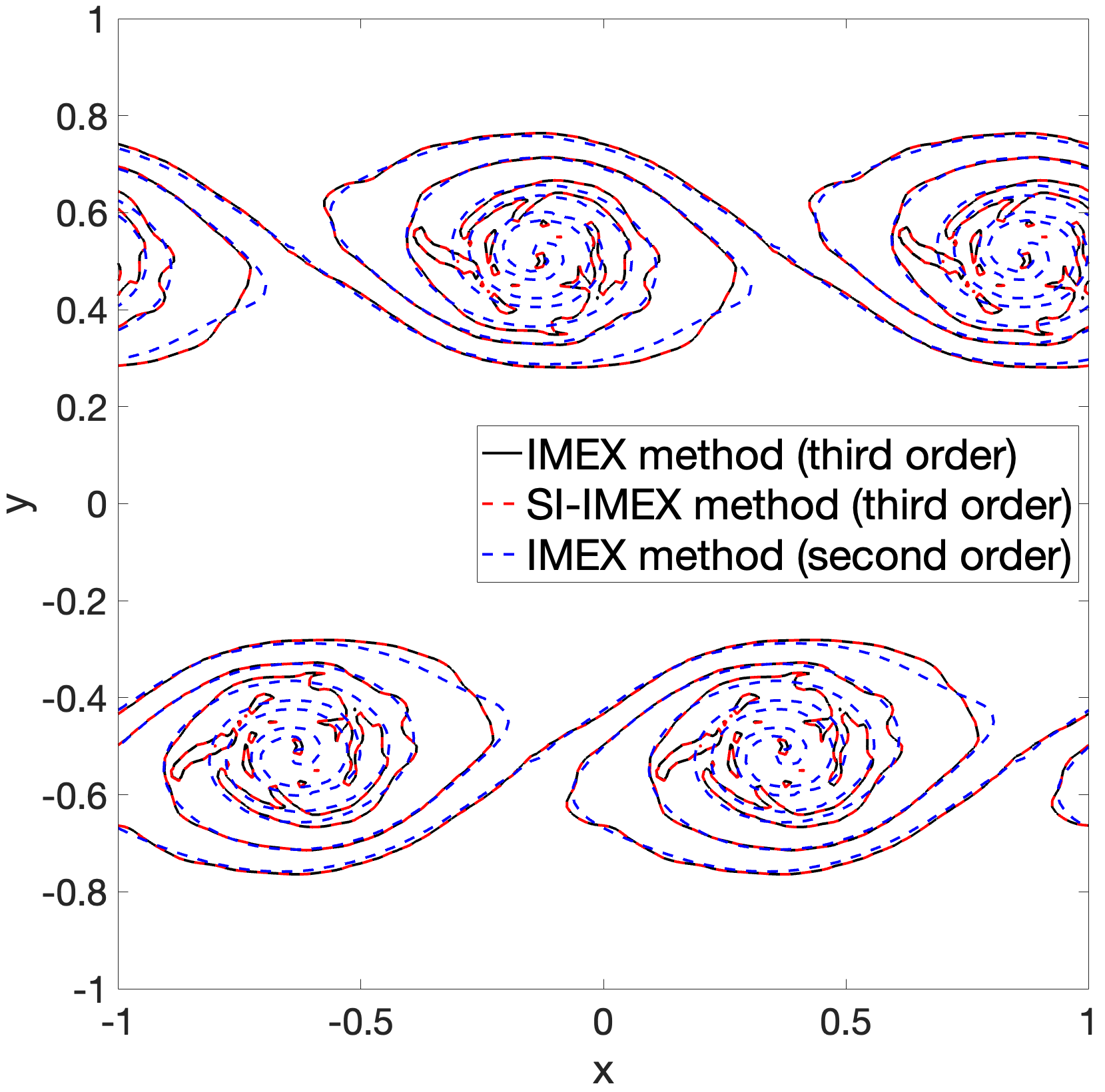}
    \end{subfigure}
    \caption{Kelvin-Helmholtz instability with the ideal gas law \eqref{eq:ideal_gas}, results using the third order scheme of type II (Table \ref{tab:rk3_butch_type_II}) and polynomial degree $r = 2$ with $N_{el} = 80$ at $t = T_{f} = 5$. Left: contour plots of the density field obtained with the IMEX method. Right: comparison between the IMEX method, the SI-IMEX method, and the second order IMEX method for the isoline equal to $1$. The continuous black line shows the results with third order IMEX method, the dashed red line represents the results with the third order SI-IMEX method, whereas the dashed-dotted blue line reports the results with the second order IMEX method.}
    \label{fig:KH_IG_comparison} 
\end{figure}

Next, we employ the third order scheme of type II reported (Table \ref{tab:rk3_butch_type_II}) using polynomial degree $r = 2$ and $N_{el} = 80 \times 80 = 6400$ elements, so that the number of degrees of freedom does not change. Recall indeed that the total number of degrees of freedom per scalar variable is $N_{el}\left(r + 1\right)^{2}$. One can easily notice that the use of higher order methods is beneficial to resolve the fine details of the solution (Figure \ref{fig:KH_IG_comparison}). An excellent agreement is once more established between the IMEX method and the SI-IMEX method (Figure \ref{fig:KH_IG_comparison}).

\begin{figure}[h!]
    \centering
    \begin{subfigure}{0.475\textwidth}
	\centering
        \includegraphics[width = 0.85\textwidth]{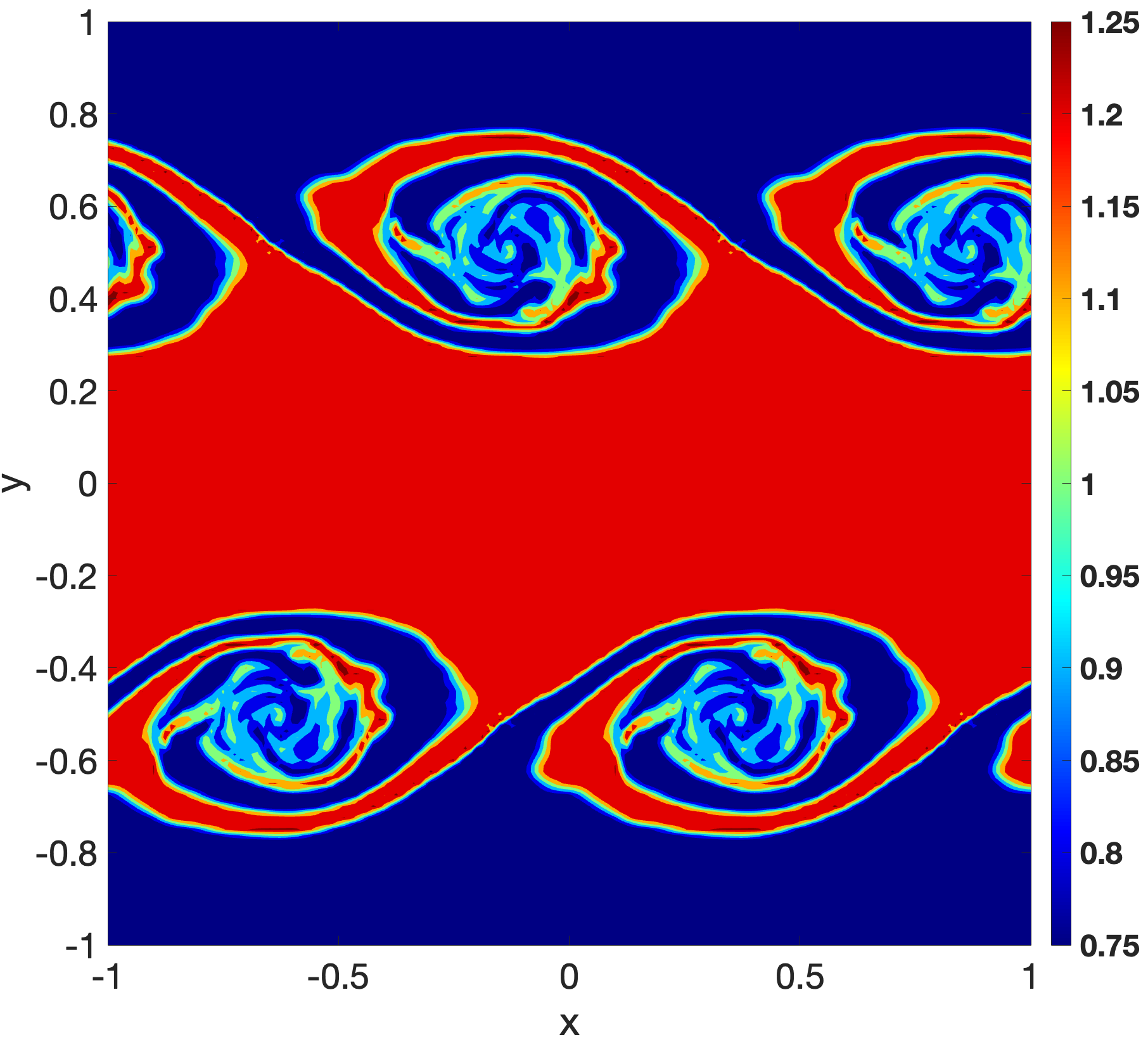}
    \end{subfigure}
    \begin{subfigure}{0.475\textwidth}
	\centering
        \includegraphics[width = 0.8\textwidth]{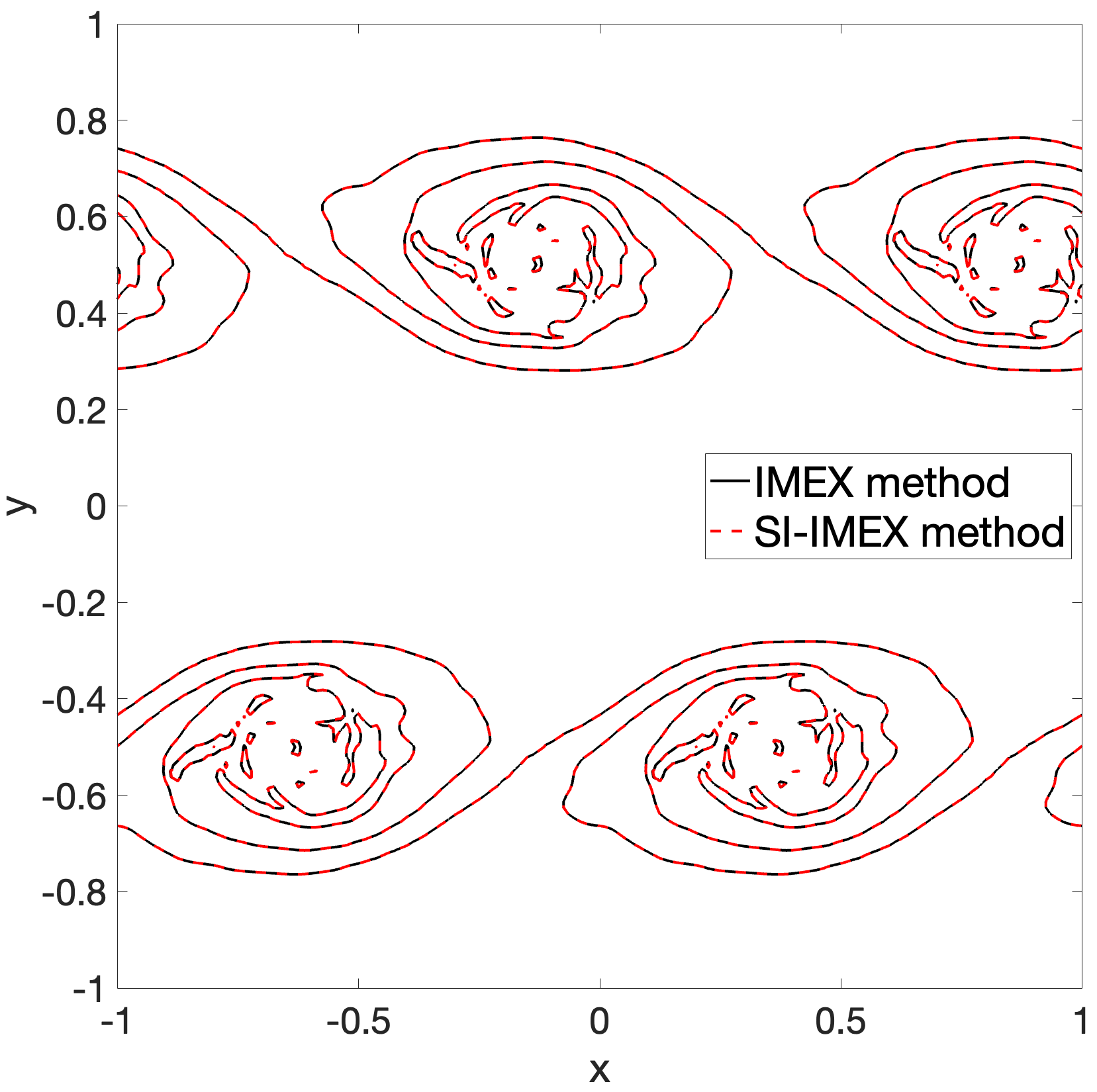}
    \end{subfigure}
    \begin{subfigure}{0.475\textwidth}
	\centering
        \includegraphics[width = 0.8\textwidth]{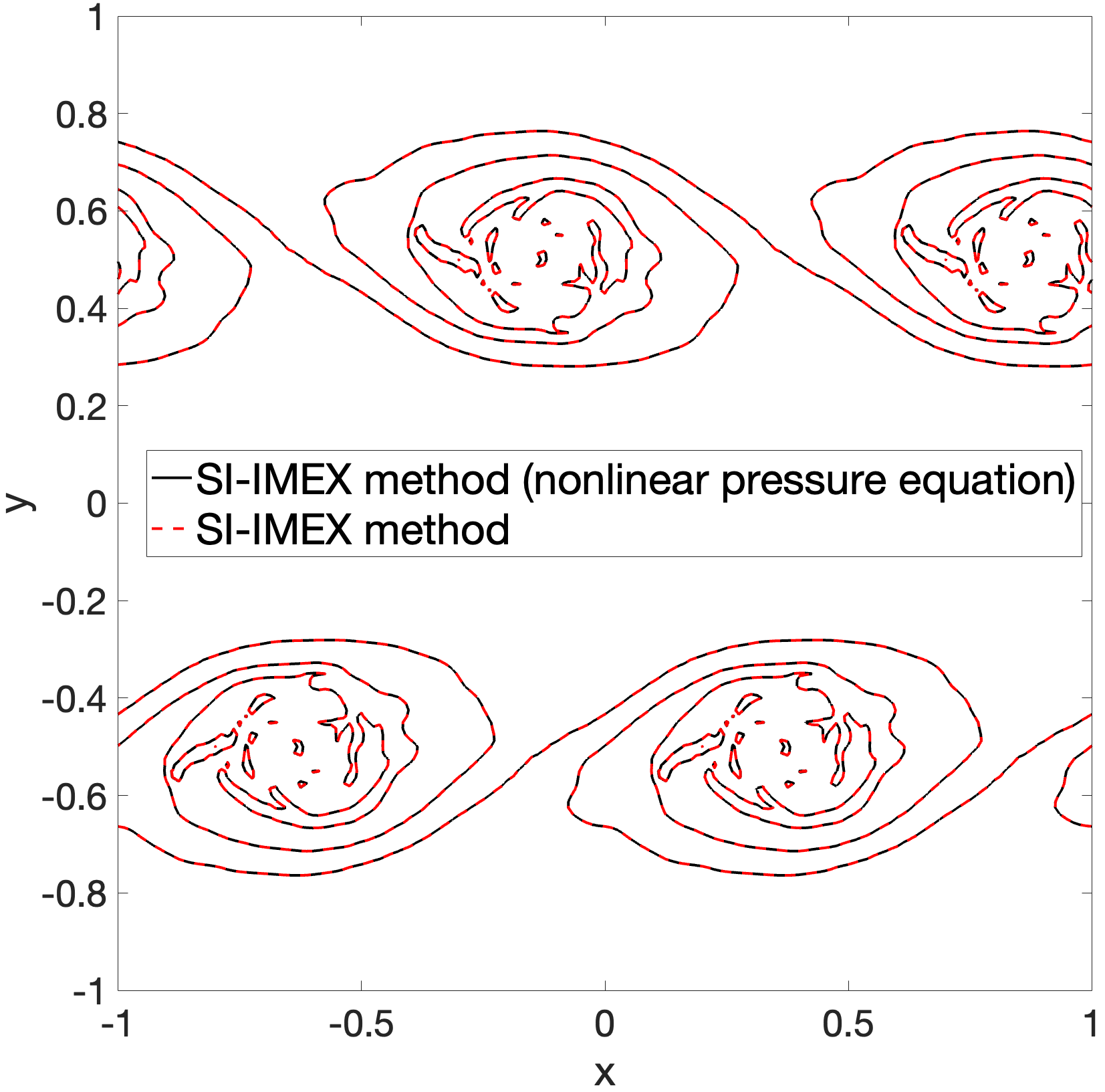}
    \end{subfigure}
    \caption{Kelvin-Helmholtz instability with the Peng-Robinson EOS \eqref{eq:ideal_gas}, results using the third order scheme of type II (Table \ref{tab:rk3_butch_type_II}) and polynomial degree $r = 2$ with $N_{el} = 80$. Top-left: contour plots of the density field obtained with the IMEX method at $t = T_{f} = 5$. Top-right: comparison between the IMEX method and the SI-IMEX method for the isoline equal to $1$. The continuous black lines show the results with IMEX method, while the dashed red lines represent the results with the SI-IMEX method. Bottom: comparison between the SI-IMEX method solving the nonlinear pressure equation and the SI-IMEX method using the linearization proposed in Section \ref{ssec:num_eos} for the isoline equal to $1$. The continuous black lines show the results with SI-IMEX method, while the dashed red lines represent the results with the linearized SI-IMEX.}
    \label{fig:KH_PR} 
\end{figure}

Next, we consider an extension of this test case employing the general cubic EOS \eqref{eq:int_energy_general_cubic_eos}. Following \cite{boscheri:2021a}, we set 
\begin{equation}
    a(T) = \frac{1}{2\sqrt{T}}
\end{equation}
for the attraction term in \eqref{eq:int_energy_general_cubic_eos}. Moreover, we set $b = 5 \times 10^{-3}$. Finally, we assume that \eqref{eq:int_energy_general_cubic_eos_constant} is valid, with $c_{v} \approx 742.0$, and $R_{g} \approx 296.8$ in \eqref{eq:general_cubic_eos}, which correspond to the thermodynamic properties of the nitrogen. An excellent agreement is obtained between the IMEX method and the SI-IMEX method, in spite of the linearization described in \eqref{eq:pressure_elliptic_SI_IMEX_mod} (Figure \ref{fig:KH_PR}). For the sake of completeness, we have also computed the solution obtained with the SI-IMEX method without linearizing the relation between internal energy and pressure, solving therefore a nonlinear equation for the pressure. No difference arises in the development of the Kelvin-Helmholtz instability (Figure \ref{fig:KH_PR}) and a computational time saving of around $70\%$ is obtained thanks to the linearization proposed in \eqref{eq:pressure_elliptic_SI_IMEX_mod}.

Finally, we consider two different configurations. First, we modify the initial condition, taking
\begin{equation}
    p\left(\mathbf{x}, 0\right) = \frac{10^{5}}{\gamma},
\end{equation}
so that a more realistic maximum temperature $T \approx 320$ is obtained. In this case, we need to decrease the time step of the SI-IMEX method to $\Delta t = 5 \times 10^{-4}$ in order to achieve a stable solution. Next, we consider the proper expression of the attraction term $a(T)$ and of the co-volume $b$ for the Peng-Robinson EOS \cite{fernandez:2009}, \cite[p.~263]{sandler:2017}:
\begin{equation}
    \begin{cases}
	a(T) &= 0.45724\frac{R_{g}^{2}T_{c}^{2}}{p_{c}}\alpha(T)^{2} \\
	\alpha(T) &= 1 + \Gamma\left(1 - \sqrt{\frac{T}{T_{c}}}\right) \\
	\Gamma &= 0.37464 + 1.54226\omega - 0.26992\omega^{2} \\
	b &= 0.07780\frac{R_{g}T_{c}}{p_{c}},
    \end{cases}
\end{equation}
where $T_{c}$ denotes the critical temperature, $p_{c}$ the critical pressure, and $\omega$ the acentric factor. For what concerns the nitrogen, we find $T_{c} = \SI{126.19}{\kelvin}, p_{c} = \SI[parse-numbers=false]{3.3978 \cdot 10^{6}}{\pascal}$, and $\omega = 0.0372$ \cite{lias:2010, jacobsen:1986}. Moreover, we consider the following relation for $e^{\#}(T)$ in \eqref{eq:int_energy_general_cubic_eos} \cite{chase:1996, lias:2010}:
\begin{eqnarray}\label{eq:cv_N2}
    e^{\#}\left(T\right) &=& \left[A\frac{T}{1000} + \frac{1}{2}B \left(\frac{T}{1000}\right)^2 + \frac{1}{3}C\left(\frac{T}{1000}\right)^3 + \frac{1}{4}D\left(\frac{T}{1000}\right)^4 - E\frac{1000}{T}\right]\frac{10^6}{M_{w}} - R_{g}T,
\end{eqnarray}
with $M_{w} = \SI{28.0134}{\gram\per\mole}$ and $A,B,C,D,E$ denoting suitable coefficients whose values are reported in Table \ref{tab:cv_N2}. Notice that the polynomial expansion employed in \cite{lias:2010} provides results expressed in \SI[parse-numbers=false]{}{\kilo\joule\per\mole}. Hence, a proper conversion to obtain results in \SI[parse-numbers=false]{}{\kilo\joule\per\kilo\gram} has to applied. The same consideration holds for the factor $1000$, since the argument of the polynomial is expressed in thousandths of Kelvin. The IMEX method can achieve a stable solution for this configuration, whereas severe issues arise for the SI-IMEX method in the development of the instability.

\begin{table}[h!]
    \begin{center}
	   \begin{tabular}{|c|c|}
		\hline
		\texttt{A} & $28.98641$   \\
		\hline
		$B$ & $1.853978$ \\ 
		\hline
		$C$ & $-9.647459$ \\
		\hline
		$D$ & $16.63537$ \\
		\hline
		$E$ & $0.000117$ \\
		\hline
	\end{tabular}
    \end{center}
    \caption{Values of the coefficients in \eqref{eq:cv_N2}}
    \label{tab:cv_N2}
\end{table}

\section{Conclusions}
\label{sec:conclu}

Based on the experience of \cite{boscarino:2022} and \cite{orlando:2025}, we have performed a quantitative comparison between two different Implicit-Explicit Runge-Kutta (IMEX-RK) approaches for the Euler equations of gas dynamics. The two methods are particularly well suited for low Mach number flows, but keep their full accuracy for moderate values of the Mach number. The spatial discretization is based on the Discontinuous Galerkin (DG) method, which naturally allows for high-order accuracy, even though it is characterized by some limitations in the very low Mach limit on quadrilateral cells. The two schemes, namely the IMEX-DG method and the SI-IMEX-DG method, have been compared in a number of relevant benchmarks for ideal gases and on their non-trivial extension for non-ideal gases. The stiff dependence has been carefully analyzed in order to avoid the solution of a nonlinear pressure equation for a general class of equations of state (EOS).

First, we have assessed the convergence properties of the two methods. We have shown that they are asymptotic-preserving (AP) and asymptotically-accurate (AA) in the range of low Mach accuracy established by the spatial discretization. The SI-IMEX method provides a sizeable computational time saving ensuring the same level of accuracy, in particular for moderate values of the Mach number. Moreover, we have noticed an impact of the time discretization strategy and also of the specific time discretization scheme in the activation of spurious modes for low Mach numbers. More specifically, schemes of type I provide a superior stability for low values of the Mach number. A detailed analysis of the spatial discretization will be matter of future work.

Next, we have considered non-trivial and less standard configurations. First, we have analyzed the case in which a time-dependent pressure is imposed at the boundary, for which the asymptotic limit does not coincide with the incompressible Euler equations. Notice that this configuration is outside of the theoretical framework of the SI-IMEX method and, indeed, only some IMEX-RK schemes of type II are suitable. The SI-IMEX method requires a significant smaller time step with respect to that needed by the IMEX method to achieve a stable solution. 

Finally, we have focused on a Kelvin-Helmholtz instability at low Mach, that is in the incompressible regime, but for which the density is not constant. Here we have tested for non-ideal gases a novel linearization in the relation between internal energy and pressure, so as to avoid the solution of a nonlinear equation for the pressure. The proposed linearization is applicable to a general EOS and automatically recovers the linear systems obtained using the ideal gas law. No evident loss of accuracy occurs and a significant computational time saving is established.

In future work, as already mentioned, we aim to employ a spatial discretization based on simplices or Voronoi meshes that has been shown to be low Mach accurate for steady flows. Moreover, we aim to employ a spatial discretization based on compatible finite elements so as to improve the scaling properties with respect to the Mach number $M$. Finally, we aim to further analyze the stability properties of the two methods and to consider an extension of these approaches for the compressible Navier-Stokes equations and for two-phase flows. 

\section*{Acknowledgements}

We thank the two anonymous reviewers for their very useful and constructive comments and remarks, which have greatly helped in improving the quality of the presentation of our results. G. Orlando would like to acknowledge Vincent Perrier for useful discussions on related topics. G. Orlando, S. Boscarino, and G. Russo are part of the INdAM-GNCS National Research Group. 

The simulations have been partly run at CINECA thanks to the computational resources made available through the ISCRA-C projects FEM-GPU - HP10CQYKJ1 and DGNWP - HP10C121HQ. We acknowledge the CINECA award, for the availability of high-performance computing resources and support. 

G. Russo and S. Boscarino thank the Italian Ministry of Instruction, University and Research (MIUR) to support this research with funds coming from PRIN Project 2022, 2022KA3JBA, entitled ``Advanced numerical methods for time dependent parametric partial differential equations and applications'' and from PRIN 2022 PNRR ``FIN4GEO: Forward and Inverse Numerical Modeling of hydrothermal systems in volcanic regions with application to geothermal energy exploitation.'', (No. P2022BNB97). Both authors also have been supported for this work by the Spoke 1“FutureHPC\&BigData” of the Italian Research Center on High-Performance Computing, Big Data and Quantum Computing (ICSC) funded by MUR Missione 4 Componente 2 Investimento 1.4: Potenziamento strutture di ricerca e creazione di ``campioni nazionali di R\&S (M4C2-19 )'' - Next Generation EU (NGEU).

\appendix

\section{Coefficients of employed IMEX-RK schemes}
\label{app:IMEX_coeffs}

We report here for the convenience of the reader some information concerning the IMEX-RK schemes employed in the numerical simulations. We consider the second order IMEX-RK scheme proposed in \cite{giraldo:2013} and also employed in \cite{orlando:2022, orlando:2025}, whose coefficients are reported in the following Butcher tableaux

\begin{table}[pos=H]
    \begin{minipage}{0.45\textwidth}
	\begin{center}
		\begin{tabular}{c|ccc}
			$0$ & $0$ & $0$ & $0$  \\ [0.05cm]
			$2 - \sqrt{2}$ & $2 - \sqrt{2}$ & 0 & 0 \\ [0.05cm]
			$1$ & $\frac{1}{2}$ & $\frac{1}{2}$ & 0 \\
			\hline \\ [-0.35cm]
                & $\frac{\sqrt{2}}{4}$ & $\frac{\sqrt{2}}{4}$ & $1 - \frac{\sqrt{2}}{2}$
		  \end{tabular}
	\end{center}
    \end{minipage}
    \begin{minipage}{0.45\textwidth}
	\begin{center}
		\begin{tabular}{c|ccc}
			$0$ & $0$ & $0$ & $0$ \\ [0.05cm]
                $2 - \sqrt{2}$ & $1 - \frac{\sqrt{2}}{2}$ & $1 - \frac{\sqrt{2}}{2}$ & $0$ \\ [0.05cm]
                $1$ & $\frac{\sqrt{2}}{4}$ & $\frac{\sqrt{2}}{4}$ & $1 - \frac{\sqrt{2}}{2}$ \\
			\hline \\ [-0.35cm]
                & $\frac{\sqrt{2}}{4}$ & $\frac{\sqrt{2}}{4}$ & $1 - \frac{\sqrt{2}}{2}$
		\end{tabular}
	\end{center}
    \end{minipage}
    \caption{\it Butcher tableaux of the IMEX-RK(2,2,2) scheme in \cite{giraldo:2013}. Left: explicit method. Right: implicit method.}
    \label{tab:rk2_butch}
\end{table}
\noindent
For what concerns third order time discretization methods, we consider the following method of type I \cite{boscarino:2016}:

\begin{table}[pos=H]
    \begin{minipage}{0.9\textwidth}
	\begin{center}
		\begin{tabular}{c|cccc}
			$0$ & $0$ & $0$ & $0$ & $0$ \\
			$0.435866521508$ & $0.435866521508$ & $0$ & $0$ & $0$ \\
                $0.717933260754$ & $0.435866521508$ & $0.282066739245$ & $0$ & $0$ \\
                $1$ & $-0.733534082748750$ & $2.150527381100$ & $-0.416993298352$ & $0$ \\
                \hline \\ [-0.35cm]
                & $0$ & $1.208496649176$ & $-0.644363170684$ & $0.435866521508$
		\end{tabular}
	\end{center}
    \end{minipage}
    \vskip 0.3cm
    \begin{minipage}{0.9\textwidth}
	\begin{center}
		\begin{tabular}{c|cccc}
                $0.435866521508$ & $0.435866521508$ & $0$ & $0$ & $0$ \\
                $0.435866521508$ & $0$ & $0.435866521508$ & $0$ & $0$ \\
                $0.717933260754$ & $0$ & $0.282066739245$ & $0.435866521508$ & $0$ \\
                $1$ & $0$ & $1.208496649176$ & $-0.644363170684$ & $0.435866521508$ \\
			\hline \\ [-0.35cm]
			& $0$ & $1.208496649176$ & $-0.644363170684$ & $0.435866521508$
		\end{tabular}
	\end{center}
    \end{minipage}
    \caption{\it Butcher tableaux of the IMEX-RK(4,3,3) scheme in \cite{boscarino:2022}. Top: explicit method. Bottom: implicit method.}
    \label{tab:rk3_butch_type_I}
\end{table}
\noindent
and the following scheme of type II \cite{kennedy:2003}:

\begin{table}[pos=H]
    \begin{minipage}{0.9\textwidth}
	\begin{center}
		\begin{tabular}{c|cccc}
			$0$ & $0$ & $0$ & $0$ & $0$ \\ [2mm]
                $\frac{1767732205903}{2027836641118}$ & $\frac{1767732205903}{2027836641118}$ & $0$ & $0$ & $0$ \\ [2mm]
                $\frac{3}{5}$ & $\frac{5535828885825}{10492691773637}$ & $\frac{788022342437}{10882634858940}$ & $0$ & $0$ \\ [2mm]
                $1$ & $\frac{6485989280629}{16251701735622}$ & $\frac{-4246266847089}{9704473918619}$ & $ \frac{10755448449292}{10357097424841}$ & $0$ \\ [2mm]
                \hline \\ [-0.35cm]
                & $\frac{1471266399579}{7840856788654}$ & 
                $\frac{-4482444167858}{7529755066697}$ & $\frac{11266239266428}{11593286722821}$ & $\frac{1767732205903}{4055673282236}$
		\end{tabular}
	\end{center}
    \end{minipage}
    \vskip 0.3cm
    \begin{minipage}{0.9\textwidth}
	\begin{center}
		\begin{tabular}{c|cccc}
			$0$ & $0$ & $0$ & $0$ & $0$ \\ [2mm]
                $\frac{1767732205903}{2027836641118}$ & $\frac{1767732205903}{4055673282236}$ & $\frac{1767732205903}{4055673282236}$ & $0$ & $0$ \\ [2mm]
				$\frac{3}{5}$ & $\frac{2746238789719}{10658868560708}$ &
                $-\frac{640167445237}{6845629431997}$ & $\frac{1767732205903}{4055673282236}$ & $0$ \\ [2mm]
                $1$ & $\frac{1471266399579}{7840856788654}$ & 
                $\frac{-4482444167858}{7529755066697}$ & $\frac{11266239266428}{11593286722821}$ & $\frac{1767732205903}{4055673282236}$ \\[2mm]
			\hline \\ [-0.35cm]
			& $\frac{1471266399579}{7840856788654}$ & 
                $\frac{-4482444167858}{7529755066697}$ & $\frac{11266239266428}{11593286722821}$ & $\frac{1767732205903}{4055673282236}$
		\end{tabular}
	\end{center}
    \end{minipage}
    \caption{\it Butcher tableaux of the IMEX-RK(3,3,3) scheme. Top: explicit method. Bottom: implicit method.}
    \label{tab:rk3_butch_type_II}
\end{table}
\noindent
Finally, we employ the fourth order time discretization method of type ARS proposed in \cite{calvo:2001}

\begin{table}[pos=H]
    \begin{minipage}{0.9\textwidth}
	\begin{center}
		\begin{tabular}{c|rrrrrr}
			$0$ & $0$ & $0$ & $0$ & $0$ & $0$ & $0$ \\[1.5mm]
                $\frac{1}{4}$ & $\frac{1}{4}$ & $0$ & $0$ & $0$ & $0$ & $0$ \\[1.5mm]
                $\frac{3}{4}$ & $-\frac{1}{4}$ & $1$ & $0$ & $0$ & $0$ & $0$ \\[1.5mm]
                $\frac{11}{20}$ & $-\frac{13}{100}$ & $\frac{43}{75}$ & $\frac{8}{75}$ & $0$ & $0$ & $0$ \\[1.5mm]
                $\frac{1}{2}$ & $-\frac{6}{85}$ & $\frac{42}{85}$ & $\frac{179}{1360}$ & $-\frac{15}{272}$ & $0$ & $0$ \\[1.5mm]
                $1$ & $0$ & $\frac{79}{24}$ & $-\frac{5}{8}$ & $\frac{25}{2}$ & $-\frac{85}{6}$ & $0$ \\[1.5mm]
                \hline \\[-0.3cm]
                & $0$ &	$\frac{25}{24}$ & $-\frac{49}{48}$ & $\frac{125}{16}$ & $-\frac{85}{12}$ & $\frac{1}{4}$
		  \end{tabular}
	\end{center}
    \end{minipage}
    \vskip 0.3cm
    \begin{minipage}{0.9\textwidth}
	\begin{center}
		\begin{tabular}{c|rrrrrr}
                $0$ & $0$ & $0$ & $0$ & $0$ & $0$ & $0$\\[1.5mm]
                $\frac{1}{4}$ & $0$ & $\frac{1}{4}$ & $0$ & $0$ & $0$ & $0$ \\[1.5mm]
                $\frac{3}{4}$ & $0$ & $\frac{1}{2}$ & $\frac{1}{4}$ & $0$ & $0$ & $0$ \\[1.5mm]
                $\frac{11}{20}$ & $0$ & $\frac{17}{50}$ & $-\frac{1}{25}$ & $\frac{1}{4}$ & $0$ & $0$ \\[1.5mm]
                $\frac{1}{2}$ & $0$ & $\frac{371}{1360}$ & $-\frac{137}{2720}$ & $\frac{15}{544}$ & $\frac{1}{4}$ & $0$ \\[1.5mm]
                $1$ & $0$ &	$\frac{25}{24}$ & $-\frac{49}{48}$ & $\frac{125}{16}$ & $-\frac{85}{12}$ & $\frac{1}{4}$ \\[1.5mm]
			\hline \\[-0.3cm]
                & $0$ &	$\frac{25}{24}$ & $-\frac{49}{48}$ & $\frac{125}{16}$ & $-\frac{85}{12}$ & $\frac{1}{4}$
		\end{tabular}
	\end{center}
    \end{minipage}
    \caption{\it Butcher tableaux of the IMEX-RK(5,5,4) scheme. Top: explicit method. Bottom: implicit method.}
    \label{tab:rk4_butch_ARS}
\end{table}
\noindent
and the fourth order time discretization method of type II presented in \cite{kennedy:2019}

\begin{table}[pos=H]
    \begin{minipage}{0.9\textwidth}
	\begin{center}
		\begin{tabular}{c|ccccccc}
			$0$ & $0$ & $0$ & $0$ & $0$ & $0$ & $0$ & $0$ \\
                $2\gamma$ & $2\gamma$ & $0$ & $0$ & $0$ & $0$ & $0$ & $0$ \\
                $c_{3}$ & $\tilde{a}_{31}$ & $\tilde{a}_{32}$ & $0$ & $0$ & $0$ & $0$ & $0$ \\
                $c_{4}$ & $\tilde{a}_{41}$ & $\tilde{a}_{42}$ & $\tilde{a}_{43}$ & $0$ & $0$ & $0$ & $0$ \\
                $c_{5}$ & $\tilde{a}_{51}$ & $\tilde{a}_{52}$ & $\tilde{a}_{53}$ & $\tilde{a}_{54}$ & $0$ & $0$ & $0$ \\
                $c_6$ & $\tilde{a}_{61}$ & $\tilde{a}_{62}$ & $\tilde{a}_{63}$ & $\tilde{a}_{64}$ & $\tilde{a}_{65}$ & $0$ & $0$ \\
                $1$ & $\tilde{a}_{71}$ & $\tilde{a}_{72}$ & $\tilde{a}_{73}$ & $\tilde{a}_{74}$ & $\tilde{a}_{75}$ & $\tilde{a}_{76}$ & $0$ \\
			\hline
			& $0$ &	$0$ & $b_3$ & $b_4$ & $b_5$ & $b_6$ & $\gamma$
		  \end{tabular}
	\end{center}
    \end{minipage}
    \vskip 0.3cm
    \begin{minipage}{0.9\textwidth}
	\begin{center}
		\begin{tabular}{c|ccccccc}
                $0$ & $0$ & $0$ & $0$ & $0$ & $0$ & $0$ & $0$ \\
                $2\gamma$ & $ \gamma$ & $ \gamma$ & $0$ & $0$ & $0$ & $0$ & $0$ \\
                $c_{3}$ & $a_{32}$ & $a_{32}$ & $\gamma$ & $0$ & $0$ & $0$ & $0$ \\
                $c_{4}$ & $a_{42}$ & $a_{42}$ & $a_{43}$ & $ \gamma$ & $0$ & $0$ & $0$ \\
                $c_{5}$ & $a_{52}$ & $a_{52}$ & $a_{53}$ & $a_{54}$ & $ \gamma$ & $0$ & $0$ \\
                $c_{6}$ & $a_{62}$ &	$a_{62}$ & $a_{63}$ & $a_{64}$ & $a_{65}$ & $\gamma$ & $0$ \\
			$1$ & $0$ &	$0$ & $b_3$ & $b_4$ & $b_5$ & $b_6$ & $\gamma$ \\
			\hline
			& $0$ &	$0$ & $b_3$ & $b_4$ & $b_5$ & $b_6$ & $\gamma$
		\end{tabular}
	\end{center}
    \end{minipage}
    \caption{\it Butcher tableaux of the IMEX-RK(6, 6, 4) scheme. Top: explicit method. Bottom: implicit method.}
    \label{tab:rk4_butch_type_II}
\end{table}
with
\begin{align*}
    \begin{array}{ll}
        \gamma = 247/2000, & c_{2} = 2\gamma,\\ c_{3} = (2 + \sqrt{2})\gamma, & c_{4} = 67/200, \\
        c_{5} = 3/40, & c_{6} = 7/10, 
    \end{array}
\end{align*}
\begin{align*}
    \begin{array}{ll}
        b_{3}\phantom{_2} = \phantom{-} 9164257142617/17756377923965, & b_{4}\phantom{_2} = -10812980402763/74029279521829, \\
        b_{5}\phantom{_2} = \phantom{-}1335994250573/5691609445217, & b_{6}\phantom{_2} = \phantom{-}2273837961795/8368240463276, \\
	b_{7}\phantom{_2} = \phantom{-}247/2000, & \\
	a_{32} = \phantom{-}624185399699/4186980696204, &  \\
        a_{42} = \phantom{-}1258591069120/10082082980243, & a_{43} = -322722984531/8455138723562, \\
        a_{52} = -436103496990/5971407786587, & a_{53} = -2689175662187/11046760208243, \\
        a_{54} = \phantom{-}4431412449334/12995360898505, & \\
        a_{62} = -2207373168298/14430576638973, & a_{63} = \phantom{-}242511121179/3358618340039, \\
        a_{64} = \phantom{-}3145666661981/7780404714551, & a_{65} =   \phantom{-}5882073923981/14490790706663.
    \end{array}
\end{align*}
and
\begin{align*}
    \begin{array}{ll}
        \tilde{a}_{31} = 247/4000, & \tilde{a}_{32} = 2694949928731/7487940209513 \\
        \tilde{a}_{41} = 464650059369/8764239774964, & \tilde{a}_{42} = 878889893998/2444806327765, \\
        \tilde{a}_{43} = -952945855348/12294611323341, & \\
        \tilde{a}_{51} = 476636172619/8159180917465, & \tilde{a}_{52} = -1271469283451/7793814740893, \\
        \tilde{a}_{53} = -859560642026/4356155882851, & \tilde{a}_{54} = 1723805262919/4571918432560,\\ 
        \tilde{a}_{61} = 6338158500785/11769362343261, &\tilde{a}_{62} = -4970555480458/10924838743837, \\
        \tilde{a}_{63} = 3326578051521/2647936831840, & \tilde{a}_{64} = -880713585975/1841400956686, \\
        \tilde{a}_{65} = -1428733748635/8843423958496, & \\
        \tilde{a}_{71} = 760814592956/3276306540349, &
        \tilde{a}_{72} = 760814592956/3276306540349, \\
        \tilde{a}_{73} = -47223648122716/6934462133451, & \tilde{a}_{74} = 71187472546993/9669769126921, \\
        \tilde{a}_{75} = -13330509492149/9695768672337, & \tilde{a}_{76} = 11565764226357/8513123442827. 
    \end{array}
\end{align*}

\bibliographystyle{cas-model2-names}
\bibliography{SI_IMEX_Euler.bib}
	
\end{document}